\documentclass[10pt]{article}
\usepackage{hyperref}
\usepackage{amssymb}
\usepackage{amsmath}
\input{liemacs10.sty} 
\usepackage{color}

\renewcommand\up{{\uparrow}}

\newcommand{\sH}{{\sf H}}

\newcommand{\sV}{{\tt V}}
\newcommand{\sW}{{\tt W}}
\newcommand{\sE}{{\tt E}}

\newcommand{\csp}{{\mathfrak{csp}}} 
\newcommand{\hcsp}{{\mathfrak{hcsp}}}

\newcommand{\Inv}{\mathop{{\rm Inv}}\nolimits}

\renewcommand{\phi}{\varphi} 

\newcommand{\AdS}{\mathop{{\rm AdS}}\nolimits}

\newcommand{\Stand}{\mathop{{\rm Stand}}\nolimits}

\renewcommand\mlabel{\label} 
\usepackage{authblk}

\begin{document}
\title{Wedge domains in compactly causal symmetric spaces} 
 
\author{Karl-Hermann Neeb and  Gestur \'Olafsson\thanks{The research of K.-H. Neeb was partially supported
  by DFG-grant NE 413/10-1. The research of G. \'Olafsson was partially supported by Simons grant 586106.}}

\maketitle

\begin{center}
{\bf Dedicated to the memory of Ottmar Loos} 
 
\end{center}

\abstract{
Motivated by construction in Algebraic Quantum Field Theory we
introduce wedge domains in compactly causal symmetric
spaces $M=G/H$, which includes in particular anti de Sitter space 
in all dimensions and 
its coverings. Our wedge domains generalize Rindler wedges in Minkowski space. 
The key geometric structure we use is the modular flow on $M$ 
defined by an Euler element in the Lie algebra of~$G$.
Our main geometric result asserts that three seemingly different  characterizations 
of these domains coincide: the positivity domain of the modular vector 
field; 
the domain specified by a KMS like analytic extension condition for the 
modular flow; and the domain specified by a polar decomposition in terms of 
certain cones. 

In the second half of the article we show that our wedge domains 
share important properties with wedge domains in Minkowski space. 
If $G$ is semisimple, there exist unitary representations $(U,\cH)$ 
of $G$ and isotone covariant nets of real subspaces 
$\sH(\cO) \subeq \cH$, defined for any open subset $\cO \subeq M$, 
which assign to connected components of the 
wedge domains a standard subspace whose modular group 
corresponds to the modular flow on~$M$. This corresponds to the 
Bisognano--Wichmann property in Quantum Field Theory. 
We also show that the set of $G$-translates of the connected components of 
the wedge domain provides a geometric realization of the 
abstract wedge space introduced by the 
first author and V.~Morinelli.} 
\bigskip


\tableofcontents 

\section{Introduction} 
\mlabel{sec:1}

Let $M = G/H$ be a homogeneous space of the Lie group $G$ 
and $(U,\cH)$ a unitary representation of $G$. Motivated by 
their occurrence in Algebraic Quantum Field Theory (AQFT), 
in the sense of Haag--Kastler (\cite{Ha96}), one considers 
nets of real subspaces on $M$, i.e., to each open subset $\cO 
\subeq M$, we assign a closed real subspace $\sH(\cO) \subeq \cH$ 
satisfying 
\begin{itemize}
\item Isotony: $\cO_1 \subeq \cO_2 \Rarrow \sH(\cO_1)  \subeq \sH(\cO_2)$, 
\item Equivariance: $\sH(g\cO) = U(g) \sH(\cO)$ for $g \in G$. 
\end{itemize}
If $(W(v))_{v \in \cH}$ are the Weyl operators on bosonic Fock space 
$\cF_+(\cH)$, then we assign to each closed real subspace 
$\sE \subeq \cH$ the von Neumann algebra $\cR(\sE) := W(\sE)''$ 
(\cite{Si74}). 
We thus obtain on $M$ the isotone covariant net 
$(\cR(\sH(\cO)))_{\cO \subeq M}$ of von Neumann algebras, where 
covariance refers to the canonical representation of $G$ on $\cF_+(\cH)$. 
This method has been developed by Araki and Woods  in the context of 
free bosonic quantum fields (\cite{Ar63, Ar64, AW63, AW68}), 
but it can be modified to deal with fermionic Fock spaces 
(cf.\ \cite{EO73}, \cite{BJL02}) 
and there are variants for other statistics (anyons) 
(\cite{Schr97}, \cite[\S 3]{Le15}). 

We are interested in structures on homogeneous spaces which 
resemble space-times in AQFT. In this context the von Neumann 
algebra associated to an open subset $\cO \subeq M$ corresponds 
to observables measurable in the ``laboratory'' $\cO$ 
(\cite{Ha96}). One often focuses on those 
von Neumann algebras $\cR(\sV)$ for which the ``vacuum vector'' 
$\Omega \in \cF_+(\cH)$ is cyclic and separating, so that 
the Tomita--Takesaki Theorem 
applies to $(\cR(\sV),\Omega)$ (\cite[Thm.~2.5.14]{BR87}). 
It is not hard to see that $\Omega$ is 
\begin{itemize}
\item cyclic for $\cR(\sV)$ if and only if 
$\sV + i \sV$ is dense in $\cH$, 
\item separating for $\cR(\sV)$ if and only if 
$\sV \cap  i \sV = \{0\}$. 
\end{itemize}
A closed real subspace $\sV \subeq \cH$ with both properties is called 
{\it standard} (cf.\ \cite{Lo08} for the basic theory of standard subspaces). 
It is therefore of particular interest to understand how to construct 
nets for which the subspaces $\sH(\cO)$ are standard. 
The main feature of Tomita--Takesaki theory is that it 
provides a modular conjugation and a modular automorphism group, 
and these can already be associated to a standard subspaces $\sV \subeq \cH$ 
(and pass through second quantization to the von Neumann context). 
Concretely, we associate to $\sV$ a {\it pair of modular 
objects} $(\Delta_\sV, J_\sV)$: 
\begin{itemize}
\item the {\it modular operator} $\Delta_\sV$ is a positive 
selfadjoint operator, 
\item $J_\sV$ is a {\it conjugation} (an antiunitary 
involution),
\end{itemize}
and these two operators satisfy the modular relation 
$J_\sV \Delta_\sV J_\sV = \Delta_\sV^{-1}.$ 
The pair $(\Delta_\sV, J_\sV)$ is obtained
by the polar decomposition $\sigma_\sV = J_\sV \Delta_\sV^{1/2}$ of the closed operator 
\[ \sigma_\sV \: \sV + i \sV \to \cH, \quad 
x + i y \mapsto x- iy \] 
with $\sV = \Fix(\sigma_\sV)$. Key properties of these operators are that 
\[ J_\sV \sV = \sV' := \{ w \in \cH \: (\forall v\in \sV)\ \Im \la v,w \ra = 0\}
\quad \mbox{ and } \quad 
\Delta_\sV^{it} \sV = \sV \quad \mbox{ for }  \quad t \in \R.\] 
So we obtain a one-parameter group of automorphisms of $\sV$ 
(the modular group) 
and a symmetry between $\sV$ and its commutant $\sV'$, implemented by $J_\sV$. 

The current interest in standard subspaces arose 
in the 1990s from the work of 
Borchers 
and 
Wiesbrock (\cite{Bo92, Wi93}). 
This led to the concept of modular localization 
in Quantum Field Theory introduced by  
Brunetti, Guido and Longo in \cite{BGL02, BGL93}; see also \cite{BDFS00} 
and \cite{Le15, LL15} for important applications of this technique. 

For a net on a homogeneous space $G/H$, 
the domains $\cO$ for which $\sH(\cO)$ is standard 
are of particular relevance. Here one would like to know when 
the modular group $(\Delta_\cO^{it})_{t \in \R}$ is ``geometric'' 
in the sense that it is implemented by a one-parameter subgroup of~$G$, hence 
corresponds to a one-parameter group of symmetries of~$M$ 
preserving the domain~$\cO$. For the modular 
conjugation $J_\cO$, we may likewise ask for the existence of an involution 
$\tau^M$ on $M$ satisfying $\sH(\tau^M \cO) = \sH(\cO)'$. 
Building on \cite{NO21a}, which deals with the case $H = \{e\}$, 
i.e., left invariant nets on Lie groups, 
we turn in this paper to nets on a class of symmetric spaces to 
which the tools and methods developed in \cite{NO21a}  apply particularly well; 
see also the companion paper \cite{NO21b}. 

First we explain how to find 
natural standard subspaces for unitary representations. 
So let $(U,\cH)$ be a unitary representation of $G$ which extends  
to an {\it antiunitary representation} of the  extended group~
\[ G_\tau  = G \rtimes \{\id_G,\tau^G\},\] 
where $\tau^G$ is an involutive automorphism of $G$. 
Then $J := U(\tau^G)$ is a conjugation satisfying 
$J U(g) J = U(\tau^G(g))$ for $g \in G$. 
We then obtain for each pair 
$(h,\tau^G)$, for which $h$ is fixed by the Lie algebra 
involution $\tau$ induced by $\tau^G$, a standard subspace 
$\sV := \sV_{(h,\tau^G,U)} = \Fix(J_\sV \Delta_\sV^{1/2})$, specified by 
\begin{equation}
  \label{eq:bgl}
J_\sV = U(\tau^G) \quad \mbox{ and } \quad 
\Delta_{\sV}^{-it/2\pi} = U(\exp th) \quad \mbox{ for } \quad t \in \R.
\end{equation}
This assignment is called the {\it Brunetti--Guido--Longo 
(BGL) construction} (see \cite{BGL02}). 
As a consequence, standard subspaces can be associated to antiunitary 
representations in abundance, but only a few of them 
carry interesting geometric 
information. In particular, we would like to understand 
when a standard subspace of the form $\sV_{(h,\tau^G,U)}$ arises from a 
net on some homogeneous space~$G/H$.

A structural property with strong impact in this context is the {\it spectrum condition} 
on the infinitesimal generators $\partial U(x)$ of the one-parameter groups 
$(U(\exp tx))_{t \in \R}$. This is the requirement that the closed convex cone 
\[ C_U := \{ x \in \g \: -i \partial U(x) \geq 0 \} \] 
(the {\it positive cone of $U$}) 
is ``large'' in the sense that the ideal $\g_{C} := C_U - C_U$ 
satisfies \break $\g = \g_{C_U} + \R  h$. 
Let $\sV = \sV_{(h,\tau^G,U)}$ be a standard subspace obtained from the BGL 
construction as in \eqref{eq:bgl}. Then it is a natural requirement 
that the semigroup 
\begin{equation} \label{eq:sv}
  S_\sV = \{ g \in G \: U(g)\sV \subeq \sV\} 
\end{equation}
of endomorphisms of $\sV$ is large in the sense that its Lie wedge 
\[ \L(S_\sV) = \{ x \in \g \: \exp(\R_+ x) \subeq S_\sV \} \] 
(the set of infinitesimal generators of one-parameter subsemigroups of $S_\sV$) 
spans the Lie algebra~$\g$. 
In \cite{Ne21} it is shown that 
if $\L(S_\sV)$ spans $\g$, then 
$\ad h$ defines a $3$-grading in the sense that 
\begin{equation}
  \label{eq:3grad}
\g = \g_1(h) \oplus \g_0(h) \oplus \g_{-1}(h) \quad \mbox{ for } \quad 
\g_j(h) := \ker(\ad h - j\id_\g),
\end{equation}
and 
\begin{equation}
  \label{eq:euler-inv} 
\Ad(\tau^G) = \tau_h,
\end{equation}
so that $\g_0(h)  = \g^{\tau_h}$ and $\g_1(h) + \g_{-1}(h) = \g^{-\tau_h}$. 
We call elements $h$ with \eqref{eq:3grad}  {\it Euler elements} 
because they correspond to the Euler vector field on an open subset 
of the homogeneous spaces with tangent space $\g/(\g_0(h) + \g_{-1}(h))$. 
Assuming that $h$ is an Euler element and $\tau = \tau_h$, 
the semigroup $S_\sV$ has been completely 
determined in \cite{Ne19}. To describe the structure of $S_\sV$, 
let 
\[ C_\pm := \pm C_U \cap \g_{\pm 1}(h)\quad \mbox{ and write } \quad 
G_\sV = \{ g \in G \: U(g) \sV = \sV \}\]
for the stabilizer group of $\sV$ in~$G$. Then  
\begin{equation}
  \label{eq:SE}
S_\sV = \exp(C_+) G_\sV \exp(C_-) = G_\sV \exp(C_+ + C_-).
\end{equation}

Combining Euler elements with symmetric spaces 
leads to the following concept:
A~{\it modular causal symmetric Lie algebra} is a 
quadruple $(\g,\tau, C, h)$, where 
\begin{itemize}
\item $(\g,\tau)$ is a {\it symmetric Lie algebra}, i.e., 
$\g$ is a finite-dimensional real Lie algebra 
and $\tau$ an involutive automorphism of~$\g$. 
\item $(\g,\tau,C)$ is a {\it causal symmetric Lie algebra}, i.e., 
$C \subeq \g^{-\tau}$ is a pointed generating 
closed convex cone 
invariant under the group $\Inn_\g(\fh) = \la e^{\ad \fh}\ra$. 
\item $h \in \g^\tau$ is an {\it Euler element}, i.e., 
$\ad h$ is diagonalizable with eigenvalues $\{-1,0,1\}$, 
and the involution $\tau_h := e^{\pi i \ad h}\in \Aut(\g)$ satisfies 
$\tau_h(C) = - C$.
\end{itemize}

A causal symmetric Lie algebra $(\g,\tau,C)$ is called 
{\it compactly causal} (cc for short) if the cone $C$ is {\it elliptic}, i.e., 
if its interior consists of elements $x$ which are {\it elliptic} 
in the sense that $\ad x$ is semisimple with purely imaginary spectrum. 
Typical examples arise from invariant pointed generating cones $C_\g \subeq \g$. Then 
\begin{equation}
  \label{eq:gt}
(\g \oplus \g, \tau_{\rm flip}, C_\g) \quad \mbox{ with } \quad 
\tau_{\rm flip}(x,y) = (y,x), 
\end{equation}
is a compactly causal symmetric Lie algebra; called {\it of group type (GT)}. 
Another interesting class of compactly causal symmetric 
Lie algebras arises from the involution $\tau_h = e^{\pi i \ad h}$ 
defined by an Euler element. If $C_\g \subeq \g$ satisfies 
$-\tau_h(C_\g) = C_\g$, then 
$(\g,\tau_h, C_\g^{-\tau_h}, h)$ is a modular compactly causal 
symmetric Lie algebra; called of {\it Cayley type (CT)}. 

Irreducible compactly causal symmetric Lie algebras 
$(\g,\tau,C)$ are either of group type with 
$\g^\tau$ simple hermitian, or $\g$ is simple hermitian, and 
$\fz(\fk) \subeq \g^{-\tau}$ holds for a $\tau$-invariant 
Cartan decomposition $\g = \fk \oplus \fp$. 
We refer to \cite{HO97} for a classification. 
The additional assumption that $\fh$ contains an Euler element 
leaves us with the following four types: 
\begin{itemize}
\item[(GT)] Group type: $\g = \fh \oplus \fh$, $\fh$ simple hermitian 
of tube type 
\[ \fh = \su_{r,r}(\C),\ \  \fsp_{2r}(\R), \ \ \so_{2,d}(\R), \ \ 
\so^*(4r),\ \  \fe_{7(-25)}.\] 
\item[(CT)] The Euler element $h$ is central in $\fh$ and $\tau = \tau_h$: 
\[ \begin{tabular}{|l|l|l|l|l|l|}\hline
$\g$ & $\su_{r,r}(\C)$ & $\fsp_{2r}(\R)$ &  $\so_{2,d}(\R), d > 2$&$\so^*(4r)$ 
& $\fe_{7(-25)}$  \\ \hline
$\fh$ &$\R \oplus \fsl_r(\C)$  &$\R \oplus \fsl_r(\R)$ &$\R \oplus \so_{1,d-1}(\R)$  & $\R \oplus \fsl_r(\H)$  & $\R \oplus \fe_{6(-26)}$ \\
\hline 
\end{tabular} \]
\item[(ST)] Split types:  $\tau \not=\tau_h$ and 
${\rm rk}_\R \fh = {\rm rk}_\R \g$:\begin{footnote}{The real rank $\rk_\R(\g)$ 
of a reductive Lie algebra $\g$ is the dimension 
of a maximal abelian $\ad$-diagonalizable subspace.}  
\end{footnote}
\[ \begin{tabular}{|l|l|l|l|l|}\hline
$\g$ & $\su_{r,r}(\C)$ &  $\so_{2,p+q}(\R)$  &$\so^*(4r)$ & $\fe_{7(-25)}$ \\ \hline
$\fh$ &$\so_{r,r}(\R)$ & $\so_{1,p}(\R) \oplus \so_{1,q}(\R)$ 
 & $\so_{2r}(\C)$ & $\fsl_4(\H)$ \\
\hline 
\end{tabular} 
\]
\item[(NST)] Non-split types: $\tau \not=\tau_h$, $\rk_\R \fh 
=\frac{ \rk_\R \g}{2}$:

\[ \begin{tabular}{|l|l|l|l|l|l|l|}\hline
$\g$ & $\su_{2s,2s}(\C)$\hspace{-2mm} & $\fsp_{4s}(\R)$ 
&  $\so_{2,d}(\R)$ 
 \\ \hline
$\fh$ & $\fu_{s,s}(\H)$ & $\fsp_{2s}(\C)$ 
& $\so_{1,d}(\R)$ \\
\hline 
\end{tabular} \]
\end{itemize}

We introduce ``wedge domains''  
in compactly causal symmetric spaces $M = G/H$, 
specified by a modular causal symmetric Lie algebra 
$(\g,\tau,C,h)$ and a connected Lie group~$G$ 
with Lie  algebra $\g$ on which $\tau$ 
induces an involutive automorphism $\tau^G$ 
for which $H \subeq G^{\tau^G}$ is an open subgroup. 
We then obtain a one-parameter group 
$\alpha_t := e^{t \ad h}$ of automorphisms of $\g$,
inducing automorphisms of $G$ and $M$.

Let $V_+(gH) := g.C^0 \subeq T_{gH}(M)$ denote the open cones 
defining the causal structure on $M$ and let 
$X_h^M \in \cV(M)$ be the {\it modular vector field} defined by 
\begin{equation}
  \label{eq:xhdef}
 X_h^M(gH) 
:= \frac{d}{dt}\Big|_{t = 0} \alpha_t(g)H 
=  \frac{d}{dt}\Big|_{t = 0} \exp(th) g H. 
\end{equation}
We introduce the following domains associated to this data: 
\begin{itemize}
\item The {\it positivity domain of the vector filed $X_h^M$ in $M$:} 
\[ W_M^+(h) := \{ m \in M \: X_h^M(m) \in V_+(m) \}.\]

\item The {\it  KMS wedge domain} 
\[ W_M^{\rm KMS}(h) := \{ m \in M \: 
(\forall z \in \cS_\pi)\ \alpha_z(m) \in \cT_M\},\] 
where $(\alpha_z)_{z \in \C}$ denotes the holomorphic extension 
of the modular flow to the complex symmetric space~$M_\C$ 
and 
\[ \cT_M := G \times_H i C^0 \] 
is the {\it tube domain of $M$}, a complex manifold containing 
$M \cong  G \times_H \{0\}$ in its ``boundary''. 
\item The {\it polar wedge domain} 
\[ W_M(h) 
:= \bigcup_{m \in M^\alpha} \Exp_m((C_m^c)^0),\] 
where $M^\alpha$ is the submanifold of fixed points of the modular flow on $M$, 
\[ C_m := \oline{V_+(m)} \subeq T_m(M), \]
and 
\[ C_m^c := C_m \cap T_m(M)_1 - C_m \cap T_m(M)_{-1} \subeq T_m(M),\] 
where $T_m(M)_j$, $j =-1,0,1$ denote the eigenspaces of the 
generator of the modular flow on $T_m(M)$. 
\end{itemize}
Our main geometric result, which is proved 
after various preparations in Section~\ref{sec:5}
(Theorem~\ref{thm:6.5}), asserts that 
these three open subsets coincide: 
\[ W_M^+(h) = W_M^{\rm KMS}(h)  = W_M(h).  \] 
For this result we assume 
that the causal symmetric Lie algebra $(\g,\tau,C)$ 
is {\it extendable}, i.e., there exists a pointed closed convex invariant 
cone $C_\g \subeq \g$ satisfying 
\[ C_\g \cap \fq = C \quad \mbox{ and } \quad -\tau(C_\g) = C_\g.\] 
If $\g$ is reductive and $\fh$ contains no non-zero ideals 
(which is a very natural assumption as these ideals act trivially on $M$), 
then $(\g,\tau,C)$ is extendable (Theorem~\ref{thm:extend}), 
but there are also interesting non-reductive examples. 
As we show  in \cite{NO21b}, the companion paper dealing with wedge domains 
in non-compactly causal symmetric spaces,
the geometry of wedge domains in a non-compactly causal space 
is more complicated because these domains 
are generated from pieces of closed geodesics, the corresponding polar map  
has singularities and it is defined only on a subset of an open~cone.

To connect this geometric insight with nets of standard subspaces, 
we show in our second main result (Theorem~\ref{thm:7.5}) 
that for any antiunitary representation 
$(U,\cH)$ of $G_{\tau_h} = G \rtimes \{\1, \tau_h^G\}$ 
whose positive cone $C_U$ is pointed and generating and satisfies 
$C = C_U \cap \fq$, every cyclic distribution vector 
$\eta \in \cH^{-\infty}$ fixed by $H$ and the conjugation $J = U(\tau_h^G)$ 
leads to a natural covariant isotone net 
$\sH^M$ on $M$ assigning to the connected component 
$W_M(h)_{eH}$ of the wedge domain the standard subspace 
$\sV_{(h,\tau^G_h,U)}$. For classification results concerning this class 
of representations, we refer to \cite{KNO97}. 

Our third main result (Theorem~\ref{thm:8.12}) provides a concrete 
construction of such representations and thus implies that 
such nets exist on all reductive compactly causal modular symmetric spaces. 
Starting with a $C_\g$-positive antiunitary 
representation $(\rho,\cK)$ of $G_{\tau_h}$ which is a 
finite sum of irreducible ones, we obtain on 
the space $B_2(\cK)$ of Hilbert--Schmidt operators on $\cK$ a representation 
by $U(g)A := \rho(g) A \rho(\tau^G(g))^{-1}$ and 
the trace defines a distribution vector which generates a 
subspace $\cH_\rho$, so that Theorem~\ref{thm:7.5} applies to 
$(U,\cH_\rho)$. We even obtain a rather direct description of 
$\sV$ as the closed real subspace generated by 
$\rho(C^\infty_c(S_{\sV,e}^0,\R)) \subeq B_1(\cK) \subeq B_2(\cK)$. 
Here we use that $\rho$ is a trace class representation 
(see Section~\ref{sec:6} for details). 

Motivated by the analysis of abstract wedge spaces  
in \cite{MN21}, we prove in Section~\ref{sec:9} that the wedge 
space $\cW(M,h)$ of all $G$-translates of the connected component 
$W := W_M(h)_{eH}$ actually is a non-compactly causal parahermitian 
symmetric space. It is a covering of the adjoint orbit 
$\cO_h = \Ad(G)h \subeq \g$. The most difficult part of the proof 
is to get hold of the stabilizer group  $G_W = \{ g \in G \: gW = W\}$. 
It is an interesting open problem to connect the covering theory 
and the corresponding twisted duality developed in \cite{MN21} 
to the geometry of coverings of compactly causal symmetric spaces~$G/H$.  

We conclude this paper with a section on the anti de Sitter space 
$\AdS^d$, the prototypical Lorentzian example of a compactly causal symmetric space. 
This section can be read independently. 
For this example we evaluate all geometric data explicitly and show directly 
that the three types of wedge domains coincide, without referring to the elaborate 
structure theory used in the rest of the paper. 
\\

\nin {\bf Structure of this paper:} 
We start in Section~\ref{sec:2} by recalling basic facts on 
causal symmetric spaces 
on the Lie algebra level (Subsections~\ref{subsec:2.1} to 
~\ref{subsec:redliealg}) and also on the global level 
(Subsection~\ref{subsec:2.3}). 

In Subsection~\ref{sec:3new.1} we classify for a 
reductive compactly causal symmetric reductive Lie algebra $(\g,\tau)$ 
the orbits of $\Inn_\g(\fh)$ in the set 
of Euler elements in $\fh$. As this set describes the different choices 
of the Euler element $h \in \fh$, this can be seen as a classification 
of all ``modular structures'' on the associated symmetric spaces. 
In Subsection~\ref{sec:3new.2} we show that the orbits 
of the centralizer group $G^h$ in $M^\alpha$ are open and that this 
leads to a natural bijection between the orbit spaces $M^\alpha/G^h$ 
and $(\cO^h \cap \fh)/H$. Actually both sets correspond to the set of 
double cosets $G^h g H \subeq G$ for which $\Ad(g)^{-1}h \in \fh$. 

In Section~\ref{sec:3} we introduce the three types of 
wedge domains in a compactly causal symmetric space~$M =  G/H$. 
In our analysis of wedge domains we 
follow the strategy to first study spaces of group 
type and then use embeddings into these spaces to derive corresponding 
results in general. In Section~\ref{sec:4} we start with 
causal symmetric spaces of group type, i.e., pairs 
$(G,C_\g)$ of a connected Lie group $G$ and a 
pointed generating invariant cone $C_\g \subeq \g$. 

In Section~\ref{sec:5} we turn to wedge domains in compactly 
causal symmetric spaces. In this context we 
assume that $(\g,\tau,C)$ is extendable in the sense of
Subsection~\ref{subsec:2.1}, keeping in mind that this is always 
the case if $\g$ is reductive. The main result of this section 
is Theorem~\ref{thm:6.5}, asserting that the three wedge domains in $M$ coincide. 
This is first proved for the special case 
$H = G^{\tau^G}$ in Theorem~\ref{thm:6.4}. 
The proof of this theorem builds heavily on the group case 
(Theorem~\ref{thm:4.2}). We conclude this section with a brief 
discussion of the assumption that the Euler element $h$ is contained 
in $\fh$, i.e., that the corresponding modular flow on $M$ has a fixed point. 
Note  that the definition of $W_M(h)$ makes no sense if $M^\alpha = \eset$. 

In Section~\ref{sec:6a} we turn to representation theoretic aspects of 
compactly causal symmetric spaces. We start with 
introducing standard subspaces, the Brunetti--Guido--Longo construction 
and recall the construction of nets of closed real subspaces 
from distribution vectors,  introduced in \cite{NO21a}. 
In Subsection~\ref{subsec:7.4} we describe in the general 
Theorem~\ref{thm:7.5} how the methods from \cite{NO21a} can be used 
to construct covariant nets of standard subspaces on 
compactly causal symmetric spaces. 
In Section~\ref{sec:6} we construct such representations explicitly 
in spaces of Hilbert--Schmidt operators $\cH_\rho \subeq B_2(\cK)$, where 
$(\rho, \cK)$ is an antiunitary representation of $G_{\tau^G_h}$ which is a 
finite sum of irreducible representations. This is done in three steps: 
First we recall from \cite{Ne00, Ne19} some results on 
the representations $(\rho,\cK)$ of $G_{\tau_h} = G \rtimes \{\1,\tau^G_h\}$, 
then we use these representations to construct nets of standard subspaces 
on the symmetric space of group type $G \cong (G \times G)/\Delta_G$, 
and finally we use the twisted embedding 
$G \to G \times G, g \mapsto (g, \tau^G(g))$ to obtain by 
pullbacks representations of $G_{\tau_h}$ that can be used to obtain with 
Theorem~\ref{thm:7.5} nets of standard subspaces on $M = G/H$. 

In Section~\ref{sec:9} we return to a geometric topic. 
We show that, under rather natural assumptions, 
the wedge space $\cW(M,h)$ of all $G$-translates 
of the connected  component $W_M(h)_{eH}$ of the wedge domain in $M$  
carries the structure of an ordered symmetric space.  
We conclude this paper with a short section discussing 
some open problems and two appendices  on some useful facts on 
symmetric Lie groups and on distribution vectors of unitary representations. \\

\noindent 
\nin {\bf Notation:} 
\begin{itemize}
\item $\g$, $\fh$ and $\fk$ will denote Lie algebras 
of the Lie groups $G$, $H$, $K$.
\item An element $h\in\g$ is called an {\it Euler element} if
$\ad h$ is non-zero and diagonalizable 
with $\Spec(\ad h) \subeq \{-1,0,1\}$. The set of Euler elements in $\g$ is
denoted by $\cE (\g)$.
\item  For $h\in \cE (\g)$ we write 
$\tau_h := e^{\pi i \ad h}$ for the corresponding involution on $\g$ and 
note that $\kappa_h := e^{-\frac{\pi i}{2} \ad h}$ is an automorphism of 
$\g_\C$ of order~$4$ with $\kappa_h^2 = \tau_h$ which 
implements the ``Wick rotation'' (partial Cayley transform). 
\item   For $h \in \g$, $\lambda \in \R$, and $\fl\subset \g$ an $\ad h$ invariant
subspace  we write 
$\fl_\lambda(h) := \ker(\ad h|_{\fl} - \lambda \1)$ for the corresponding eigenspace 
in the adjoint representation.
\item We write $e \in G$ for the identity element in the Lie group~$G$, 
$G_e$ for its identity component and $Z(G)$ for its center. 
\item For the left and right translations on the tangent bundle $T(G)$, 
we write $g.v$ and $v.g$ for $g \in G, v \in T(G)$, respectively. 
\item For $x \in \g$, we write $G^x := \{ g \in G \: \Ad(g)x = x \}$ 
for the stabilizer of $x$ in the adjoint representation 
and $G^x_e := (G^x)_e$ for its identity component. 
\item If $\tau^G \in \Aut(G)$ is an involution, then 
we write $G_\tau := G \rtimes \{\id_G, \tau^G\}$, 
$g^\sharp := \tau^G(g)^{-1}$, and $G^\sharp := \{ g \in G \: g^\sharp = g\}$. 
\item For a Lie subalgebra $\fh \subeq \g$, we write 
$\Inn_\g(\fh) = \la e^{\ad \fh}\ra \subeq \Aut(\g)$ for the corresponding group of inner 
automorphism and put $\Inn(\g) := \la e^{\ad \g} \ra$. 
\item For a (continuous) unitary representation $(U, \cH)$ of a Lie group 
$G$ with Lie algebra $\g$, we write 
$\partial U(x)$ for the skew-adjoint infinitesimal generator of 
the unitary one-parameter group 
$U(\exp t x)$, so that we have 
$U(\exp t x) = e^{t \partial U(x)}$ for $t \in \R$. 
The {\it positive cone of $U$} is the $\Ad(G)$-invariant closed convex cone 
\begin{equation}
  \label{eq:wpi-intro} 
C_U := \{ x \in \g \: -i \partial U(x) \geq 0\}.
\end{equation}
\item For a complex Hilbert space, we write $\AU(\cH)$ for the 
group of unitary and antiunitary operators. 
\item For a function $f \: G \to \C$, we write 
$f^\vee(g) := f(g^{-1})$. 
\end{itemize}

\section{Causal symmetric Lie algebras} 
\mlabel{sec:2}

In this section we collect some basic facts and definitions concerning 
causal symmetric spaces. Our main reference for these spaces is \cite{HO97} for the
semisimple situation and \cite{HN93} for the general case.

\begin{defn} \mlabel{def:symliealg}
For a symmetric Lie algebra $(\g,\tau)$, we write 
\[ \fh := \g^\tau = \{ x \in \g \: \tau(x) = x\}\quad \mbox{ and } \quad 
  \fq := \g^{-\tau} = \{ x \in \g \: \tau(x) = -x\}.\] 
If $\g$ is semisimple, then 
there exists a Cartan involution $\theta$ commuting with $\tau$, 
and all these Cartan involutions are mutually conjugate under 
$\Inn_\g(\fh)$ (\cite[Prop.~I.5(iii)]{KN96}). 
If $\theta$ is a Cartan involution commuting with~$\tau$, 
then we write $\fk := \g^\theta$ and $\fp := \g^{-\theta}$ for the 
$\theta$-eigenspaces. Now 
$\g$ decomposes as follows into simultaneous $(\tau,\theta)$-eigenspaces: 
\begin{equation}
  \label{eq:foursubspaces}
\g = \fh \oplus \fq = \fk \oplus \fp= \fh_\fk \oplus \fh_\fp \oplus \fq_\fk \oplus \fq_\fp
\end{equation}
with
\[
\fh_\fk := \fh \cap \fk, \quad 
\fh_\fp := \fh \cap \fp, \quad 
\fq_\fk := \fq \cap \fk, \quad  
\fq_\fp := \fq \cap \fp.
\]
We call 
$(\g^c,\tau^c)$ with $\g^c = \fh + i \fq$ and 
$\tau^c(x + iy) = x - iy$ for $x \in \fh, y \in \fq$, the {\it $c$-dual 
symmetric Lie algebra}. 
\end{defn}

\begin{defn} Let $(\g,\tau)$ be a symmetric Lie algebra. 
  \begin{footnote}{
Here we extend concepts typically used for irreducible semisimple 
symmetric Lie algebras to a more general context. This is motivated by 
our heavy use of embeddings, which is most nicely done in a suitable 
functorial framework.}
  \end{footnote}
  \begin{itemize}
  \item[\rm(a)] A triple $(\g,\tau_h,h)$ is called a {\it parahermitian 
symmetric Lie algebra} if $h \in \cE(\g)$ is an Euler element and 
$\tau_h = e^{\pi i \ad h}$. 
\item[\rm(b)] A triple $(\g,\tau,C)$ is called a {\it 
compactly causal (cc) symmetric Lie algebra} if 
$C \subeq \fq$ is a pointed generating elliptic 
cone invariant under $\Inn_\g(\fh) := \la e^{\ad \fh}\ra$.
\item[\rm(c)] A triple $(\g,\tau,C)$ is called a {\it 
non-compactly causal (ncc) symmetric Lie algebra} if 
$C \subeq \fq$ is a pointed generating hyperbolic 
cone invariant under $\Inn_\g(\fh)$, i.e., the interior of $C$ consists of 
$\ad$-diagonalizable elements. 
\item[\rm(d)] A triple $(\g,\tau,C)$ is called a {\it causal symmetric Lie algebra} if it is either a cc or an ncc symmetric Lie algebra. 
\item[\rm(e)] 
A {\it modular causal symmetric Lie algebra} is a 
quadruple $(\g,\tau, C, h)$, where 
$(\g,\tau,C)$ is a causal symmetric Lie algebra, 
$h \in \g^\tau$ is an Euler element,  
and the involution $\tau_h = e^{\pi i \ad h}\in \Aut(\g)$ satisfies 
$\tau_h(C) = - C$.
  \end{itemize}
Note that the $c$-dual of a cc symmetric Lie algebra $(\g,\tau,C)$ is 
ncc and vice versa.  
\end{defn}

For the classification 
of $\Inn(\fh)$-invariant elliptic cones $C \subeq \fq$, i.e., causal structures on the symmetric Lie algebra 
$(\g,\tau)$,
see Sections 4.4 and 4.5 in \cite{HO97} and  \cite[Thm.~X.3]{KN96}. Note that the
classification in the cc case and ncc case are the same via   $c$-duality.

\begin{rem} \mlabel{rem:cpm}
Suppose that $C \subeq \fq$ is a closed convex $\Inn_\g(\fh)$-invariant 
cone with $-\tau_h(C) = C$. Since $\ad h$ preserves 
$\fq$, the space $\fq$ decomposes as 
\[ \fq = \fq_{-1}(h) \oplus \fq_0(h) \oplus \fq_1(h).\]
As $C$ is invariant under $e^{\R \ad h}$ and $-\tau_h$, we see that 
$x= x_1+x_0+x_{-1}\in C$ with $x_j \in \fq_j(h)$ implies 
\[ e^{t \ad h} x = e^tx_1 +x_0 +e^{-t}x_{-1} \in C.\]
This shows that $x_{\pm 1} = \lim_{t\to \infty} e^{-t} e^{\pm t \ad h} x\in C$. 
Therefore the two 
closed convex cones 
\begin{equation}
  \label{eq:cpm}
 C_\pm := \pm C \cap \fq_{\pm 1}(h) 
\end{equation}
satisfy 
\begin{equation}
  \label{eq:coneproj}
p_{\fq^{-\tau_h}}(C) = C^{-\tau_h} =  C_+ - C_-. 
\end{equation}

Furthermore, the automorphism $\kappa_h =  e^{-\frac{\pi i}{2} \ad h} \in \Aut(\g_\C)$ 
takes the form 
\[ \kappa_h(z_{-1} + z_0 + z_1) = i z_{-1} + z_0 -i z_1 \]
with respect to the $3$-grading defined by $\ad h$ and
we have 
\begin{equation}
  \label{eq:zetarel2}
 i \kappa_h(C^{-\tau_h}) = C_+ + C_-  =: C^c
\end{equation}
\end{rem}

\subsection{Embedding into group type} 
\mlabel{subsec:2.1}

In this section we discuss the embedding of compactly causal triples
$(\g,\tau ,C)$ into a natural causal triple of group
type $( (\fg \oplus \fg)_C, \tau_{\rm flip}, C_\fg)$. This embedding will play
a crucial role in the following.

We call a compactly causal symmetric Lie algebra  
$(\g,\tau,C)$ 
{\it extendable} if there exists a {\it pointed} closed convex invariant 
cone $C_\g \subeq \g$ satisfying 
\[ C_\g \cap \fq = C \quad \mbox{ and } \quad -\tau(C_\g) = C_\g.\] 
We then call $(\g,\tau,C_\g)$ an {\it extension of $(\g,\tau,C$)}. 
The following theorem (\cite[Thm.~X.7]{KN96}, \cite[Thm.~4.5.8]{HO97}, 
\cite[Thm.~7.8]{O91}) shows that, for reductive Lie algebras, 
extensions always exist, even with generating cones $C_\g$. 

\begin{thm}{\rm(Extension Theorem)} \mlabel{thm:extend} 
Let $(\g,\tau,C)$ be a reductive compactly causal symmetric Lie algebra 
for which $\fh$ contains no non-zero ideals of $\g$. 
Then there exists a pointed generating invariant closed 
convex cone $C_\g \subeq \g$ with 
\begin{equation}
  \label{eq:dagg-145}
 -\tau(C_\g) = C_\g \quad \mbox{ and } \quad C_\g \cap \fq = C.
\end{equation}
\end{thm}

Suppose that $(\g,\tau,C_\g)$ is an extension 
of the compactly causal symmetric Lie algebra $(\g,\tau,C)$. 
If $C_\g$ is not generating, then 
\[ \g_C := C_\g - C_\g \trile \g \] 
is a $\tau$-invariant ideal containing $C - C = \fq$, so that 
$\g_C = \fh_C \oplus \fq$ for $\fh_C = \fh \cap \g_C$. 
Then 
\[ (\g \oplus \g)_C = \{ (x,y) \in \g \oplus \g \: x -y \in \g_C \} 
\quad \mbox{ with } \quad \tau_e(x,y) = (y,x)\] 
is a symmetric Lie algebra with 
\[ (\g \oplus \g)_C^{\tau_e}=\{(x, x)\: x\in \g\}\cong \g   
\quad\text{and}\quad (\g \oplus \g)_C^{-\tau_e} = \{ (x,-x) \: x \in \g_C \} \cong \g_C
.\] 
We thus obtain an embedding of causal symmetric Lie algebras 
\begin{equation}
  \label{eq:grptypemb}
 (\g,\tau,C) \into ((\g \oplus \g)_C, \tau_e, C_\g), \quad 
x  \mapsto (x, \tau x) \quad \mbox{ for }\quad 
x \in \fg .
\end{equation}
This embedding is also compatible with the corresponding Euler elements. 
Any Euler element $h \in \fh = \g^\tau$ is mapped 
to the Euler element $(h,h) \in (\g \oplus \g)_C^{\tau_e}$.

On the global level, the embedding 
\eqref{eq:grptypemb} corresponds to the 
quadratic representation of the symmetric space $G/H$ in the group 
$G_C$ with Lie algebra $\g_C$ (see \eqref{eq:quademb} below). 
We shall use such embeddings to derive 
properties of compactly causal symmetric spaces from those of 
spaces of group type. 

Since $\tau$ commutes with $\ad h$ on $\g$, it preserves all eigenspaces 
$\g_j(h)$, $j =-1,0,1$. 
We conclude from \eqref{eq:dagg-145} that the cones 
\[ C_{\g,\pm} := \pm C_\g \cap \g_{\pm 1}(h) \] 
satisfy $- \tau(C_{\g,\pm}) = C_{\g,\pm}$, which shows in particular that 
the cone 
\begin{equation}
  \label{eq:cgc}
 C_\g^c := C_{\g,+} + C_{\g,-} 
\end{equation}
is $-\tau$-invariant with 
\begin{equation}
  \label{eq:doubledcone2}
(C_\g^c)^{-\tau} = 
(C_{\g,+} + C_{\g,-})^{-\tau} = C_+ + C_- = C^c. 
\end{equation}

\begin{rem}
In \cite[Ex.~X.4]{KN96} one finds examples of 
compactly causal symmetric Lie algebras $(\g,\tau,C)$ 
which are not extendable, but \cite[Thm.~X.7]{KN96} describes sufficient 
conditions for general Lie algebras that cover all cases where 
$C_\g$ is generating. 

We also note that, if $\g$ is not reductive 
and $C_\g \subeq \g$ is a pointed generating invariant cone, 
then $C_\g \cap \fz(\g) \not=\{0\}$. As $\tau_h$ fixes the center 
pointwise, this condition is incompatible with 
$-\tau_h(C_\g) = C_\g$. Therefore it would be too much to hope 
for extensions of modular compactly causal symmetric Lie algebras 
$(\g,\tau,C,h)$ for which $C_\g$ is also generating. 
We refer to \cite{Oeh20, Oeh20b} 
for more details on the structure of non-reductive modular 
compactly causal symmetric Lie algebras $(\g,\tau, C,h)$ 
and corresponding classification results. 
\end{rem}

\begin{ex} \mlabel{ex:2.6} 
(A non-reductive example; cf.\ \cite[Ex.~3.7]{Ne19}) 
We consider the Lie algebra 
\[ \g = \hcsp(V,\omega) := \heis(V,\omega) \rtimes \csp(V,\omega), \] 
where 
$(V,\omega)$ is a symplectic vector space, 
$\heis(V,\omega) = \R \oplus V$ is the corresponding Heisenberg algebra 
with the bracket $[(z,v),(z',v')] = (\omega(v,v'),0)$, and 
\[ \csp(V,\omega) := \sp(V,\omega) \oplus \R \id_V \] 
is the {\it conformal symplectic Lie algebra} of $(V,\omega)$. 

The hyperplane ideal $\fj := \heis(V,\omega) \rtimes \sp(V,\omega)$ 
(the {\it Jacobi algebra})  
can be identified by the linear isomorphism 
\[ \phi \colon \fj \to \Pol_{\leq 2}(V), \qquad 
\phi(z,v,x)(\xi) := z + \omega(v,\xi) + \frac{1}{2} \omega(x\xi,\xi), \quad 
\xi \in V \] 
with the Lie algebra of real polynomials 
$\Pol_{\leq 2}(V)$ of degree $\leq 2$ on $V$, 
endowed with the Poisson bracket (\cite[Prop.~A.IV.15]{Ne00}).  
The set 
\[ C := \{ f \in \Pol_{\leq 2}(V) \colon f \geq 0 \}  \] 
is a pointed invariant cone in $\g$ generating the ideal $\g_C = \fj$. 
The element $h_0 := \id_V$ defines a derivation on $\fj$ by 
$(\ad h_0)(z,v,x) = (2z,v,0)$ for $z \in \R, v \in V, x \in \sp(V,\omega)$. 
Any involution $\tau_V$ on $V$ satisfying $\tau_V^*\omega = - \omega$ 
defines by 
\begin{equation}
  \label{eq:invtau}
\tilde\tau_V(z,v,x) := (-z,\tau_V(v), \tau_V x \tau_V) 
\end{equation}
an involution on $\g$ with $\tilde\tau_V(h_0) = h_0$,  
and $-\tilde\tau_V(C) = C$ follows from 
\[ \phi(\tilde \tau_V(z,v,x)) = - \phi(z,v,x) \circ \tau_V.\]
Considering $\tau_V$ as an element of $\sp(V,\omega)$, 
we obtain an Euler element 
\[ h := \shalf(\id_V + \tau_V)\in \csp(V,\omega).\] 
Writing $V = V_1 \oplus V_{-1}$ for 
the $\tau_V$-eigenspace decomposition, the $\ad h$-eigenspaces in $\g$ are 
\[ \g_{-1}  = 0 \oplus 0 \oplus \sp(V,\omega)_{-1}, \quad 
\g_0 = 0 \oplus V_{-1}\oplus \sp(V,\omega)_0 
\cong V_{-1} \rtimes \gl(V_{-1}), \quad 
\g_1 = \R \oplus V_1\oplus \sp(V,\omega)_1.\] 
The eigenspace $\g_1$ can be identified with the space 
$\Pol_{\leq 2}(V_{-1})$ of polynomials of 
degree $\leq 2$ on $V_{-1}$ and 
\[ C_+ = C \cap \g_1 = \{ f \in \Pol_{\leq 2}(V_{-1}) \colon f \geq 0\}.\] 
This  cone is invariant 
under the natural action of the 
affine group 
\[ G^0 \cong \Aff(V_{-1})_e \cong V_{-1} \rtimes \GL(V_{-1})_e \]
whose Lie algebra is $\g_0$. We also note that 
\[ \g_{-1} = \sp(V,\omega)_{-1} \cong \Pol_2(V_1) \quad \mbox{ and } \quad 
C_- = - C \cap \g_{-1} = \{ f \in \Pol_2(V_1) \colon f \leq 0\}.\] 

Finally, we observe that $\tau_h =  e^{\pi i \ad h} = (-\tau_V)\,\tilde{}$\ \ implies in particular $-\tau_h(C) = C$, so that 
\[ (\g,\tau_h, C^c, h)\quad \mbox{ with } \quad C^c = C_+ - C_- = C^{-\tau_h} \] 
is a modular causal symmetric Lie algebra for which 
$(\g, \tau_h, C^c)$ is compactly causal and extendable to $(\g,\tau_h,C)$. 
\end{ex}

\subsection{Hermitian simple Lie algebras} 
\mlabel{subsec:2.2}

A real simple Lie algebra $\g$ contains a 
pointed generating invariant cone $C_\g$ if and only if 
$\g$ is hermitian (\cite{Vi80}). If this is the case, then 
there exist two such cones $C_\g^{\rm min} \subeq C_\g^{\rm max}$ with the property 
 that 
any other pointed generating invariant cone $C_\g$  satisfies 
\begin{equation}
  \label{eq:min-max-g}
 C_\g^{\rm min} \subeq C_\g\subeq C_\g^{\rm max} 
\quad \mbox{ or }\quad 
 C_\g^{\rm min} \subeq -C_\g\subeq C_\g^{\rm max}.
\end{equation}

We start by collecting some information on the relevant class of simple 
Lie algebras. 

\begin{prop} \mlabel{prop:cc2}
Let $\g$ be a simple hermitian Lie algebra. Then 
$\g$ contains an Euler element if and only if $\g$ is of tube type. 
If this is the case, then the following assertions hold: 
\begin{itemize}
\item[\rm(a)] $\Inn(\g)$ acts transitively on the set $\cE(\g)$ of Euler elements. 
\item[\rm(b)] If $(\g,\tau   , C)$ is  compactly causal , then $\cE(\g) \cap \fh \not= \eset$. 
\item[\rm(c)] For every Euler element $h\in \g$, 
the 
cone $C^{\rm max}_\g \subeq \g$ satisfies 
$\tau_h(C^{\rm max}_\g) = - C^{\rm max}_\g$. In particular $(\g,\tau_h)$ is 
compactly causal. 
\item[\rm(d)] For any pointed generating invariant cone $C_\g \subeq C^{\rm max}_\g$, 
we have 
$C_\g \cap \g^{-\tau_h} = C^{\rm max}_\g \cap \g^{-\tau_h}$. 
\end{itemize}
\end{prop}

\begin{prf} (a) follows from \cite[Prop.~3.11]{MN21}. 

\nin (b) Assume that $(\g,\tau)$ is compactly causal.  
In \cite[App.~D]{NO21b} we show that 
$\tau$ preserves a $3$-grading of $\g$, and this means that the 
corresponding Euler element $h$ is fixed by $\tau$, i.e., $h \in \fh$. 

\nin (c) Since $\tau_h(C_\g^{\rm max})$ is one of the two maximal 
invariant cones $\pm C_\g^{\rm max}$, 
(c) follows from \cite[Lemma~2.28]{Oeh20}. 


\nin (d) Clearly, $C^{\rm max}_\g \cap \g^{-\tau_h} \supeq C_\g \cap \g^{-\tau_h}$ 
by \eqref{eq:min-max-g}. 
Writing $\g^{-\tau_h} = \g_1(h) \oplus \g_{-1}(h)$ and using the invariance 
under $e^{\R \ad h}$, it follows as in Remark~\ref{rem:cpm} that 
\begin{equation}
  \label{eq:pmdec}
 C_\g \cap \g^{-\tau_h} = C_{\g,+} + - C_{\g,-}, \quad \mbox{ where } \quad  
C_{\g,\pm} = \pm C_\g \cap \g_{\pm 1}(h).
\end{equation}
The cones $C_{\g,\pm}$ and $C^{\rm max}_{\g,\pm}$ 
coincide by \cite[Lemma~2.28]{Oeh20}, so that 
(d) follows from \eqref{eq:pmdec}. 
\end{prf}

\begin{prop} \mlabel{prop:eul-tt} 
Let $\g$ be simple and $h \in\cE(\g)$.  
Then $(\g,\tau_h)$ is causal if and only if 
it is of Cayley type, i.e., $\g$ is simple hermitian  of tube type.
\end{prop} 

In this case any causal structure represented 
by a cone $C \subeq \fq = \g^{-\tau_h}$ 
leads to two cones $C = C_+ - C_-$ and $C^c = C_+ + C_-$. 
One is hyperbolic and the other is elliptic. 

\begin{prf} We choose a Cartan involution $\theta$ with 
$\theta(h) = -h$. 

If $(\g,\tau_h)$ is causal and 
$C \subeq \fq$ is a pointed generating 
$\Inn_\g(\fh)$-invariant cone, then 
$C^0$ contains fixed points of the compact group 
$\Inn_\g(\fh_\fk)$, so that $\fq^{\fh_\fk} \not=\{0\}$. 
Hence the conclusion follows from \cite[Thm.~1.3.11]{HO97}. 

If, conversely, $\g$ is simple hermitian and $h \in \cE(\g)$, then 
Proposition~\ref{prop:cc2}(c) implies that $(\g,\tau_h)$ is causal.   
\end{prf}

\subsection{Reductive Lie algebras} 
\mlabel{subsec:redliealg}

The last section was devoted to the semisimple case. In this section
we discuss the reductive case. {\it We assume that $\fh$ does not contain any
non-zero ideal.} This is no restriction because the integral subgroup corresponding to 
such an ideal would act trivially on $G/H$ and thus 
can be factorized out. This in particular implies
that $\z (\g) \cap \fh =\{0\}$ because $\z (\g)\cap \fh$ is an ideal of~$\g$. Thus
$H\cap Z(G)$ is discrete. 

The following proposition extends Proposition~\ref{prop:cc2}(c) to 
the larger class of reductive Lie algebras. 

\begin{prop} If $(\g,\tau)$ is a 
reductive compactly causal symmetric Lie algebra 
containing an Euler element and $\fh$ contains no non-zero ideal, 
then $\cE(\g) \cap \fh \not=\eset$. 
\end{prop}

\begin{prf} If $h \in \g$ is an Euler element and 
$h = h_0 + h_1$ with $h_0 \in \fz(\g)$ and $h_1 \in [\g,\g]$, 
then $h_1 \in \cE([\g,\g])$. As $(\g,\tau)$ is compactly causal 
and $\fh$ contains no non-zero ideal, all simple 
ideals of $\g$ are either hermitian or compact. 
In fact, the Extension Theorem~\ref{thm:extend} implies the 
existence of a pointed generating invariant cone in $\g$, 
hence $\g$ is quasihermitian by \cite[Cor.~VII.3.9]{Ne00}, 
i.e., $\fk = \fz_\g(\fz(\fk))$ holds for a maximal compactly embedded subalgebra. 
Since $\fk$ is adapted to the decomposition into simple ideals $\g_j$, 
the ideal $\g_j$ is compact if $\fz(\fk_j) = \{0\}$ and 
it is hermitian if $\fz(\fk_j) \not=\{0\}$ (by definition).

As $\ad h_1$ is  diagonalizable, $h_1$ is contained in the sum of all 
hermitian simple ideals. Any such ideal is either $\tau$-invariant, 
or corresponds to an irreducible summand of group type.
For group type Lie algebras $(\g \oplus \g, \tau_{\rm flip})$ 
(cf.\ \eqref{eq:gt}), the assertion 
is trivial because $h = (h_1, h_2) \in \cE(\g \oplus \g)$ 
implies $(h_j,h_j) \in \fh$ for $j =1,2$, and at least one of these two 
elements is non-central, hence an Euler element. 
For the $\tau$-invariant ideals, the assertion follows from 
Proposition~\ref{prop:cc2}(b).
\end{prf}

Let $(\g,\tau,C)$ be a reductive compactly causal symmetric Lie algebra 
with $\fz(\g) \subeq \fq$ (which follows from $\fh \cap \fz(\g) = \{0\}$).
Then $\g = \fz(\g) \oplus \bigoplus_{j = 1}^n \g_j$, 
where the $\g_j$ are simple ideals that are either compact or simple hermitian. 
If $\g_j$ is compact, we put $C_{\g_j}^{\rm min} = \{0\}$ and 
$C_{\g_j}^{\rm max} = \g_j$. With \eqref{eq:min-max-g} 
this leads to invariant cones 
\begin{equation}
  \label{eq:min-max-g-red}
 C_\g^{\rm min} 
:= \bigoplus_{j =1}^n C_{\g_j}^{\rm min} 
\subeq C_\g^{\rm max} := \fz(\g) \oplus \sum_{j =1}^n C_{\g_j}^{\rm max}.
\end{equation}
Every Euler element $h \in \g$ decomposes as 
$h = h_0 + \sum_{j = 1}^n h_j$ with 
$h_0 \in \fz(\g)$ and $h_j = 0$ if $\g_j$ is compact, 
and if $h_j\not=0$, then it is an Euler element in $\g_j$. 
From Proposition~\ref{prop:cc2} we now obtain: 

\begin{cor} \mlabel{cor:red-minmax}
Let $\g$ be a reductive quasihermitian Lie algebra, 
i.e., all simple ideals are compact or hermitian, 
and $h \in \cE(\g)$.Then the following assertions hold: 
\begin{itemize} 
\item[\rm(a)] If $\g^{\tau_h} = \ker(\ad h)$ contains no simple hermitian ideal, 
then $\tau_h(C^{\rm max}_\g) = - C^{\rm max}_\g$. 
\item[\rm(b)] For any pointed generating invariant cone 
$C_\g \subeq C^{\rm max}_\g$, we have 
$C_\g \cap \g^{-\tau_h} = C^{\rm max}_\g \cap \g^{-\tau_h}$. 
\end{itemize}
\end{cor}

\begin{prf} (a) Our assumption on $h$ means that 
$h_j \in \cE(\g_j)$ for every simple hermitian ideal~$\g_j$ because 
$h_j \not=0$. 
Hence (a) follows immediately from \eqref{eq:min-max-g-red} 
and the corresponding assertion for hermitian 
simple Lie algebras (Proposition~\ref{prop:cc2}). 
  
\nin (b) As $\g^{\tau_h}$ contains the center and all compact ideals, 
$\g^{-\tau_h}$ is contained in the sum 
\[ \g_{\rm herm} := \g_{j_1} \oplus \cdots \oplus \g_{j_k} \] 
of all hermitian simple ideals. Further, the invariance of all ideals 
under $\tau_h$ shows that 
\[ \g^{-\tau_h} =  \g_{j_1}^{-\tau_h} \oplus \cdots \oplus \g_{j_k}^{-\tau_h}. \] 
The classification of invariant  cones 
(see \cite[Thm.~VIII.3.21]{Ne00} for more details) implies that 
$C_{\g, {\rm herm}} := C_\g \cap \g_{\rm herm}$ is a 
pointed generating invariant cone in $\g_{\rm herm}$ 
with 
\[ \sum_{\ell = 1}^k C_{\g_{j_\ell}}^{\rm min} 
\subeq C_{\g, {\rm herm}} \subeq 
C_{\g, {\rm herm}}^{\rm max} 
= \sum_{\ell = 1}^k C_{\g_{j_\ell}}^{\rm max}.\] 
Therefore Proposition~\ref{prop:cc2}(d) leads to  
\[ C_\g^{\rm max} \cap \g^{-\tau_h} \supeq 
C_\g^{-\tau_h} 
\supeq \sum_{\ell = 1}^k C_{\g_{j_\ell}}^{\rm min} \cap \g_{j_\ell}^{-\tau_h}
= \sum_{\ell = 1}^k C_{\g_{j_\ell}}^{\rm max} \cap \g_{j_\ell}^{-\tau_h}
= C_\g^{\rm max} \cap \g^{-\tau_h}.\qedhere\] 
\end{prf}

\subsection{The global setting for symmetric spaces} 
\mlabel{subsec:2.3} 

In the following $G$ denotes a connected real Lie group 
and $\eta_G \: G \to G_\C$ its universal 
complexification. {\it We assume that $\ker(\eta_G)$ is discrete,} 
which is always the case if $\g$ is semisimple 
(because  $\ker\eta_G \subeq Z(G)$),  if $G$ is simply connected, 
or if $G$ is a matrix group. 
The groups $G$ and $G_\C$ carry natural structures 
 of Loos symmetric spaces (\cite{Lo69}), defined by 
 \begin{equation}
   \label{eq:loosprod}
 g \bullet h := s_g(h) := g h^{-1} g, 
 \end{equation}
Then automorphisms, antiautomorphisms and left and right translations define 
automorphisms of the symmetric space $(G,\bullet)$. In particular, we 
obtain a transitive action of the product group $G \times G$ on $G$ 
by automorphisms of $(G,\bullet)$ via 
\[ (g_1,g_2).g := g_1 g g_2^{-1}.\] 
Then $T_e(G) \cong \{(x,-x) \: x \in \g\} \subeq \g \oplus \g$, and 
the exponential function of the pointed symmetric space 
$(G,\bullet, e)$ is given by 
\begin{equation}
  \label{eq:Exp}
 \Exp_e \:  T_e(G) \to G, \quad (x,-x) \mapsto \exp(2x)= (\exp x, \exp -x).e.
\end{equation}

If $\tau^G$ is an involutive automorphism of $G$ and 
$H \subeq G^{\tau^G}$ an open subgroup, then 
\begin{equation}
  \label{eq:mgh}
 M := G/H, \qquad g_1 H \bullet g_2 H := g_1 \tau^G(g_2)^{-1} g_1 H 
\end{equation}
is the corresponding {\it symmetric space}. 
Its {\it quadratic representation} is the map 
\begin{equation}
  \label{eq:quademb}
Q \: M \to G, \quad Q(gH) := gg^\sharp, \qquad g^\sharp := \tau^G(g)^{-1}. 
\end{equation}
It is a covering of the identity component $G^\sharp_e \cong G/G^{\tau^G}$ 
of the symmetric subspace 
\[ G^\sharp := \{ g \in G \: g^\sharp =g \} \subeq G\] 
(see Lemma~\ref{lem:a.1} in Appendix~\ref{app:1}). 

We write the natural action of $G$ on the tangent bundle 
$T(M)$ by $(g,v) \mapsto g.v$ and identify the tangent space 
$T_{eH}(M)$ with $\fq$. 
If $(\g,\tau,C)$ is a causal symmetric Lie algebra for which 
$\Ad(H) C = C$, then 
\[ V_+(gH) = g.C^0 \subeq T_{gH}(M) \] 
defines on $M = G/H$ a $G$-invariant cone field, i.e., a causal structure. 
\begin{itemize}
\item We call the pair $(M,C)$ a {\it causal symmetric space}. 
Note that $H$ need not be connected, so that the 
requirement $\Ad(H) C = C$ is stronger than the invariance 
under $\Ad(H_e) = \Inn_\g(\fh)$. 
\item We call $(M,C)$ {\it compactly causal} if  $(\g,\tau,C)$ is 
compactly causal. 
\end{itemize}
 
\begin{example}  The following example shows that the causality of $G/H$ can depends
on $\pi_0(H)$. For that let $G =\Ad(\SL_2(\R)) \cong \SO_{2,1}(\R)_e$. Let
$h=\frac{1}{2}\begin{pmatrix} 1 & 0 \\ 0 & -1\end{pmatrix}\in \mathfrak{sl} (2,\R)$. Then
$h$ is an Euler element and the involution $\tau_h$ is given by
\[\tau_h\begin{pmatrix} x & y \\ z & -x\end{pmatrix}=
\begin{pmatrix} x& -y \\ -z & -x  \end{pmatrix}.\]
Hence 
\[\fh = \R h \quad\text{and}\quad \fq = \left\{\begin{pmatrix} 0 & y \\ z & 0\end{pmatrix}\:
y,z\in\R\right\}\cong \R x_+ + \R x_-
\quad \mbox{ with } \quad 
x_+= \begin{pmatrix} 0 & 1\\ 0 & 0\end{pmatrix}, 
x_-= \begin{pmatrix} 0 & 0\\ 1 & 0\end{pmatrix}. \] 
We
also have $\g_{+1}(h)= \R x_+$, $\g_{-1}(h) =\R x_-$ and $\g_0(h) =\R h$.  The map 
$\theta(x) =  -x^\top$ is
a Cartan involution commuting with $\tau_h$. Let
\[ C_+ :=[0,\infty) x_+\subset \g_1 \quad \mbox{ and } \quad C_- := [0,\infty) x_-\subset \g_-.\]
 Then $(\fg,\tau_h,C_+-C_-)$ is
compactly causal and $(\fg,\tau_h,C^c=C_+ +C_-)$ is non-compactly
causal. Let $z=x_+-x_-\in\fk = \so_2(\R)$. Then $\theta =e^{\pi \ad z}\in G^{\tau_h}=:H$
as $\theta $ commutes with $\tau_h$, but $\theta (C^c)=-C^c$. As $H =(H\cap K)H_e=\{\1,\theta\}H_e$ it follows that $C$ is $H$-invariant but 
$\Ad (H)C^c \supeq -C^c$, so that $C^c$ is not. It is only 
invariant under the connected group~$H_e$. 

In general, if $\g$ is simple  hermitian of tube type 
with Euler element $h$ 
and involution $\tau = \tau_h$, then if $C= C_+ - C_-$ defines a compactly
causal structure, then the cone $C^c= C_+ +C_-$ defines a non-compactly causal 
structure on $(\g, \tau_h)$. As above
we see that $C^c$ is invariant under $\Ad(G)^{\tau_h}_e$ but not under
$\Ad (G)^{\tau_h}$. On the other hand $C$ is invariant under the full group $\Ad (G)^{\tau_h}$. 
\end{example}

If $(\g,\tau,C,h)$ is a modular compactly causal symmetric Lie algebra, 
then we further assume that the involution $\tau_h \in \Aut(\g)$ integrates 
to an involution $\tau_h^G$ of $G$ 
satisfying 
\begin{equation}
  \label{eq:tauHinv}
 \tau_h^G(H) = H,
\end{equation}
so that it induces an involution $\tau_h^M$ on~$M$.

\begin{rem} (Integrability of $\tau_h$) \mlabel{rem:2.12}
If $G$ is a connected Lie group with Lie algebra 
$\g$ and $h \in \cE(\g)$ an Euler element, then 
$\tau_h$ does not always integrate to an automorphism of $G$. 

Here is an example where $\g$ is simple. 
Let $G := \SU_{p,p}(\C)$ and $K := {\rm S}(\U_p(\C) \times \U_p(\C))$ 
be the canonical maximal compact subgroup. 
For the Euler element $h := \frac{1}{2}\pmat{ 0 & \1_p \\ \1_p & 0}$, 
the corresponding involution $\tau_h$ acts on $\g$ by 
\[ \tau_h\pmat{a & b \\ c & d} = \pmat{d & c \\ b & a}, \] 
so that we have on the group level  $\tau_h^G(k_1,k_2) = (k_2, k_1)$ on 
$K \subeq \U_p(\C) \times \U_p(\C)$. On 
\[ Z(\tilde G) = \{ (n,z_1, z_2) \in \Z \times C_p \times C_p \: 
e^{\frac{2\pi i n}{p}} z_1 = z_2 \} \cong \Z \times C_p\] 
(see \cite[p.~28]{Ti67}), this leads to the non-trivial involution on 
$\Z \times C_p$, given by 
\[ \tau_h^{\tilde G}(n,z_1) = (-n, \zeta_p^n z_1), 
\quad \mbox{ where } \quad \zeta_p  = e^{\frac{2\pi i}{p}}.\]
It fixes the subgroup $\{0\}\times C_p$ pointwise 
and maps $(-2,\zeta_p)$ to its inverse. 
In particular, the cyclic subgroup $\Gamma = \la (1,1) \ra$ 
is not invariant, so that $\tau_h$ does not integrate to 
the group $\tilde G/\Gamma$. 
\end{rem}

\begin{rem} \mlabel{rem:9.2} 
(Symmetric space structures) 
Let $G$ be a connected Lie group and $H \subeq G$ be a closed subgroup 
with $\L(H) = \fh$, where $(\g,\tau)$ is a symmetric Lie algebra. 
Let $N \subeq H$ be the largest normal subgroup of $G$, contained in $H$, i.e., 
the kernel of the $G$-action on $M$. We put 
\[ G_1 := G/N \quad \mbox{ and } \quad H_1 := H/N.\] 

\nin (a) $M = G/H \cong G_1/H_1$ carries a symmetric space structure 
if and only if $\tau$ integrates to an involution $\tau^{G_1}$ on $G_1$ 
leaving $H_1$ invariant. Then 
$H_1 \subeq G_1^{\tau^{G_1}}$ is an open subgroup because it 
preserves the base point and the corresponding involution 
$\tau^M$ on $M$ is determined by its differential in the base point.  
If $\tau$ integrates to an involution 
$\tau^G$ on $G$, then we need not have $H \subeq G^{\tau^G}$, we 
 only have the weaker condition 
\[ H \subeq \{ g \in G \: g \tau^G(g)^{-1} \in N\}.\]

\nin (b) If $\g$ is semisimple, then $Z(G) \cap H$ acts trivially on $M$, 
hence is contained in $N$. It follows that, for $N = \{e\}$, 
the adjoint action restricts to an injective representation of $H$. 

If, in addition, $\fh$ contains no non-zero ideal of $\g$, 
then $N$ is a discrete 
normal subgroup of the connected group $G$, hence central, and thus 
\begin{equation}
  \label{eq:nhzg}
 N = H \cap Z(G).
\end{equation}

\nin (c) From $\L(H) = \fh$ we obtain $\Ad(H)\fh = \fh$. 
If $\g$ is semisimple, so that $\fq = \fh^\bot$ with respect to the 
Cartan--Killing form implies $\Ad(H)\fq= \fq$. This already implies 
\begin{equation}
  \label{eq:adhcent}
H \subeq \Ad^{-1}(\Ad(G)^\tau) = 
\{ g \in G \: gg^\sharp \in Z(G)\}.
\end{equation}

\nin (d) If $\tau^G(H) = H$ and $\Ad$ is faithful on $H$, 
then $\Ad(H) \subeq \Ad(G)^\tau$ implies $H \subeq G^{\tau^G}$. 
\end{rem}

\section{Fixed points of the modular flow in $M$} 
\mlabel{sec:3new}

In Subsection~\ref{sec:3new.1} we classify for a 
reductive compactly causal symmetric reductive Lie algebra $(\g,\tau)$ 
the orbits of $\Inn_\g(\fh)$ in the set $\cE(\g) \cap \fh$ 
of Euler elements in $\fh$. As this set describes the choices 
of the Euler element $h \in \fh$, this provides a classification 
of all modular structures on the associated symmetric spaces. 
The classification is done by reduction to the simple case.

On $M = G/H$, an element $m = gH$ is fixed by the modular flow 
\begin{equation}\label{eq:alphat}
\alpha_t^M(gH) = \exp(th)gH = g\exp(t\Ad(g)^{-1}h)H 
\end{equation}
if and only if 
$\Ad(g)^{-1}h \in \fh$. Therefore elements of the fixed point set $M^\alpha$ 
correspond to Euler elements in $\cO_h \cap \fh$, where 
$\cO_h = \Ad(G)h$ is the adjoint orbit of $h$. 

In Subsection~\ref{sec:3new.2} we show that the orbits 
of the group $G^h$ in $M^\alpha$ are open and that this 
leads to a natural bijection between the orbit spaces $M^\alpha/G^h$ 
and $(\cO^h \cap \fh)/H$. Both sets correspond to the set of 
double cosets $G^h g H \subeq G$ for which $\Ad(g)^{-1}h \in \fh$. 

\subsection{$\Inn_\g(\fh)$-orbits in $\cE(\g) \cap\fh$}
\mlabel{sec:3new.1}

Let  $(\g,\tau,C)$ be a reductive compactly causal symmetric Lie algebra. 
Then $(\g,\tau)$ decomposes as 
\[ (\g,\tau) 
\cong (\g_0, \tau_0) \oplus \bigoplus_{j = 1}^N (\g_j, \tau_j), \]
where $\g_0$ contains the center and all compact ideals 
and the $\tau$-invariant ideals $\g_j$, $j \geq 1$, are 
either simple hermitian or irreducible of group type. 
This decomposition leads to a corresponding decomposition for $\fh$:
\[ \fh = \fh_0 \oplus \bigoplus_{j = 1}^N \fh_j,\] 
and an element $h \in  \g \setminus \fz(\g)$ 
is an Euler element in $\g$ if and only if 
all its components $h_0, h_1, \ldots, h_N$ either vanish or are Euler elements in $\g_j$. 
 In particular, if $h_j\not= 0$, then $\fg_j$ is simple hermitian 
of tube type (Proposition~\ref{prop:eul-tt}).
Furthermore,  $h_0 \in \fz(\g)$. 

The group 
$\Inn_\g(\fh)$ is a product of the subgroups 
$\Inn_\g(\fh_j)$, $j =0,\ldots, N$, acting on the ideals~$\g_j$. 
This reduces the problem to classify $\Inn_\g(\fh)$-orbits in the 
subset $\cE(\g) \cap \fh$ to the case where 
$(\g,\tau)$ is irreducible. In the group case 
$\g = \fh \oplus \fh$, and $\Inn(\fh)$ acts transitively on 
 $\cE(\g) \cap  \fh = \cE(\fh)$ (Proposition~\ref{prop:cc2}). 
Therefore it suffices 
to  analyze the  case where $\g$ is simple hermitian.

We assume for the moment that $\g$ is simple hermitian and
that $(\fg,\tau, C)$ is compactly causal. We do not assume that
$\cE (\fg) \not=\emptyset$. Let  $\g = \fk \oplus \fp$ a 
Cartan decomposition invariant under $\tau$.
Then the compact causality of $(\g,\tau)$ is equivalent to 
$\fz(\fk) \subeq \fq$. Let 
$\fc \subeq \fh_\fp$ be maximal abelian 
and $\fa \supeq \fc$ maximal abelian in $\fp$, so that 
\[ s := \dim \fc = \rk_\R(\fh) \quad \mbox{ and } \quad 
r := \dim \fa = \rk_\R(\g).\] 

If $h \in \cE(\g) \cap \fh$, then 
$\tau$ preserves the $3$-grading defined by $h$. 
On the euclidean simple Jordan algebra 
$E := \g_1(h)$ (\cite[Rem.~2.24]{Oeh20}), 
$\gamma := -\tau\res_E$ is a Jordan automorphism 
(\cite[Prop.~3.12]{Oeh20}). 
If $\gamma = \id_E$, then $(\g,\tau)$ is of {\it Cayley type} (CT) and 
$\tau = \tau_h$, 
and if this is not the case, then $\gamma$ is non-trivial 
and there are two cases. For split type (ST) 
we have $r = s$, and for non-split type (NST) we have 
$r =  2s$ 
(cf.\ \cite[\S 3]{O91}, \cite[\S D.2.2, Lemma~D.5]{NO21b}, \cite{BH98}).
If $r = s$, then $\Delta(\g,\fc) \cong C_r$ by Moore's Theorem 
as $\g$ is of tube type. 


To understand the $\Inn_\g(\fh)$-orbits in 
$\cE(\g) \cap \fh$, we have to  compare the root systems 
\[ \Delta(\fh,\fc) \subeq \Delta(\g,\fc)\] 
and the inclusions $\cW(\fh,\fc) \into \cW(\g,\fc)$ of the corresponding 
Weyl groups. 
Any Euler element $h \in \cE(\g) \cap \fh$ is conjugate under 
$\Inn_\g(\fh)$ to one in $\fc$. The space 
$\cE(\g)$ consists of a single $\Inn(\g)$-orbit 
(Proposition~\ref{prop:cc2}(a)), but 
$\Inn_\g(\fh)$ may not act transitively on $\cE(\g) \cap \fh$ 
(see Theorem~\ref{thm:classif-H-orb}). 

For the table below, we need the following two examples: 
\begin{ex} \mlabel{ex:twoexs}
(a) Assume that $\tau =\tau_h$  and that
$(\fg,\tau,C)$ is compactly causal of Cayley type.
Then $\fa = \fc$, and the root system $\Delta (\fg,\fa)$ is of type $C_r$:
\[\Delta(\g,\fa) 
\cong C_{r} = \{ \pm 2 \eps_j, \pm \eps_i \pm \eps_j \: 
1 \leq i \not=j \leq r\} \]
and
\[ \Delta (\fh,\fa)\cong A_{r-1}
=\{\pm (\eps_i - \eps_j)\: 1\leq i \not= j\leq r\}.\]
The  Weyl group $\cW (C_r )$ consists of all products of 
permutations of $\{1,\ldots,r\}$ with 
sign changes. Its contains the Weyl group $\cW (A_{r-1})$  as the subgroup of 
permutations.   Let $w_0 =\1$ and $w_j (x_1,\ldots , x_r) = (-x_1,\cdots , -x_j,x_{j+1},
\ldots , x_r)$. Then 
\[\cW (A_{r-1}) \backslash \cW (C_r) /\cW (A_{r-1}) \cong \{w_0,\ldots ,w_r\}.\]
It follows that, up the conjugation by $\Inn_\g(\fh)$, 
the Euler elements are given by  
\[h_k = \frac{1}{2}(\underbrace{1,\ldots ,1}_{k},-1,\ldots , -1) 
\quad \text{for}\ 0\le k\le r.\]
To see this, we first note that $h=(x_1,\ldots ,x_r) \in\fa$ is an Euler element if and only if 
$x_j\in \{ \frac{1}{2}, -\frac{1}{2} \}$. Applying a suitable 
permutation we can always order the $x_j$ so that
$x_1,\ldots , x_k=1/2$ and $x_{k+1},\ldots ,x_r= -1/2$. Note that all 
of those elements are conjugate under $\cW (C_r)$ as this group 
also contains the sign changes. 

\nin 
(b) For $\g = \su_{2s,2s}(\C)$ and $\fh = \fu_{s,s}(\bH)$ we have 
$\Delta(\g,\fa) \cong C_{2s}$ and $\Delta(\fh,\fc) \cong C_{s}$. 
To determine $\Delta(\g,\fc)$, we realize $\g$ with respect to the form 
defined by the hermitian $2\times 2$-block matrix 
\[ B = \pmat{0 & \1_{2s} \\ \1_{2s} & 0}\quad \mbox{ as } \quad 
\g = \{ X \in \fsl_{2r}(\C) \: X^*B = - BX\} \cong \su_{r,r}(\C).\] 
Then $\fa \subeq \g$ can be realized by diagonal matrices of the form 
\[ \diag(x_1, \ldots, x_{2s}, - x_1, \ldots, - x_{2s}),\] 
and $\fc \subeq \fa$ corresponds to diagonal matrices 
satisfying $x_1 =x_2, x_3 = x_4, \ldots, x_{2s-1} = x_{2s}$. 
In these coordinates, 
\[ \Delta(\g,\fa) 
\cong C_{2s} = \{ \pm 2 \eps_j, \pm \eps_i \pm \eps_j \: 
1 \leq i \not=j \leq 2s\}
\quad \mbox{ and }\quad 
\Delta(\g,\fc) \cong C_{s}.\] 
Restricting from $\fa$ to $\fc$, maps 
$\eps_{2j-1}$ and $\eps_{2j}$ (in $C_{2s}$) to $\eps_j$ (in $C_s$) for $j = 1,\ldots, s$.
Now it is easy to see that 
\[  \Delta(\g,\fc) = \Delta(\fh,\fc) \cup \{0\}.\]

\nin (c) For $\g = \sp_{4s}(\R)$ and $\fh = \sp_{2s}(\C)$ we have exactly 
the same pattern and therefore 
\[  \Delta(\g,\fc) = \Delta(\fh,\fc) \cup \{0\}.\]
\end{ex}

\hspace{-15mm}
\begin{tabular}{||l|l|l|l|l|l|l||}\hline
$\g$ &  $\g^c = \fh + i \fq$  & $\fh$ 
&\mbox{ type } & $\Delta(\g,\fa)$  & $\Delta(\g,\fc)$ & 
$\Delta(\fh,\fc)$  \\ 
\hline
\hline 
$\su_{r,r}(\C)^{\oplus 2}$ & $\fsl_{2r}(\C)$  \phantom{\Big(}& $\su_{r,r}(\C)$&\mbox{(GT)} 
& $C_r \oplus C_r$ & $C_r$ & $C_r$   \\
 $\sp_{2r}(\R)^{\oplus 2}$ &  $\sp_{2r}(\C)$  &$\sp_{2r}(\R)$ &\mbox{(GT)} 
& $C_r \oplus C_r$ & $C_r$ & $C_r$     \\
 $\so_{2,d}(\R)^{\oplus 2}, d > 2$ & $\so_{2+d}(\C)$ & $\so_{2,d}(\R)$ & \mbox{(GT)} 
&  $C_2 \oplus C_2$ &$C_2$ & $C_2$    \\
 $\so^*(4r)^{\oplus 2}$ & $\so_{4r}(\C)$ & $\so^*(4r)$& \mbox{(GT)}& 
$C_r \oplus C_r$ & $C_r$ & $C_r$      \\
 $(\fe_{7(-25)})^{\oplus 2}$ & $\fe_7(\C)$ & $\fe_{7(-25)}$ &\mbox{(GT)}&  
$C_3 \oplus C_3$ & $C_3$ & $C_3$     \\
\hline
 $\su_{r,r}(\C)$ & $\su_{r,r}(\C)$ & $\R \oplus \fsl_r(\C)$&\mbox{(CT)} & $C_r$ & $C_r$ & $A_{r-1}$ \\
 $\sp_{2r}(\R)$ &  $\sp_{2r}(\R)$  &$\R \oplus \fsl_r(\R)$ &\mbox{(CT)} & $C_r$ & $C_r$ & $A_{r-1}$    \\
 $\so_{2,d}(\R), d > 2$ & $\so_{2,d}(\R)$ & $\R \oplus \so_{1,d-1}(\R)$ 
& \mbox{(CT)} &  $C_2$ &$C_2$ & $A_1$  \\
 $\so^*(4r)$ & $\so^*(4r)$ & $\R \oplus \fsl_r(\H)$& \mbox{(CT)}& $C_r$ & $C_r$ & $A_{r-1}$   \\
 $\fe_{7(-25)}$ & $\fe_{7(-25)}$ & 
$\R \oplus \fe_{6(-26)}$ &\mbox{(CT)}&  $C_3$ & $C_3$ & $A_2$   \\
\hline
$\su_{r,r}(\C)$ & $\fsl_{2r}(\R)$ & $\so_{r,r}(\R)$ & \mbox{(ST)}& $C_r$ & $C_r$ & $D_r$  \\
$\so^*(4r)$ & $\so_{2r,2r}(\R)$ & $\so_{2r}(\C)$ & \mbox{(ST)}& $C_r$ & $C_r$ & $D_r$  \\
$\fe_{7(-25)}$ & $\fe_7(\R)$ & $\fsl_4(\H)$ & \mbox{(ST)}& $C_3$ & $C_3$ & $A_3=D_3$ 
  \\
\hline
$\so_{2,d}(\R)$& $\so_{p+1,q+1}(\R)$ & $\so_{1,p}(\R) \oplus \so_{1,q}(\R)$
& \mbox{(ST)} & $C_2$ & $C_2$ &$A_1 \oplus A_1 = D_2$    \\
$d =p+q> 2$& $1 < q <d-1$ & & &  & &  \\
\hline
$\su_{2s,2s}(\C)$ & $\fsl_{2s}(\H)$ &$\fu_{s,s}(\H)$ & \mbox{(NST)} 
& $C_{2s}$ & $C_s$ & $C_s$   \\
$\sp_{4s}(\R)$ & $\fu_{s,s}(\H)$ & $\sp_{2s}(\C)$ & \mbox{(NST)}
& $C_{2s}$ & $C_s$ & $C_s$  \\
$\so_{2,d}(\R)$ & $\so_{1,d+1}(\R)$  & $\so_{1,d}(\R)$             & \mbox{(NST)}&  & $A_1$ & $A_1$ \\
\hline
\end{tabular} \\[2mm] 
{\rm Table 1: Irreducible compactly causal symmetric Lie algebras 
$(\g,\tau)$ with $\cE(\g) \cap \fh \not=\eset$} \\ 

From Table 1 we get the following classification theorem. 
Here we use that, as in \cite[Thm.~3.10]{MN21}, orbits of 
Euler elements can  be classified by 
representatives in a positive Weyl chamber
in terms of a root basis. They correspond to $3$-gradings 
of the root system $\Delta(\g,\fa)$ for $h \in \fa \cap \cE(\g)$ and to 
$3$-gradings of $\Delta(\fg,\fc)$ for $h \in \fc \cap \cE(\g)$. 
We also recall from Proposition~\ref{prop:cc2} that 
the existence of an Euler element in a simple hermitian Lie algebra 
$\g$ implies that $\g$ is of tube type.

\begin{thm} \mlabel{thm:classif-H-orb}
{\rm(Classification of $\Inn_\g(\fh)$-orbits in $\cE(\g) \cap \fh$)} 
Let $(\g,\tau)$ be an irreducible compactly causal symmetric Lie algebra. 
Then we have the following situations: 
\begin{itemize}
\item[\rm(GT)] For  group type 
$\Delta(\g,\fc) = \Delta(\fh,\fc) \cong C_r$, and 
$\Inn_\g(\fh)$ is transitive on $\cE(\g) \cap \fh$. 
\item[\rm(CT)] For Cayley type 
\[ \Delta(\fh,\fc) \cong A_{r-1} \subeq C_r = \Delta(\g,\fc).\]  
In the canonical identification of $\fc$ with $\R^r$, we have 
\[ \cE(\g) \cap \fc = 
\{ \shalf(\pm 1, \cdots, \pm 1) \} \] 
and the orbits of $\Inn_\g(\fh)$  
in $\cE(\g) \cap \fh$ are represented by the elements 
\begin{equation}
  \label{eq:hk-class}
 h_k = \shalf(\underbrace{1,\ldots, 1}_{k}, -1, \ldots, -1), \quad 
k = 0,1, \ldots, r. 
\end{equation}
\item[\rm(ST)] For split type 
\[ \Delta(\fh,\fc) \cong D_r \subeq C_r = \Delta(\g,\fc), \]  
where we identify $D_3 \cong A_3$ and $D_2 \cong A_1 \oplus A_1$.
Then $\Inn(\g)^\tau$ acts transitively on 
$\cE(\g) \cap \fh$ and its identity component 
$\Inn_\g(\fh)$ has two orbits represented by $h_{r-1}$ and $h_r$. 
\item[\rm(NST)] For non-split type 
\[ \Delta(\fh,\fc) = \Delta(\g,\fc) = C_r\]  
and $\Inn_\g(\fh)$ acts transitively on $\cE(\g) \cap \fh$.
\end{itemize}
\end{thm}

\begin{prf} (GT) Here $(\g,\tau) \cong (\fh \oplus \fh,\tau_{\rm flip})$, where 
$\fh$ is simple hermitian of tube type. 
Accordingly $\fa = \fc \oplus \fc$, and the assertion follows. 

\nin (CT) 
The first assertion follows from Table~1 and the remainder 
from the discussion in Example~\ref{ex:twoexs}(a).

\nin (ST) For split type a look at Table 1 shows that 
\[ \Delta(\fh,\fc) \cong D_r \subeq C_r = \Delta(\g,\fc), \]  
where we identify $D_3 \cong A_3$ and $D_2 \cong A_1 \oplus A_1$.
Then $\cW(D_r) \trile \cW(C_r)$ is an index $2$ subgroup 
and $\cW(C_r)$ acts by automorphisms on $D_r$. Therefore 
$\Inn(\g)^\tau$ induces on $\fc$ the full group $\cW(C_r)$. 
Hence the  $\cW(D_r)$-orbits in $\cE(\g) \cap \fc$ are represented by 
$h_{r-1}$ and $h_r$, and both are conjugate under $\cW(C_r)$. 

\nin (NST) Table 1 shows that, for non-split type, we have  
$\Delta(\fh,\fc) = \Delta(\g,\fc) = C_r.$ 
As $\cW(C_r)$ acts transitively on $\cE(\g) \cap \fc$, 
it follows that $\Inn_\g(\fh)$ acts transitively on $\cE(\g) \cap \fh$.
\end{prf}

\begin{cor} \mlabel{cor:conjofhs}
If $(\g,\tau)$ is irreducible and not of Cayley type, 
then $\Inn(\g)^\tau$ acts transitively on~$\cE(\g) \cap \fh$. 
\end{cor}

\subsection{Connected components  of $M^\alpha$} 
\mlabel{sec:3new.2}

Recall our global context from Subsection~\ref{subsec:2.3}.  Our goal in
this section is to prove the following proposition:

\begin{prop} \mlabel{prop:5.3} 
The orbits of the identity component $G^h_e$ on the 
fixed point set 
\[ M^\alpha := \{ gH \in M\: \Ad(g)^{-1}h \in \fh \} 
= \{ m \in M \: X_h^M(m) = 0 \} \]  
of the modular flow coincide with its connected components. 
\end{prop}

Example~\ref{ex:ct} below shows  
that, for Cayley type involutions on hermitian 
simple Lie algebras, the connected components of $M^\alpha$ may have 
different dimensions. Our discussion of de Sitter space 
in Subsection~\ref{subsec:subgrpdata} reveals a situation where 
$G^h$ does not preserve the connected components of $M^\alpha$. 

\begin{prf} Clearly, $M^\alpha$ is invariant 
under the action of the group $G^h$ which commutes with $\alpha$. 
Therefore $G^h_e$ preserves all connected components of $M^\alpha$. 
We now show that all $G^h_e$-orbits in $M^\alpha$ are open, hence coincide 
with its connected components. 

First we consider the base point $eH$ which is contained in $M^\alpha$ 
because  $h \in \fh$. 
Using the exponential function $\Exp_{eH} \: \fq \to M$ as a chart around $eH$, 
its equivariance with respect to $\alpha$ and the one-parameter group 
$e^{t \ad h}$ on $\fq$, it follows that $M^\alpha$ is a symmetric subspace of $M$ 
with $T_{eH}(M^\alpha) = \fq_0(h) \subeq \g^h$. Therefore 
$G^h_e.eH \supeq \Exp_{eH}(\fq_0(h))$ contains a neighborhood of $eH$ in $M^\alpha$, 
and this implies that the $G^h_e$-orbit of $eH$ in $M^\alpha$ is open. 

Now let $m = gH \in M^\alpha$ be arbitrary. 
Then 
\[ \alpha_t^M(gH)=\exp(th)gH = g \exp(t \Ad(g)^{-1}h) H \] 
shows that 
$m \in M^\alpha$ is equivalent to $\Ad(g)^{-1}h \in \fh$ and 
\[ \Phi \: \fq \to M, \quad \Phi(x) := g \Exp_{eH}(x) \] 
defines a chart around $m$ satisfying 
\begin{equation}
  \label{eq:alphacov}
   \alpha_t \Phi(x) = \exp(th) g \Exp_{eH}(x) = \Phi(e^{t \Ad(g)^{-1}h}x).
\end{equation}
Hence $\Phi$ maps 
\begin{equation}
  \label{eq:5.1}
 \{x \in \fq \: [\Ad(g)^{-1}h,x] = 0\} 
= \fq \cap \Ad(g)^{-1}\g^h 
\end{equation}
onto a neighborhood of $m$ in $M^\alpha$. Finally 
\[ \Phi(\fq \cap \Ad(g)^{-1}\g^h) 
\subeq g \exp(\fq \cap \Ad(g)^{-1} \g^h) H
\subeq \exp(\g^h) g H \subeq G^h_e.m\] 
shows that $G^h_e.m$ is open in $M^\alpha$. 
\end{prf} 

\begin{ex} (The group case) In the group case the set $M^\alpha = G^\alpha$ 
is rather simple to describe. 
Here $M = G$, $h = (h_0, h_0)$ and 
$\Ad(g_1, g_2)h = (\Ad(g_1)h_0, \Ad(g_2)h_0) \in\fh$ is equivalent to 
$\Ad(g_1)h_0 = \Ad(g_2)h_0$, i.e., $g_1^{-1} g_2 \in G^{h_0}$. 
If this is the case, then $\Ad(g_1, g_2)h = \Ad(g_1, g_1)h \in \Ad(\Delta_G)h$. 

We also note that $(G \times G)^h = G^{h_0} \times G^{h_0}$ 
acts transitively on the submanifold $M^\alpha = G^{h_0}$. 
\end{ex}

\begin{ex} (An example of Cayley type) \mlabel{ex:ct} 
As in Example~\ref{ex:twoexs}, we 
consider the hermitian form on $\C^{2r}$ defined by the block 
diagonal matrix $B = \pmat{ 0 & \1 \\ \1 & 0}$, so that 
\[ \g := \{ X \in \fsl_{2r}(\C) \: X^* B = - B X \} \cong \su_{r,r}(\C) 
\quad \mbox{ and } \quad 
\tau\pmat{a & b \\ c & d} = \pmat{a & -b \\ -c & d}.\] 
Then 
\[ \fh = \Big\{ \pmat{ X & 0 \\ 0 & - X^*} \: \tr(X) \in \R \Big\} 
\cong \R \oplus \fsl_r(\C). \] 
The corresponding groups are 
\[ G \cong \SU_{r,r}(\C) \supeq H \cong \{ g \in \GL_r(\C) \: \det(g) \in \R \}.\] 
Then 
\[ \fa = \{ \diag(x_1, \ldots, x_r, - x_1, \ldots, -x_r \} \: x_j \in \R\}
\cong \R^r,\] 
and in these coordinates, we obtain Euler elements 
\[ h_k = \frac{1}{2}(\1_k,-\1_{r-k}, -\1_k, \1_{r-k} ), \quad 
k = 0,\ldots, r\] 
(Theorem~\ref{thm:classif-H-orb}). 
Then 
\[ \dim \cO^H_{h_k} = \dim [h_k,\fh] 
= \dim_\R(M_{k,r-k}(\C) \oplus M_{r-k,k}(\C)) = 2k (r-k) \] 
depends on $k$. For $k = 0,r$, the orbit is trivial. 
This corresponds to the fact that the base point $eH$ in 
$M = G/H \cong \cO^G_h$ is an isolated fixed point of the modular 
flow. The components 
\[ M_k^\alpha := \{ gH \: \Ad(g)^{-1}h \in \cO^H_{h_k}\}, \quad k = 0,\ldots, r,\] 
are of different dimensions (for $r > 1$). 
\end{ex}

\begin{lem} We have a bijection 
\[ \Gamma \: M^\alpha/G^h \to (\cO_h\cap \fh)/H, \quad 
G^h g H \mapsto \Ad(H) \Ad(g)^{-1} h.\] 
 \end{lem}

 \begin{prf} As $gH \in M^\alpha$ is equivalent to $\Ad(g)^{-1}h \in \fh$, 
we obtain a well-defined surjective map sending the $G^h$-orbit $G^h gH$ to the 
$H$-orbit $\Ad(H) \Ad(g)^{-1}h \subeq \cE(\g) \cap \fh$. 
   
If $g_1 H, g_2 H \in M^\alpha$ map to the same $H$-orbit in $\cE(\g) \cap \fh$, 
then there exists an element $k \in H$ with 
$\Ad(g_1)^{-1}h = \Ad(k^{-1}) \Ad(g_2)^{-1}h$, i.e., 
$g_2 k g_1^{-1} \in G^h$. This implies that 
$g_2 H = g_2 k H \in G^h g_1 H.$ 
As the surjectivity of $\Gamma$ follows from the fact that 
$gH \in M^\alpha$ is equivalent to $\Ad(g)^{-1}h \in \fh$, 
$\Gamma$ is bijective. 
 \end{prf}

 \begin{rem} (a) If $H$ is connected, then the preceding lemma shows that the 
orbits of $\Ad(H) = \Inn_\g(\fh)$ in $\cE(\g) \cap \fh$ 
are in one-to-one correspondence with the 
$G^h$-orbits in $M^\alpha$ for $M = G/H$. 
Therefore Theorem~\ref{thm:classif-H-orb}
provides a classification of $G^h$-orbits in $M^\alpha$ whenever 
$H$ is connected. 

\nin (b) As $G^h$-orbits in $M^\alpha$ are open by 
Proposition~\ref{prop:5.3},  we may also think of 
$\pi_0(G^h)$-orbits in $\pi_0(M^\alpha)$ being classified by 
$(\cO_h \cap \fh)/H$. 
 \end{rem}

\begin{rem}
The involution $\tau_h^M$ acts on the fixed point manifold $M^\alpha$, but not 
necessarily trivially. A typical example is the group type space 
$M = \tilde\SL_2(\R)$, 
where $Z(G) \subeq M^\alpha$ and $Z(G) \cap M^{\tau_h^G} = \{e\}$. 
\end{rem}

\begin{rem} \mlabel{rem:7.6new}
($G^h$ and $G^{\tau^G}$ for simple Lie algebras) 
As $G^h = \{ g \in G \: \Ad(g)h = h\}$, we have 
\[ \Ad(G^h) = \Ad(G)^h = \Ad(G)^{\ad h} 
\subeq \Ad(G)^{\tau_h}.\] 
We claim that, if $\g$ is simple, then 
 $\Ad(G^h)$ is of index $2$ in $\Ad(G)^{\tau_h}$. 
In fact, in this case $\g^{-\tau_h}$ is a direct sum of two irreducible
$\g_0(h)$-modules $\g_{\pm 1}(h)$ (\cite[Cor.~1.3.13]{HO97}), and any element 
$\Ad(g)$ commuting with $\tau_h$ either 
preserves both subspaces or exchanges them. If it preserves both, 
then $\Ad(g)$ commutes with $\ad h$, i.e., $g \in G^h$. 
If not, then $\Ad(g)h = - h$, and as such an element 
exists (Proposition~\ref{prop:cc2}(a)), 
the index of $\Ad(G^h)$ in $\Ad(G)^{\tau^G}$ is 
two. 
\end{rem}

\section{Wedge domains in compactly causal  symmetric spaces} 
\mlabel{sec:3}

In this section we introduce the wedge domains 
$W_M^+(h)$, $W_M^{\rm KMS}(h)$ and $W_M(h)$ 
in a compactly causal symmetric 
space~$M =  G/H$ with the infinitesimal data $(\g,\tau,C,h)$ 
(under the global assumption from 
Subsection~\ref{subsec:2.3} that the universal complexification 
$\eta_G$ has discrete kernel). 
One of our main results asserts that 
these three subsets coincide, which is far from obvious from 
their definitions. Moreover, their connected components are 
orbits of an open real Olshanski subsemigroup of~$G$. 
We therefore start recalling the construction of Olshanski subsemigroups.

\subsection{Tube domains of compactly causal symmetric spaces} 

\begin{defn} \mlabel{def:cplxols} 
(Complex Olshanski semigroups) 
Let $G$ be a connected Lie group with 
simply connected covering group $q_G \: \tilde G \to G$. 
For a pointed 
closed convex $\Ad(G)$-invariant cone ${C_\g \subeq \g}$, 
Lawson's Theorem (\cite{La94},\cite[Thm.~IX.1.10]{Ne00}) implies the existence 
of a semigroup $\Gamma_{\tilde G}(C_\g)$ which is a covering of the subsemigroup 
$\eta_{\tilde G}(\tilde G) \exp(i C_\g)$ of the universal complexification 
$\tilde G_\C$ (the simply connected group with Lie algebra $\g_\C$). 
Then the exponential function $\exp \: \g + i C_\g \to \tilde G_\C$ 
lifts to an exponential function $\Exp \: \g + i C_\g \to \Gamma_{\tilde G}(C_\g)$,  
and the polar map 
\[  \tilde G \times C_\g \to \Gamma_{\tilde G}(C_\g), \quad 
(g,x) \mapsto g \Exp(ix) \] 
is a homeomorphism. 
We now define the {\it closed complex Olshanski semigroup} 
corresponding to the pair $(G,C_\g)$ by 
\begin{equation}
  \label{eq:3.1a}
 \Gamma_G(C_\g) := \Gamma_{\tilde G}(C_\g)/(\ker q_G),
\end{equation}
where $q_G \: \tilde G \to G$ is the universal covering map 
(cf.~\cite[Thm.~XI.1.12]{Ne00}). 
Then the polar map $G \times C_\g \to \Gamma_G(C_\g)$ is a 
homeomorphism, and if 
$C_\g$ has interior points, then it restricts to a 
diffeomorphism  from 
$G \times C_\g^0$ onto the open subsemigroup 
\begin{equation}
  \label{eq:3.1}
\Gamma_G(C_\g^0) := \Gamma_{\tilde G}(C_\g^0)/(\ker q_G) 
\subeq \Gamma_G(C_\g),
\end{equation}
\end{defn}

Complex Olshanski semigroups are non-abelian generalizations 
of complex tube domains defined by open cones in real vector spaces; 
the domains $\g + i C_\g^0$ which are the tangent objects of 
the open complex Olshanski semigroups $\Gamma_G(C_\g^0)$ 
are typical examples. We now turn to similar 
objects for symmetric spaces, which on the tangent level correspond to the 
tube domain $\fq + i C^0 \subeq \fq_\C$.

\begin{defn} (The complex tube domain of a compactly causal symmetric space) 
Let $(G,\tau,H, C)$ be a {\it compactly causal symmetric Lie group}, i.e., 
$G$ is a connected Lie group, $\tau^G$ an involutive automorphism of~$G$, 
$H \subeq G^{\tau^G}$ an open subgroup, and $C \subeq \fq$ a pointed 
generating $\Ad(H)$-invariant closed convex cone such that 
$C^0$ consists of elliptic elements. 
We assume that the universal complexification $\eta_G \: G \to G_\C$ has 
discrete kernel (which is the case if $G$ is simply connected or semisimple). 
The involution $\tau^G \in \Aut(G)$ induces by the universal property 
a holomorphic automorphism $\tau^{G_\C}$ of $G_\C$ with 
$\eta_G(H) \subeq G_\C^{\tau^{G_\C}}$, 
and 
\[ G_\C^\sharp = \{ g \in G_\C \: \tau^{G_\C}(g) = g^{-1} \} \] 
is a complex symmetric subspace of $G_\C$. 

We call the fiber product 
\[ \cT_M(C) := G \times_H i C^0 \] 
the {\it tube domain of the pair $(M,C)$}. 
To obtain a complex manifold structure on this domain, 
we observe that the quadratic representation 
$Q \: M = G/H\to G$ extends to a  map 
\[ Q_\cT \: \cT_M(C) = G \times_H i C^0 \to G_\C^\sharp,  \quad 
[g, ix] \mapsto g.\Exp(ix) = g \exp(2ix) g^\sharp\] 
which is a covering of an open subset. 
Therefore $\cT_M(C)$ carries a unique complex manifold 
structure for which $q_\cT$ is holomorphic.
Clearly, $\cT_M(C) \subeq G \times_H i C$ contains $M$ as 
$G \times_H \{0\}$ in its ``boundary''  (cf.\ \cite{HOO91},\cite{KNO97}, \cite{Ne99}). 
\end{defn}

\begin{rem} \mlabel{rem:exatube} 
(a) Suppose that $(\g,\tau,C)$ is extendable 
and that $C_\g \subeq \g$ is a $-\tau$-invariant $\Ad(G)$-invariant 
pointed cone 
with $C = C_\fg \cap \fq$ (see Subsection~\ref{subsec:2.1}). 
Let $G_C \trile G$ be the normal integral subgroup 
with Lie algebra $\g_C = C_\g - C_\g$. 
Then we have the corresponding open 
complex Olshanski semigroup 
\[ \cT_{G_C}(C_\g) := G_C \exp(i C_\g^0) \] 
and the quadratic representation 
$Q \: M \to G_C$ extends to a map 
\[ Q \: \cT_M(C) = G_C \times_H i C^0 \to \cT_{G_C}(C_\g), \quad 
(g, ix) \mapsto g \exp(2ix) g^\sharp.\] 
By Lemma~\ref{lem:real-olsh} and Remark~\ref{rem:a.3}, its range 
is the identity component of the complex submanifold 
$\cT_{G_C}(C_\g)^\sharp$: 
\[ Q(\cT_M(C)) = \cT_{G_C}(C_\g)_e^\sharp 
= \{ ss^\sharp \: s \in \cT_{G_C}(C_\g)\}. \] 
Since $Q \: \cT_M(C) \to Q(\cT_M(C))$ is a covering map 
and the $G_C$-action on its range lifts to the $G_C$-action on $M$, 
the holomorphic action of the semigroup $\cT_{G_C}(C_\g)$ 
on $Q(\cT_M(C))$ by $s.m = sms^\sharp$ lifts to a holomorphic action 
\begin{equation}
  \label{eq:holsemact}
\cT_{G_C}(C_\g) \times \cT_M(C) \to \cT_M(C),   
\end{equation}
extending the $G$-action.  

\nin (b) Let $\tau^{G_\C} \in\Aut(G_\C)$ be the holomorphic involution 
with $\tau^{G_\C} \circ \eta_G = \eta_G \circ \tau^G$ 
and $H_\C \subeq G_\C^{\tau^{G_\C}}$ be an open subgroup with 
$\eta_G^{-1}(H_\C) = H$. Then we have a natural embedding 
of $G/H$ as the $G$-orbit of the base point $eH_\C \in M_\C := G_\C/H_\C$, 
and 
\begin{equation}
  \label{eq:tauM}
 \cT_M(C) = G.\Exp(iC^0) 
\end{equation}
can be identified with the orbit 
$\Gamma_G(C_G^0).eH_\C$ of the base point under the action of the 
open complex Olshanski semigroup $\cT_G(C_G)$, which is an open subset 
of~$M_\C$ (cf.~\cite[Lemma~1.3]{HOO91}). 
\end{rem}

\begin{ex} (a) Let $(G, \tau)$ be a connected symmetric Lie group 
and $(\g, \tau, C)$ be a corresponding compactly causal 
symmetric Lie algebra. 
Suppose that 
 $C_\g \subeq \g$ is a pointed $\Ad(G)$-invariant 
generating cone satisfying 
$-\tau(C_\g) = C_\g$ and $C := C_\g \cap \fq$ 
(cf.\ Theorem~\ref{thm:extend}). 
Then we have a corresponding open complex Olshanski semigroup 
$S = \cT_G(C_\fg)$ (Definition~\ref{def:cplxols}) and thus an embedding 
\[ \cT_{G/G^{\tau^G}}(C) = G \times_{G^{\tau^G}} i C^0 \cong  \bigcup_{g \in G} 
g \exp(i C^0) g^\sharp \ {\buildrel ! \over =}\  S^\sharp_e \into  S \] 
(Lemma~\ref{lem:real-olsh}(4)). 

\nin (b) Consider a symmetric space $M = (G,\bullet)$ of group 
type and an $\Ad(G)$-invariant cone $C_\fg \subeq \fg$. 
The corresponding 
symmetric Lie group is $G \times G$ with 
\[ \tau^{G \times G}(g_1, g_2) = (g_2, g_1), \quad 
(G \times G)^{\tau^{G \times G}} = \Delta_G, \quad  
\fq = \{ (x,-x) \: x \in \g\},\] and the cone 
\[ C = \{ (x,-x) \: x \in C_\fg\} = (C_\g \oplus - C_\g) \cap \fq\] 
is $\Ad(\Delta_G)$-invariant. 
We have the open complex Olshanski semigroup 
\[ \Gamma_{G \times G}(C_\g \oplus - C_\g) 
= \Gamma_G(C_\g) \times \Gamma_G(-C_\fg).\] 
Here $(g_1, g_2)^\sharp = (g_2^{-1}, g_1^{-1})$ is the 
$\sharp$-operation in $G \times G$, so that the 
$\sharp$-fixed points in the open complex Olshanski semigroup 
\[S := \Gamma_{G \times G}(C_\fg^0 \oplus -C_\fg^0) \] 
are the pairs $(s, s^{-1})$, $s \in \Gamma_G(C_\fg^0)$. We thus obtain 
\[ \cT_G(C) 
= (G \times G) \times_{\Delta_G} i C^0
\cong G \times i C^0 \cong \Gamma_G(C_\fg^0)
= \{ (s,s^{-1}) \: s \in \Gamma_G(C_\fg^0) \} \cong S^\sharp.\] 
As a complex manifold, this is a copy of the open complex Olshanski semigroup 
$\Gamma_G(C_\fg^0)$. 
\end{ex}

\subsection{The modular flow and three types of wedge domains} 

We are now ready to introduce the three different type of wedge domains 
in $M$.  We already introduced the modular flow in \eqref{eq:alphat} but
repeat the definition here:  
\begin{defn} \mlabel{def:modflow} (The modular flow) The Euler element 
$h \in \fh$ 
defines an $\R$-action by automorphisms on $G$ via 
\[ \alpha_t(g) = \exp(th) g \exp(-th), \qquad g \in G.\] 
Then $\alpha$ preserves all connected components of 
the subgroup $G^{\tau^G}$, hence in particular $H$. 
Therefore $\alpha_t$ induces a flow 
\begin{equation}
  \label{eq:modflow-M}
 \alpha_t^M(gH) = \exp(th) gH = g\exp(t \Ad(g)^{-1}h) H 
\end{equation}
on $M = G/H$. This flow is generated by the {\it modular vector field} 
\begin{equation}\label{eq:XM}
 X_h^M \in \cV(M), \quad 
X_h^M(m) = \frac{d}{dt}\Big|_{t=0}\alpha_t(m).
\end{equation}

The modular flow extends to $\R$-actions by holomorphic maps 
on the complex tube domains $\cT_G(C_\g)$ and on $\cT_M(C)$ via 
\[ \alpha_t(g\Exp(ix)) = \alpha_t(g) \Exp(ie^{t \ad h}x) 
\quad \mbox{ and } \quad 
\alpha_t^M([g, ix]) :=  [\alpha_t(g), i e^{t \ad h}x].\] 
Their infinitesimal generators are denoted $X^{\cT_G}_h$ and 
$X^{\cT_M}_h$, respectively. 

On $G_\C$ we even obtain a holomorphic flow 
by $\alpha_z(g) = \exp(zh) g \exp(-zh)$, $z \in \C$, 
but on $\cT_G(C_\g)$ and $\cT_M(C)$, the real flow does not extend to 
all of $\C$. However, for 
$s$ in the closed complex semigroup $\Gamma_G(C_\g)$, we consider the orbit map 
\[ \alpha^s \: \R \to \Gamma_G(C_\g), \quad 
t \mapsto \alpha_t(s) \] 
and define 
\[ \alpha_{x + i y}(s) := \alpha^s(x + i y) \quad \mbox{  for  } \quad 
x + i y \in \C,\] whenever the maximal local flow of the vector field 
$i X^{\cT_M}$ with initial value $\alpha_x(s)|_{s=0}$ is defined in $y \in \R$. 
This implies that 
$\alpha^s$ extends to a continuous map on the closed strip 
between $\R$ and $\R + i y$ which is holomorphic on the interior 
in the sense that its composition with the natural map 
$\Gamma_G(C_\g) \to G_\C, g \Exp(ix) \mapsto \eta_G(g)\exp(ix)$  
is a holomorphic $G_\C$-valued map (cf.\ Definition~\ref{def:cplxols}).
\end{defn}
 
Recall that {\it we assume that $\tau$ integrates to an involution 
$\tau_h^G$ on $G$ which leaves $H$ invariant,} so that 
it induces an involution $\tau_h^M$ on $M = G/H$ 
(see the end of Subsection~\ref{subsec:2.3}). 
\begin{defn} \mlabel{def:wedgedom} 
We consider the following types of {\it wedge domains}:
\footnote{In 
complex analysis ``domains'' are assumed connected. Here we use 
the term ``domain'' for an open subset.} 
\begin{itemize}
\item The {\it positivity domain of the modular vector field $X_h^M$ in $M$} is  
\[ W_M^+(h) := \{ m \in M \: X_h^M(m) \in V_+(m) \},\] 
where $V_+(m) \subeq T_m(M)$ is the open cone corresponding to the 
$G$-invariant cone field with $V_+(eH) = C^0$ in $T_{eH}(M) \cong \fq$. 
If $M$ is Lorentzian and $V_+(m)$ is the future light cone in $T_m(M)$, 
then this is the domain where the modular vector field is future oriented 
timelike.  
\item The {\it  KMS wedge domain} is 
  \begin{align*}
W_M^{\rm KMS}(h) 
:=& \{ m \in M \: (\forall z \in \cS_\pi)\ \alpha_z(m) \in \cT_M(C)\}\\
= &\{ m \in M \: (\forall y \in (0,\pi))\ \alpha_{iy}(m) \in \cT_M(C)\},
  \end{align*}
where $\alpha_z(m)$ is assumed to be defined in $m$ 
in as in Definition~\ref{def:modflow},
 and 
\[ \cS_\pi = \{ z \in \C \: 0 < \Im z < \pi\}.\]
This is the set of all points $m \in M$ whose $\alpha$-orbit map extends 
analytically  to a map $\cS_\pi \to \cT_M(C)$. Comparing boundary values on 
$\R$ and $\pi i + \R$ resembles KMS conditions; hence the name 
(see also \cite[App.~A.2]{NOO21}). 
\item Let $m \in M^\alpha$ be a fixed point of the modular flow 
and $C_m = \oline{V_+(m)} \subeq T_m(M)$. Then $\alpha$ induces a 
$1$-parameter group of linear automorphisms on $T_m(M)$ whose infinitesimal 
generator has the eigenvalues $-1,0,1$. We write 
$T_m(M)_j$, $j =-1,0,1$, for the corresponding eigenspaces and consider the pointed cone 
\[ C_m^c := C_m \cap T_m(M)_1 - C_m \cap T_m(M)_{-1} \subeq T_m(M).\] 
Note that 
$C_m^c - C_m^c \subeq T_m(M)$ is a 
is a vector space complement of $T_m(M^\alpha) = T_m(M)_0$. 

Then the {\it wedge domain in $m$} is defined as 
\[ W_M(h)_m 
:= G^h_e.\Exp_m(C_m^{c,0}),\] 
where $C_m^{c,0} := (C_m^c)^0$ 
is the relative interior of $C_m^c$ in  $T_m(M)_1 + T_m(M)_{-1}$ 
(cf.\ \eqref{eq:cpm}) and 
$G^h_e = (G^h)_e$ is the connected subgroup 
with Lie  algebra $\g^{\tau_h} = \g_0(h)$. 
On the infinitesimal level, this domain corresponds to the wedge 
$T_m(M)_0 + C_m^c \subeq T_m(M)$. 

Identifying $C$ with $C_{eH}\subeq T_{eH}(M)\cong \fq$, 
the {\it wedge domain in the base point is} 
\[ W_M(H)_{eH} := G^h_e.\Exp_{eH}( (C_+ + C_-)^0),
\quad \mbox{ where } \quad C_\pm = \pm C \cap \fq_{\pm 1}(h).\]
 We define the {\it polar wedge domain of $(M,h)$} as 
\begin{equation}
  \label{eq:polwedgedom-M}
W_M(h) 
:= \bigcup_{m \in M^\alpha} W_M(h)_m = \bigcup_{m \in M^\alpha} \Exp_m(C_m^{c,0}).
\end{equation}
\end{itemize}
\end{defn}

\begin{rem} \mlabel{rem:6.5} 
For $g \in G$ and $v \in T_m(M)$, the fact that $G$ acts by automorphisms 
on the symmetric space $M$ implies that 
\begin{equation}
  \label{eq:Exp-covar}
 g.\Exp_m(v) = \Exp_{gm}(g.v) \quad \mbox{ for }\quad g \in G, m \in M, v 
\in T_m(M).
\end{equation}
Now let $m = gH \in M^\alpha$ and identify $T_{eH}(M)$ with $\fq$. 
For $h' := \Ad(g)^{-1}h \in \fh$, we then have 
\[ C_m = g.C \quad \mbox{ and }  \quad 
T(\alpha_t)(g.x) = g.(e^{t \ad h'}x) \quad \mbox{ for }\quad 
t \in \R,x \in \fq,\]
so that 
\begin{equation}
  \label{eq:cmc}
 C_m^c = g.C^c(h') \quad \mbox{ for }\quad 
C^c(h') := C \cap \fq_1(h') - C \cap \fq_{-1}(h').
\end{equation}
This leads to the relation 
\begin{equation}
  \label{eq:exp-cmc}
\Exp_m(C_m^{c,0}) 
= \Exp_{gH}(g.C^c(h')^0)
= g.\Exp_{eH}(C^c(h')^0),
\end{equation}
which in turn entails with $ G^h g= g G^{h'}$
\begin{align} 
  \label{eq:hhprime-rel}
 W_M(h)_{gH} 
&= G^h_e.\Exp_{gH}( C_m^{c,0}) 
= G^h_e g.\Exp_{eH}(C^c(h')^0) 
= g G^{h'}_e.\Exp_{eH}(C^c(h')^0) \notag\\
&= g.W_M(h')_{eh} 
\quad \mbox{ for } \quad gH \in M^\alpha.
\end{align} 
Example~\ref{ex:ct}  shows that the dimensions of the eigenspaces of 
$\ad h'$ on $\fq$ are not always the same.  
\end{rem}

\begin{rem} (On the assumption $h \in \fh$) \\
\nin (a) The domains $W_M^+(h)$ and $W_M^{\rm KMS}(h)$ are defined for any 
Euler element $h \in \cE(\g)$. They only require the corresponding 
flows in $M$ and $\cT_M(C)$. However, these domains may be trivial 
if $h \not\in \fh$ (cf.~Proposition~\ref{prop:h0-grptype}). 

\nin (b) If  the vector field $X_h^M$ has a zero, then we may choose 
the base point accordingly to obtain $h \in \fh$. An Euler element 
$h \in \cE(\g)$ has this property if and only if its adjoint orbit 
$\cO_h$ intersects~$\fh$ (cf.\ \eqref{eq:modflow-M} in 
Definition~\ref{def:modflow}). 
\end{rem}

For a compactly causal symmetric Lie algebra $(\g,\tau,C)$, 
there may be many different cones $C'$ for which 
$(\g,\tau,C')$ is compactly causal, 
but there is a rather explicit classification of all these 
cones which is described in \cite[Thm.~3.6]{NO21a} in some 
detail in the $c$-dual context of non-compactly causal spaces. 
It implies in particular that 
$C$ is contained in a uniquely determined 
maximal $\Inn_\g(\fh)$-invariant elliptic 
cone $C_\fq^{\rm max}$ (which need not be pointed; 
see also \eqref{eq:min-max-g-red} and the 
discussion in Subsection~\ref{subsec:redliealg}). 

Recall the definition of $C_\pm := (\pm C)\cap \g_{\pm 1}$ from Remark \ref{rem:cpm}. The following 
proposition follows immediately by $c$-duality from 
\cite[Prop.~3.8]{NO21b}. It has the interesting consequence that, 
if $\g$ is reductive, then the cones $C_\pm$ remain the same when we replace $C$ by the maximal cone $C^{\rm max} \subeq \fq$ containing $C$ 
(see also Corollary~\ref{cor:red-minmax}(b) for spaces of group type). 

\begin{prop} \mlabel{prop:reductiontominmax}
Let $(\g,\tau,C,h)$ be a reductive modular compactly causal 
 symmetric Lie algebra. 
Then 
\begin{equation}
  \label{eq:cpm2}
 C_+ - C_- = C \cap \fq^{-\tau_h} = C^{\rm max} \cap \fq^{-\tau_h}.
\end{equation}
\end{prop}

This proposition implies in particular that the 
wedge domain $W_M(h)$ remains the same if we replace $C$ by $C^{\rm max}$. 
So it only depends on the ``direction'' of the cone $C$, not on its specific 
shape. This has interesting consequences for the global structure of $W_M(h)$. 

\begin{rem} \mlabel{rem:4.10} 
(Factorization of wedge domains) Let 
$(\g,\tau,C)$ be a reductive compactly causal symmetric Lie algebra. 
Then $(\g,\tau)$ decomposes as 
\[ (\g,\tau) 
\cong (\g_0, \tau_0) \oplus \bigoplus_{j = 1}^N (\g_j, \tau_j), \]
where $\g_0$ contains the center and all compact ideals 
and the $\tau$-invariant ideals $\g_j$, $j \geq 1$, are 
either simple hermitian or irreducible of group type 
(cf.\ \cite[Prop.~2.5]{NO21a} and Subsection~\ref{subsec:redliealg}). 
Decomposing $\fq = \fq_0 \oplus \bigoplus_{j = 1}^N \fq_j$ accordingly, 
\[ C^{\rm max} = \fq_0 + \sum_{j = 1}^N C^{\rm max}_j \quad \mbox{ with }\quad 
C_j^{\rm max} = C^{\rm max}\cap \fq_j,\] 
where the cones $C_j^{\rm max}$ are pointed. 
If $M$ is the simply connected symmetric space corresponding to $(\g,\tau)$, 
it follows that 
\begin{equation}
  \label{eq:m-prod}
 M \cong M_0 \times M_1 \times \cdots \times M_N,
\end{equation}
and
\begin{equation}
  \label{eq:wedge-prod}
W_M(h) = M_0 \times \prod_{j = 1}^N W_{M_j}(h_j),
\end{equation}
where $h = h_0 + \sum_{j = 1}^N h_j$ with $h_j \in \g_j$ and $h_0 \in \fz(\g)$. 

We have a similar global decomposition if $G$ is semisimple with 
$Z(G) = \{e\}$. 
Then 
\[ G = \Aut(G)_e \cong \prod_{j = 0}^N \Aut(\g_j), \quad 
H = G^\tau \cong \prod_{j = 0}^N \Inn(\g_j)^{\tau_j}, \quad 
M  \cong \prod_{j = 0}^N M_j.\] 
\end{rem}

\begin{ex} (The case of compact Lie algebras) 
If $\g$ is a compact Lie algebra, then  
$\ad x$ has purely imaginary spectrum for every $x \in \g$, 
so that $\cE(\g) = \eset$. However, we may consider central elements 
$h\in\fz(\g)$ as a degenerate kind of Euler elements with $\ad h =0$ 
and $\g = \g_0(h)$. 

There are many compactly 
causal symmetric spaces for which $\g$ is a compact Lie algebra. 
For instance, a compact Lie algebra $\g$ contains pointed generating invariant 
cones $C$ if and only if $\fz(\g) \not=\{0\}$, and then 
$\fz(\g) \cap C^0 \not=\eset$ holds for any such cone. 
The Lie algebra 
\[ \g = \fu_n(\C) \quad \mbox{ with } \quad 
 C := \{ X \in \g \: - i X \geq 0 \}\] 
is an important example. 
The corresponding causal group $(\U_n(\C),C)$ is the conformal compactification 
of the euclidean Jordan algebra $\Herm_n(\C)$, which, for $n = 2$, is 
isomorphic to the $4$-dimensional Minkowski space 
(\cite{FK94}). Involutions $\tau$ of $\fu_n(\C)$  
with $\tau(C) = -C$ are obtained from involutive 
automorphisms of $\su_n(\C)$, extending by $\tau(i\1) = -i\1$. 

For any compactly causal symmetric Lie algebra $(\g,\tau,C)$ 
and $h \in \fz(\g)$, we have $C_\pm = \{0\}$, so that 
$W_M(h) = \eset$. Further, the triviality of the modular flow 
implies $W_M^{\rm KMS}(h) = \eset$, but if 
$h \in \fz(\g)^{-\tau} \cap C^0$, then $W_M^+(h) = M$. 
\end{ex}

\section{Wedge domains in spaces of group type} 
\mlabel{sec:4}

For our analysis of wedge domains in compactly causal symmetric spaces, 
we shall follow the strategy to first study spaces of group 
type $G \cong (G \times G)/\Delta_G$, 
and then use embeddings into these spaces implemented by 
the quadratic representation to derive corresponding 
results in general. 

Let $(G,C)$ be a causal symmetric spaces of group type, i.e., 
$G$ is a connected Lie group and $C \subeq \g$ a 
pointed generating invariant cone. We further assume that the universal 
complexification $\eta_G \:  G \to G_\C$ has discrete kernel. 
We consider 
$G$ as the symmetric space $(G \times G)/\Delta_G$ with $\tau^{G \times G}(g_1, g_2) 
= (g_2, g_1)$ and $G \times G$ acting by $(g_1, g_2).g = g_1 g g_2^{-1}$. 
The associated causal symmetric Lie algebra is 
\[ (\g \oplus \g, \tau, C),\quad \mbox{ where } \quad 
\tau(x,y) = (y,x),\] 
and 
\[ C \subeq \g \cong \fq = \{ (x,-x) \: x \in \g \}  \] 
is an $\Ad(G)$-invariant pointed generating invariant cone. 
The Euler element $h \in \cE(\g)$ defines the Euler element 
$(h,h) \in \g \times \g$ which generates the flow 
$\alpha_t(g) = \exp(th) g \exp(-th)$ on $G$. 
If $\tau_h^G$ is an  involution on $G$ integrating 
$\tau_h = e^{\pi i \ad h}$ (cf.\ Remark~\ref{rem:2.12}), then 
this involution extends to the involution 
$\tau^G_h \times \tau_h^G$ on $G \times G$. 

We recall the closed/open complex Olshanski semigroup  
$\Gamma_G(C)$ 
from Definition~\ref{def:cplxols}. 
We start with the preparation of the 
analysis of wedge domains in 
modular symmetric spaces of group  type.

\begin{prop} \mlabel{prop:ne19b} 
Let $G$ be a connected Lie group with Lie algebra $\g$, 
$h \in \cE(\g)$ an Euler element and $C \subeq \g$ 
an $\Ad(G)$-invariant pointed closed convex cone such that: 
\begin{itemize}
\item[\rm(i)] $-\tau_h(C) = C$.
  \begin{footnote}
{If $\g$ is not reductive and $C$ is pointed {\it and} generating, 
then $C \cap \fz(\g) \not=\{0\}$ 
(\cite[Thm.~VII.3.10]{Ne00}).\\ As $\tau_h\res_{\fz(\g)} = \id_{\fz(\g)}$,  
the condition $-\tau_h(C) = C$ can only be satisfied if 
$C$ is {\it not} generating (see Example~\ref{ex:2.6}).}
  \end{footnote}
\item[\rm(ii)] $\tau_h$ integrates to an automorphism $\tau_h^G$ of~$G$. 
\item[\rm(iii)] The kernel of the universal complexification 
$\eta_G \: G \to G_\C$ is discrete. 
\end{itemize}
Let $g^\sharp=\tau (s)^{-1}$, $g\in G$.
Then the following assertions hold: 
\begin{itemize}
\item[\rm(a)] 
Let $C^c := C_+ + C_-$ with $C_\pm := \pm C \cap \g_{\pm 1}(h)$. Then the convex cones
$C^c, C_+$ and $C_-$ are $G^h$-invariant and
the subset 
  \begin{equation}
    \label{eq:Cc}
\Gamma_{G^h}(C^c) := G^h  \exp(C^c)
  \end{equation}
is a real Olshanski semigroup, associated to the symmetric 
Lie group $(G,\tau_h^G)$. In particular, 
the polar map $G^h \times C^c \to \Gamma_{G^h}(C^c)$ is a homeomorphism and
$\Gamma_{G^h}(C^c)$ is invariant under $s\mapsto s^\sharp$. 
\item[\rm(b)] The semigroup 
\[ S(C,h) := \{ g \in G \:  h - \Ad(g)h \in  C\} \] 
is invariant under $s\mapsto s^\sharp$ and
satisfies 
\begin{align}
  \label{eq:semeq}
S(C,h) & = \exp(C_+) G^h \exp(C_-) = G^h\exp (C_+)\exp (C_-)\nonumber \\
&=
\exp (C_-)\exp (C_+) G^h = \exp (C_- ) G^h \exp (C_+)= \Gamma_{G^h}(C^c) .
\end{align} 
If $C$ is generating,  
\begin{equation}
  \label{eq:semeq-int}
S(C,h)^0 = \Gamma_{G^h}(C^{c,0}) = 
S(C^0,h) := \{ g \in G \: h - \Ad(g)h \in C^0\}.
\end{equation}

It is invariant under  
$g \mapsto g^\sharp $. 
\item[\rm(c)] The subset 
of those elements $g \in G$ for which the orbit map 
\[ \alpha^g \: \R \to G,\quad \alpha^g(t)=  \exp(th) g\exp(-th) \] 
extends to a continuous map 
$\alpha^g \: \oline{\cS_{\pi}} \to \Gamma_G(C)$, which is 
holomorphic on $\cS_\pi$ when composed with the natural map 
$\Gamma_G(C) \to G_\C$, coincides with 
the closed semigroup $\Gamma_{G^h}(C^c)$. 
If $C$ has interior points, then 
  \begin{equation}
    \label{eq:kms}
 W_G^{\rm KMS}(h) = \Gamma_{G^h}(C^{c,0}). 
  \end{equation}
\item[\rm(d)] Let $X_h^G \in \cV(G)$ be the vector field defined by 
$X_h^G(g) :=  \frac{d}{dt}\big|_{t = 0}  \alpha_t(g).$ 
We write 
\[ G \times TG \to TG, \quad (g,x) \mapsto g.x \] 
for left translation of vectors on $G$. Then 
\begin{equation}
  \label{eq:SCh}
 S(C,h) = \{ g \in G \: X_h^G(g) \in g.C \}, 
\quad \mbox{ and } \quad 
 S(C,h)^0 = \{ g \in G \: X_h^G(g) \in g.C^0 \} \
\end{equation}
if $C$ is generating. 
\item[\rm(e)] On 
$\Gamma_G(C)$ we consider the antiholomorphic involution defined by 
\begin{equation}
  \label{eq:antiholfix}
 \oline\tau_h(g \exp(ix)) = \tau_h^G(g) \exp(-i \tau_h(x)).
\end{equation}
Its set of fixed points coincides with 
\begin{equation}
  \label{eq:poldecreal}
\Gamma_G(C)^{\oline\tau_h} = G^{\tau_h^G} \exp(i C^{-\tau_h}) 
\end{equation}
and, if $C$ has interior points, then the interior 
of the identity component $\Gamma_G(C)^{\oline\tau_h}_e$ 
satisfies 
\[ S(C,h)_e^0 = \alpha_{-\pi i/2}\big( (\Gamma_G(C)^{\oline\tau_h})_e^0\big).\] 
\end{itemize}
\end{prop}

\begin{prf} Similar results are stated in \cite{Ne19} for 
simply connected Lie groups $G$. If  
$q_G \: \tilde G \to G$ is the universal covering group of $G$, 
they apply to the group $\tilde G$. It therefore suffices 
to derive everything from \cite{Ne19}. 

\nin (a)  That $C^c, C_+$ and $C_-$ are $G^h$ invariant follows from
$\Ad (G^c)\g_j(h)=\g_j(h)$ and the invariance of $C$. The
remainder of (a) follows by applying  \cite[Prop.~2.6]{Ne19} to  $\tilde G$ 
because $\ker(q_G) \subeq \tilde G^h$ and 
$G^h = q_G(\tilde G^h)$. 

\nin (b) First we note that for $g\in G$ we have
\[-\tau (h-\Ad (g^\sharp)h )= \Ad (g^{-1})h -h =\Ad (g^{-1})(h-\Ad (g)h) \in C\]
As $C$ is $G$-invariant.
Thus $h- \Ad (g^\sharp)h \in -\tau (C)=C$.
That the four middle sets in \eqref{eq:semeq} agree now
 follows
now from the $\Ad (G^h)$ invariance of $C_\pm$.

For $\tilde S(C,h) := \{ g \in \tilde G \: h-\Ad(g)h \in C \}$ 
we have $\tilde S(C,h) = q_G^{-1}(S(C,h))$, so that the 
corresponding result for $\tilde G$ (\cite[Thm.~2.16]{Ne19}) 
and the invariance of $G^h, \tilde S(C,h)$ and $\tilde \Gamma_{G^h}(C^c) 
:= \tilde G^h \exp(C^c)$  
under multiplication with the central subgroup $\ker(q_G)$ 
show that it also holds for~$G$. 

Assume, in addition, that $C$ is generating, so that $C^0 \not=\eset$. 
Then $S(C^0,h) \subeq S(C,h)$ is an open subsemigroup, 
hence contained in 
\[ S(C,h)^0 = G^h \exp(C^{c,0}) = G^h \exp(C_+^0) \exp(C_-^0)=\exp (C_-^0)\exp (C_+^0)G^h.\] 
To show equality, write an element $g \in S(C,h)^0$ as 
\[ g =  \exp(x_-) \exp(x_+) g' \quad \mbox{ with } \quad 
g' \in G^h_C, x_\pm  \in C_\pm^0.\] 
Then 
\begin{align*}
\Ad (g) h &= e^{\ad x_-}e^{\ad x_+}h = e^{\ad x_-}(h +  [x_+,h] ) =e^{\ad x_-}(h-x_+)\\
&= h - x_+ +[x_-,h] - [x_-,x_+]  - [x_-, [x_-,x_+]]\\
&= h-\Ad (\exp x_-) (x_+ - x_- )\in h-C^0
\end{align*}
as $C^0$ is $G$-invariant. It follows that $h-\Ad (g)h \in C^0$. 
and hence $g \in S(C^0,h)$. 

\nin (c) For $\tilde G$ and the closed cone, 
this is \cite[Thm.~2.21]{Ne19}. As 
$\tilde G^h$ and the specified subset of $G$ 
are $\ker(q_G)$-invariant, the assertion also holds for~$G$.

For $z= x+iy \in \g_\C$ with $x,y\in\g$, we write $x= \Re z$ and $y = \Im z$. 
To verify \eqref{eq:kms}, for $g = g_0 \exp(x_- + x_+)$ with $g_0 \in G^h$, 
$x_\pm \in C_\pm = C \cap \g_{\pm 1}(h)$, and 
$0 < t < \pi$, we   note that $\sin (t) >0$ and
\begin{align*}
 \alpha_{it}(g_0\exp(x_++ x_-)) 
&= g_0\exp(e^{it} x_+ + e^{-it} x_-) \\
&=g_0 \exp ((\cos t + i\sin t)x_+ + (\cos t- i\sin t)x_-).\end{align*}
Thus 
\[ \Im(e^{it}x_+ e^{-it} x_-) = \sin(t)(x_+ - x_-).\] 
This element is contained in the interior $C^0$ if and only if 
$x_\pm \in C_\pm^0$, and in this case 
\[ \exp(e^{it} x_+ + e^{-it} x_-) 
\in \exp(\g + i C^0) \subeq \Gamma_G(C^0).\] 
We therefore obtain \eqref{eq:kms}. 

\nin (d) With the notation for left and right translations by $G$ 
on $TG$ (see the notation introduced in the introduction), we have 
\[ X_h^G(g) = h.g - g.h = g.(\Ad(g)^{-1}h -h).\] 
Hence \eqref{eq:SCh} follows immediately from 
$\Ad(g)C = C$, \eqref{eq:semeq} and \eqref{eq:semeq-int}.

\nin (e) Equation~\eqref{eq:poldecreal} follows from 
\eqref{eq:antiholfix}. We observe that 
\[ \alpha_{\pi i/2}(S(C,h)_e) 
= (G^h)_e \exp\big(e^{\frac{\pi i}{2} \ad h}(C_+ + C_-)\big)
= (G^h)_e \exp(i(C_+ - C_-)), \] 
so that $\alpha_{\pi i/2}(S(C,h)_e)$ is the identity component 
of the fixed point set $\Gamma_G(C)^{\oline\tau_h}$ 
in the closed Olshanski semigroup $\Gamma_G(C)$. 
If $C$ has interior points, this argument shows that 
$\alpha_{\pi i/2}(S(C,h)_e^0)$ is the connected component 
of the fixed point submanifold in the open semigroup 
$\Gamma_G(C^0)$ whose closure contains the identity. 
This proves (e). 
\end{prf}

To provide a good context to deal also with group 
type spaces that are not reductive, we have to deal 
with invariant cones $C \subeq \g$ which are not generating. 
This is mainly due to the fact that the requirement of $C$ being generating 
is not compatible with $-\tau_h(C) = C$. 
As in Subsection~\ref{subsec:2.1}, this problem can be overcome by observing that 
\[ \g_C = C - C \trile \g \] 
is an  ideal of the Lie algebra $\g$. We consider the corresponding 
integral subgroup $G_C \trile G$, endowed with its natural 
Lie group structure (it is closed if $G$ is simply connected). 
Then 
\[ (G \times G)_C := \{ (g_1, g_2) \in G \: g_1 g_2^{-1} \in G_C \} \] 
is a Lie group containing the diagonal $\Delta_G$ and isomorphic to 
$G_C \rtimes G$, where $G$ acts by conjugation on $G_C$. 
Then 
\[ G_C \cong (G \times G)_C/\Delta_G \] 
is a symmetric space on which 
\[  V_+(g) := g.C^0 \] 
defines a $(G \times G)_C$-invariant field of pointed open cones 
(a causal structure). The corresponding symmetric Lie algebra is 
\[ (\g \oplus \g)_C = \{ (x,y) \in \g \oplus \g \: x -y \in \g_C \} 
\quad \mbox{ with } \quad \tau(x,y) = (y,x).\] 
Here $\fq = \{ (x,-x) \: x \in \g_C\}\cong \g_C$ contains the pointed generating cone
$\{(x,-x)\: x\in C\} \cong
 C$. 
We now determine the wedge domains in the causal symmetric space $M = G_C$ 
for the structure defined by $((\g \oplus \g)_C, \tau, C)$. 

We do not assume that $C$ is generating in $\g$, but it is natural 
to assume that the two cones $C_\pm$ generate $\g_{\pm 1}(h)$ to 
ensure that the semigroup $S(C,h)$ has interior points. 

\begin{thm} {\rm(Wedge domains in spaces of group type)} 
\mlabel{thm:4.2}
 Suppose that $C_\pm$ generate $\g_{\pm 1}(h)$. 
Then the wedge domains in the compactly causal symmetric space 
$(G_C,C)$ corresponding to  the modular compactly causal symmetric Lie algebra 
$((\g \oplus \g)_C, \tau, C, (h,h))$ are 
\[   W_{G_C}^+(h) = W_{G_C}^{\rm KMS}(h)= W_{G_C}(h) = G_C^h \exp(C_+^0 + C_-^0).\] 
\end{thm}

\begin{prf} The fixed point set of the modular flow in $M = G_C$ is 
the centralizer $G_C^h$ of $G$. For a pair $(g_1, g_2) \in (G \times G)_C$, 
we have $m := (g_1, g_2).e = g_1 g_2^{-1} \in G_C^h$ if and only if 
$\Ad(g_1)^{-1} h= \Ad(g_2)^{-1}h$, i.e., $\Ad(g_1,g_2)^{-1}(h,h) \in \fh = \Delta_\g$. 
The invariance of $C$ under $\Ad(G)$ implies the invariance 
of the cones $C_\pm = \pm C \cap \g_{\pm 1}(h)$ under 
$\Ad(G^h)$. Therefore $C^c = C_+ + C_- \subeq \g^{-\tau_h}$ is invariant under 
$\Ad(G^h)$. 

For $m = g_1 g_2^{-1} \in G_C^h$, the modular flow 
acts on $T_m(G_C) = g_1.\g_C.g_2^{-1} = m.\g_C$ 
by 
\[ \alpha_t(m.x) = m.e^{t \ad h}x, \] 
so that $C_m^c := m.C^c$ (here we use the dot notation for left and right translations on $T(G)$). 
Therefore the polar wedge domain in $M = G_C$ is 
\begin{equation}
  \label{eq:polarwedge-group} 
 W_{G_C}(h) 
= \bigcup_{m \in G_C^h}  m \exp(C^{c,0})
= G_C^h\exp(C^{c,0}).
\end{equation}
Next we derive from Proposition~\ref{prop:ne19b}(c) that 
\begin{equation}
  \label{eq:wg2}
W_{G_C}^{\rm KMS}(h) =  G^h_C \exp(C_-^0 + C_+^0) = W_{G_C}(h).
\end{equation}
Finally, we observe that Proposition~\ref{prop:ne19b}(a),(d) imply that 
\begin{equation}
  \label{eq:wg3}
W_{G_C}^+(h) = G_C^h \exp(C_+^0 + C_-^0) = W_{G_C}(h).
\end{equation}
The assertion now follows from 
\eqref{eq:wg2} and \eqref{eq:wg3}.
\end{prf}

The open subset $W_{G_C}(h)$ neet not be connected, but its connected
components are parametrized by 
the group $\pi_0(G_C^h)$ of connected components of~$G_C^h$.

\begin{rem} \mlabel{rem:pi0h} (Some information on $\pi_0(G^h)$) 

\nin (a) As the example $G = \tilde\SL_2(\R)$ with the Cayley type 
involution
\[ \tau\pmat{a & b \\ c & d} = \pmat{a & -b \\ -c & d}\] 
shows, $G^h$ is in general not connected. It contains $Z(G) \cong \Z$ 
and 
\[ (G^h)_e = \exp (\R h) \quad \mbox{  for } \quad 
h = \frac{1}{2}\pmat{1 & 0 \\ 0 & -1}. \]
This shows that $\pi_0(G^h) \cong \Z$.  
This example also shows that, in general 
\[ G^h \not\subeq G^{\tau^G},\] 
because here $\tau^G$ acts on $Z(G)$ by inversion 
(cf.~\cite[Ex.~2.10(d)]{MN21}). 

\nin (b)  If $G$ is contained in $G_\C$, then 
$(G_\C)^h$ coincides with the centralizer of the circle group 
$\oline{\exp(ih)}$, hence is connected (\cite[Cor.~14.3.10]{HN12}). Therefore 
$G^h = (G_\C)^h \cap G$ is the fixed point group of the complex conjugation 
(with respect to $G$) on the connected complex group $(G_\C)^h$. Thus 
$\pi_0(G^h)$ is a finite elementary abelian $2$-group, i.e., isomorphic to 
$\Z_2^k$ for some $k\in \N_0$ (\cite[Thm.~IV.3.4]{Lo69}).
\end{rem}

\section{Wedge domains in 
extendable causal symmetric spaces} 
\mlabel{sec:5}

In this section we turn to the wedge domains in compactly 
causal symmetric spaces. {\it We 
assume that $(\g,\tau,C)$ is extendable in the sense of
Subsection}~\ref{subsec:2.1},  i.e., there exists a $G$-invariant
convex cone $C_\g\subset \fg$ such that $-\tau C_\g = C_\g$ and $C_\g\cap \fq = C$. 
This always holds if $\g$ is reductive by the Extension Theorem~\ref{thm:extend}.
We write $\g_C = C_\g - C_\g \trile \g$ for the ideal generated by $C_\g $ and
$G_C \trile G$ for the connected normal subgroup whose Lie algebra is $\g_C$  (see \eqref{eq:quademb}).  

The main result of this section 
is Theorem~\ref{thm:6.5}, asserting that 
all three wedge domains in $M = G/H$ are the same. 
This is first proved for the special case 
$H = G^{\tau^G}$ in Theorem~\ref{thm:6.4}. 
The proof of this theorem builds heavily on the group case 
(Theorem~\ref{thm:4.2}). We conclude this section with a brief 
discussion of the assumption that the Euler element $h$ is contained 
in $\fh$, i.e., that the corresponding modular flow on $M$ has a fixed point.

The following proposition identifies the basic connected components 
of the wedge domains in $G^\sharp_e$. 
 
\begin{prop} \mlabel{prop:5.2}
Let $(G,\tau^G)$ be a connected symmetric Lie group corresponding 
to the modular compactly causal symmetric Lie algebra 
$(\g,\tau,C,h)$, 
\[ M_G := G/G^{\tau^G} \cong G^\sharp_e = (G_C^\sharp)_e,\]
 and let $C_\g \subeq \g$ be 
a pointed invariant closed convex cone with 
\[ C_\g \cap \fq = C  \quad \mbox{ and }\quad 
\quad -\tau(C_\g) = -\tau_h(C_\g) = C_\g.\] 
For the real Olshanski semigroup 
$S(C_\g,h) = G_C^h \exp(C_\g^c)$, the following assertions hold: 
\begin{itemize}
\item[\rm(a)] The basic connected component of the 
wedge domains in $M_G$ coincide:  
\[ W_{M_G}(h)_{e} = W_{M_G}^+(h)_{e} = W_{M_G}^{\rm KMS}(h)_{e},\] 
and this domain can be obtained from the semigroup 
$S(C_\g,h) = G_C^h \exp(C_\g^c)$ as 
\begin{equation}
  \label{eq:wedge-s-fix}
 S(C_\g^0,h)^\sharp_e 
= \bigcup_{g \in (G_C^h)_e} g \exp(C_+^0 + C_-^0) g^\sharp 
\cong (G_C^h)_e \times_{(G_C^h)_e^{\tau}} (C_+^0 + C_-^0). 
\end{equation}
\item[\rm(b)] $W_{M_G}(h)_e$ is the orbit $S(C_\g^0,h)_e.e$ of the base point under 
the action of the open semigroup $S(C_\g^0,h)_e$ on ${M_G} = (G_C^{\sharp})_e$ 
by $g.x := gxg^\sharp$. 
\end{itemize}
\end{prop}

\begin{prf}  
The flow $\alpha_t$ 
commutes with the involution $\tau$. 
If $m \in {M_G} \cong (G_C^{\sharp})_e \subeq G_C^\sharp$ 
is contained in  $W_{M_G}^{\rm KMS}(h)$ 
(Definition~\ref{def:wedgedom}), then 
\[ \alpha_z(m) \in 
\cT_{M_G}(C) \subeq \cT_{G_C}(C_\fg) =  \Gamma_{G_C}(C_\g^0)
\quad \mbox{ for } \quad z \in \cS_\pi \]  
(cf.\ Remark~\ref{rem:exatube}), so that 
Theorem~\ref{thm:4.2} implies that 
\begin{equation}
  \label{eq:weq}
 W_{M_G}^{\rm KMS}(h) \subeq 
{M_G}  \cap W_{G_C}^{\rm KMS}(h) 
= {M_G}  \cap S(C_\g^0, h) \subeq S(C_\g^0,h)^\sharp,
\end{equation} 
and thus 
\begin{align}
  \label{eq:weq2}
 W_{M_G}^{\rm KMS}(h)_e 
&\subeq 
 S(C_\g^0,h)^\sharp_e = (G_C^h)_e.\Exp((C_\g^{c,0})^{-\tau}) \nonumber\\
& {\buildrel \eqref{eq:doubledcone2}\over =}\ (G_C^h)_e.\Exp((C_+ + C_-)^0)= W_{M_G}(h)_e.
\end{align}
The invariance of both sides under 
$(G_C^h)_e$ reduces the verification of the converse inclusion to showing 
that 
\begin{equation}
  \label{eq:kms-inc}
\Exp(C_+^0 + C_-^0) \subeq W_{M_G}^{\rm KMS}(h).
\end{equation}
For $x_\pm \in C_\pm^0$ and $0 < t < \pi$, we have $\sin t>0$ and hence
\[ \Im(e^{it}x_+ + e^{-it} x_-) = \sin(t)(x_+ - x_-)\in C^0.\] 
Now \eqref{eq:kms-inc}  follows
from 
$ \alpha_{it}\Exp(x_++ x_-) 
=\Exp(e^{it} x_+ + e^{-it} x_-) $, 
$\exp(e^{it} x_+ + e^{-it} x_-) \in \Gamma_{G_C}(C_\g^0)$ 
and 
\[ \cT_{M_G}(C) = \{ ss^\sharp \: s \in \Gamma_{G_C}(C_\g^0)\} 
\subeq  \Gamma_{G_C}(C_\g^0) = \cT_{G_C}(C_\g)\] 
(Remark~\ref{rem:exatube}(a)). 
We thus obtain 
\begin{equation}
  \label{eq:wedge-eq2}
 W_{M_G}^{\rm KMS}(h)_e =  (G_C^h)_e.\Exp(C_-^0 + C_+^0) = W_{M_G}(h)_e.
\end{equation}

Next we observe that $C^0 = C_\g^0 \cap \fq$ 
implies that, for $m = gg^\sharp \in {M_G} = (G_C^\sharp)_e \subeq G_C$, we have 
\begin{equation}
  \label{eq:ksm-h}
 V^{M_G}_+(m) = g.C^0.g^\sharp 
= T_m({M_G}) \cap  g.C_\g^0.g^\sharp = T_m({M_G}) \cap V_+^{G_C}(gg^\sharp).
\end{equation}
In view of Theorem~\ref{thm:4.2}, this shows that 
\begin{equation}
  \label{eq:conmg}
 W_{M_G}^+(h) = W_{G_C}^+(h) \cap {M_G} = S(C_\g, h)^0 \cap {M_G}. 
\end{equation}

The semigroup $S(C_\g,h) = G_C^h \exp(C_\g^c)$ is a real Olshanski semigroup 
corresponding to the involution $\tau_h$ on~$\g$ 
(Proposition~\ref{prop:ne19b}). 
Further, 
$\tau$ and $\tau_h$ commute and the cone $C_\g^c$ is $-\tau$-invariant 
(cf.\ \eqref{eq:doubledcone2}), 
so that Lemma~\ref{lem:real-olsh}(3) shows that 
\[ S(C_\g,h)^\sharp_e 
= (G_C^h)_e.\exp((C_\g^c)^{-\tau})
= (G_C^h)_e.\exp(C^c) = (G_C^h)_e.\exp(C_+ + C_-).\] 
We thus conclude with \eqref{eq:wedge-eq2} 
and Lemma~\ref{lem:real-olsh}(3,4) that 
\[  W_{M_G}^+(h)_e \ 
\ {\buildrel {\rm\ref{lem:real-olsh}}\over = }\ 
\  S(C_\g^0,h)^\sharp_e\ \ {\buildrel\eqref{eq:wedge-eq2} 
\over =}\ \  W_{M_G}^{\rm KMS}(h)_e.\] 
Now (a) follows from \eqref{eq:wedge-eq2}, 
so that  (b) is a direct consequence of \eqref{eq:weq2}. 
\end{prf}

The preceding proposition identifies the basic components 
of the wedge domains in the compactly causal symmetric space $M_G$. 
We now prepare the ground for the analysis of the 
other connected components. 
We recall the polar wedge domain 
\begin{equation}
  \label{eq:equivwedgedom}
W_M(h) 
= \bigcup_{m\in M^\alpha} \Exp_m(C^c_m(h)^0) 
\ {\buildrel \eqref{eq:exp-cmc}\over =}
\  \bigcup_{gH\in M^\alpha} g.\Exp_{eH}(C^c(\Ad(g)^{-1}h)^0) \subeq M.
\end{equation}
Comparing with Definition~\ref{def:wedgedom} and using that 
$M^\alpha_{eH} = G^h_e.eH$ (Proposition~\ref{prop:5.3}), we now see that 
\begin{equation}
  \label{eq:wedge-conn-comp}
 W_M(h)_{eH} = G^h_e.\Exp_{eH}(C^c(h)^0) 
= \bigcup_{m\in M^\alpha_{eH}} \Exp_m(C^c_m(h)^0).
\end{equation}

\begin{rem} As the connected components of $M^\alpha$ may have different dimensions 
(Example~\ref{ex:ct}) 
and the cones $C^c_m \subeq T_m(M)$ span complements to $T_m(M^\alpha)$, they 
may also be of different dimension.   
\end{rem}

\begin{thm} \mlabel{thm:6.4} 
{\rm(Wedge domains in extendable compactly causal symmetric spaces)} 
Let $(G,\tau^G)$ be a connected symmetric Lie group corresponding 
to the modular compactly causal symmetric Lie algebra 
$(\g,\tau,C,h)$, 
\[ M = M_G = G^\sharp_e \cong G/G^{\tau^G} \cong (G_C^\sharp)_e,\]
 and let $C_\g \subeq \g$ be 
a pointed invariant closed convex cone with 
\[ C_\g \cap \fq = C  \quad \mbox{ and }\quad 
\quad -\tau(C_\g) = -\tau_h(C_\g) = C_\g.\] 
Then 
\[  W_M(h) = W_M^+(h) = W_M^{\rm KMS}(h).\] 
\end{thm}

\begin{prf} We start with the discussion of the positivity domain 
$W_M^+(h)$. From  \eqref{eq:conmg} in the proof of 
Proposition~\ref{prop:5.2} we know that, for 
\[ S_C :=  S(C_\g,h) = G_C^h \exp(C_\g^c)
\quad \mbox{ and } \quad S_C^0 = G_C^h \exp(C_\g^{c,0}),\]  we have 
\[ W_{M}^+(h) = S_C^0 \cap M = S_C^0 \cap G^\sharp_e \subeq (S_C^0)^\sharp \] 
is the union of all connected components of $(S_C^0)^\sharp$ 
which are contained in $M = G^\sharp_e$. 
Information on these connected components comes from 
Lemma~\ref{lem:real-olsh}, 
which shows that each connected component of  
$(S_C)^\sharp$ intersects $(G_C^h)^\sharp$, the $\alpha$-fixed point set 
in $G_C^\sharp$. So let $g_0 \in M^\alpha \cap S_C$ 
and consider the involution $\gamma^G(g) := g_0^{-1} \tau(g) g_0$ of $G$ 
and the Lie algebra involution 
\[ \gamma := \Ad(g_0)^{-1} \tau 
= \Ad(\tau(g)) \tau \Ad(\tau(g))^{-1} \in \Aut(\g).\] 
By Lemma~\ref{lem:real-olsh}(5)(c), the connected component 
of $S^\sharp_C$ containing $g_0$ is 
\begin{equation}
  \label{eq:appb-semsubset}
g_0 \{ s\gamma^G(s)^\sharp \: s \in S_{C,e} \}
= g_0 \bigcup_{g_1 \in (G^h_C)_e} g_1 \exp(C_\g^{c,-\gamma}) g_0^{-1} g_1^\sharp g_0.
\end{equation}

We write $g_0 = gg^\sharp$ with $g \in G$ 
and observe that $h' = \Ad(g)^{-1}h \in \fh$, so that we also have 
$h' = \tau(h') = \Ad(\tau(g))^{-1}h$. The Lie algebra involution 
$\gamma$ satisfies 
\begin{equation}
  \label{eq:hsym}
 \Ad(g^\sharp) \g^{-\gamma} =  \Ad(\tau(g))^{-1} \g^{-\gamma}  = \g^{-\tau}= \fq.
\end{equation}
Therefore the cone 
\[ C_\g^{c,-\gamma} = C_\g^c \cap \g^{-\gamma} 
= \big(C_\g \cap \g_1(h) - C_\g \cap \g_{-1}(h)\big) \cap \g^{-\gamma} \] 
satisfies 
\begin{equation}
  \label{eq:conetrafo}
 \Ad(g^\sharp)C_\g^{c,-\gamma}
= \big(C_\g \cap \g_1(h') - C_\g \cap \g_{-1}(h')\big) \cap \fq
= C \cap \g_1(h') - C \cap \g_{-1}(h') =: C^c(h'),
\end{equation}
where we have used that the projection $p_\fq \: \g \to \fq$ commutes 
with $\ad h'$ and maps $C_\g$ into itself. 

We now arrive with \eqref{eq:appb-semsubset} at 
\begin{align*}
(S_C)_{g_0}^\sharp 
&=  \bigcup_{g_1 \in (G_C^h)_e} g_0 g_1 \exp(C_\g^{c,-\gamma}) g_0^{-1} g_1^\sharp  g_0 \\
&=  \bigcup_{g_1 \in (G_C^h)_e} (g_0 g_1 g_0^{-1}) g_0 \exp(C_\g^{c,-\gamma}) (g_0 g_1 g_0^{-1})^\sharp  \\
&=  \bigcup_{g_2 \in (G^h_C)_e} g_2 g_0 \exp(C_\g^{c,-\gamma}) g_2^\sharp  
= (G^h_C)_e.(g_0\exp(C_\g^{c,-\gamma}))\\
&= (G^h_C)_e.\big(gg^\sharp\exp(C_\g^{c,-\gamma})\big)
=  (G^h_C)_e.\Big(g\exp\big(\Ad(g^\sharp) C_\g^{c,-\gamma}\big)g^\sharp\Big)\\
&\ {\buildrel\eqref{eq:conetrafo}\over=}\ (G^h_C)_e g.\exp(C^c(h')). 
\end{align*}
For the interior of this connected component, we obtain 
\begin{equation}
  \label{eq:concompfib}
(S_C^0)_{g_0}^\sharp = (G^h_C)_e g.\exp(C^c(h')^0). 
\end{equation}
The boundary of this domain contains in particular the connected 
component 
\[ M^\alpha_{g_0} = (G^h_C)_e g.e 
= \{g_1 gg^\sharp g_1^\sharp \: g_1 \in (G^h_C)_e\} \] 
(cf.\ Proposition~\ref{prop:5.3}). 

We are now ready to use this information to identify the positivity domain. 
As $X^M_h(g.m) \in V_+(g.m) = g.V_+(m)$ is equivalent to 
$g^{-1}.X^M_h(g.m) = X^M_{h'}(m) \in V_+(m)$, it follows that 
\begin{equation}
  \label{eq:posdom-trafo}
  g^{-1}.W_M^+(h) = W_M^+(h') = W_M^+(\Ad(g)^{-1}h).
\end{equation}
Combining Remark~\ref{rem:6.5} with 
Proposition~\ref{prop:5.2}, we conclude that 
\begin{align} \label{eq:inc-wm-in-wplus}
W_M(h) 
&= \bigcup_{m \in M^\alpha} W_M(h)_m 
\ {\buildrel \ref{rem:6.5} \over =}\ 
 \bigcup_{gH \in M^\alpha} g.W_M(\Ad(g)^{-1}h)_{eH} \notag \\
&\subeq  \bigcup_{gH \in M^\alpha} g.W_M^+(\Ad(g)^{-1}h)
\ {\buildrel\eqref{eq:posdom-trafo}\over=}\  W_M^+(h).
\end{align}
To obtain the converse inclusion, we use \eqref{eq:conmg} to obtain 
\begin{align} \label{eq:hatwm-in-wmplus}
 W_M^+(h) 
&\ {\buildrel\eqref{eq:conmg}\over=} \  S_C^0 \cap M
= S_C^0 \cap G^\sharp_e 
= \bigcup_{g \in G, gg^\sharp \in G^h} (S_C^0)^\sharp_{gg^\sharp}\notag \\
&\ {\buildrel \eqref{eq:concompfib}\over = }\ 
 \bigcup_{g \in G, \Ad(g)^{-1}h \in \fh} 
(G^h_C)_eg. \exp(C^c(\Ad(g)^{-1}h)) = W_M(h).
\end{align}
This completes the proof of the equality of $W_M(h)$ and $W_M^+(h)$. 

We now turn to the KMS wedge domain $W_M^{\rm KMS}(h)$. 
Since $G$ acts naturally on the tube domain $\cT_M(C)$, the modular 
flow $\alpha'$ defined by $h' = \Ad(g)^{-1}h$ on $M$ satisfies 
\[ \alpha_z(g.m) = g.\alpha_z'(m)  \quad \mbox{ for } \quad z \in \cS_\pi. \] 
This shows that 
\begin{equation}
  \label{eq:kms-trafo}
 g^{-1}.W_M^{\rm KMS}(h) 
= W_M^{\rm KMS}(h')
= W_M^{\rm KMS}(\Ad(g)^{-1} h).
\end{equation}
As for the positivity domain, we now obtain with 
Proposition~\ref{prop:5.2},
\begin{align} \label{eq:inc-wm-in-wkms} 
W_M(h) 
&= \bigcup_{m \in M^\alpha} W_M(h)_m 
= \bigcup_{gH \in M^\alpha} g.W_M(\Ad(g)^{-1}h)_{eH} \notag \\
&\subeq  \bigcup_{gH \in M^\alpha} g.W_M^{\rm KMS}(\Ad(g)^{-1}h)
\ {\buildrel \eqref{eq:kms-trafo} \over = }\  W_M^{\rm KMS}(h).
\end{align}
Combining \eqref{eq:weq} with \eqref{eq:hatwm-in-wmplus}, we finally get  
\[  W_M^{\rm KMS}(h)\ {\buildrel\eqref{eq:weq} \over \subeq }\ 
  S(C_\g^0, h) \cap M
=  S_C^0  \cap G^\sharp_e = W_M(h),\] 
and therefore the equality $W_M(h) = W_M^{\rm KMS}(h)$. 
\end{prf}

\begin{rem} The relation in \eqref{eq:posdom-trafo} shows that 
any $h' = \Ad(g)^{-1}h \in \cO_h$ leads to a wedge domain $W_M(h') = g^{-1}.W_M(h)$ 
that is a translate of the wedge domain $W_M(h)$. In this sense 
the geometric structure of the set $W_M(h)$ does not depend 
on the choice of the Euler element in $\cO_h \cap \fh$. 
Nevertheless the structure of the basic connected component 
$W_M(h)_{eH}$ of the wedge domain $W_M(h)$ does (Example~\ref{ex:ct}). 
\end{rem}

Eventually, we turn to the general case, where 
$M = G/H$ does not necessarily embed into~$G$.

\begin{thm} \mlabel{thm:6.5}
{\rm(Wedge domains in extendable compactly causal symmetric spaces)} 
Let $(G,\tau^G)$ be a connected symmetric Lie group corresponding 
to the modular compactly causal symmetric Lie algebra 
$(\g,\tau,C,h)$, assume that $\eta_G \: G \to G_\C$ has discrete kernel,  
let $H \subeq G^{\tau^G}$ be an open subgroup, 
$M := G/H$,  and let $C_\g \subeq \g$ be 
a pointed invariant closed convex cone with 
\[ C_\g \cap \fq = C  \quad \mbox{ and }\quad 
\quad -\tau(C_\g) = -\tau_h(C_\g) = C_\g.\] 
Then 
\begin{equation}
  \label{eq:allequal}
 W_M(h) = W_M^+(h) = W_M^{\rm KMS}(h).
\end{equation}
\end{thm}

\begin{prf} Let $\tilde\tau \in \Aut(\tilde G)$ be the involution 
integrating $\tau$. 
In the simply connected covering $\tilde G$ of~$G$, 
the subgroup $\tilde G^{\tilde\tau}$ is connected 
by \cite[Thm.~IV.3.4]{Lo69}, so that
\cite[Cor.~11.1.14]{HN12} implies that 
\[ M_{\tilde G} := \tilde G^\sharp_e \cong \tilde G/\tilde G^{\tilde\tau}  \] 
is simply connected, hence can be identified with the simply connected 
covering $\tilde M$ of~$M$. Accordingly, the universal covering map is 
\[ q_M \: \tilde M = \tilde G/\tilde G^{\tilde\tau}\to M = G/H, \quad 
g \tilde G^{\tilde\tau} \mapsto q_G(g) H,\] 
where $q_G \: \tilde G \to G$ is the universal covering of~$G$. 
From the isomorphism $\tilde M \cong M_{\tilde G}$ and 
Theorem~\ref{thm:6.4}, we obtain the equalities 
\begin{equation}
  \label{eq:tildem-equl}
  W_{\tilde M}^+(h) = W_{\tilde M}^{\rm KMS}(h) = W_{\tilde M}(h).
\end{equation}

The equivariance of $q_M$ with respect to the modular flow implies that 
\begin{equation}
  \label{eq:+domcov}
 W_{\tilde M}^+(h) = q_M^{-1}(W_M^+(h)), 
\quad \mbox{ hence } \quad 
 W_M^+(h) = q_M(W_{\tilde M}^+(h)).
\end{equation}
We also have coverings of the tube domains 
$\cT_{\tilde M}(C) \to \cT_M(C)$ 
which are  equivariant with respect to the modular flow. This 
entails 
\begin{equation}
  \label{eq:kms-dom-cover}
q_M^{-1}(W_{M}^{\rm KMS}(h))=  W_{\tilde M}^{\rm KMS}(h), 
\quad \mbox{ and thus } \quad 
 W_M^{\rm KMS}(h) = q_M(W_{\tilde M}^{\rm KMS}(h)).
 \end{equation}
In fact, the inclusion $\subeq$ follows directly from the  equivariance. 
For the converse inclusions  we use the existence of lifts of trajectories of 
the imaginary vector field $i \cdot X_h^{\cT_M}$ in $\cT_{M}(C)$ to 
$\cT_{\tilde M}(C)$. 

Next we observe that 
\begin{equation}
  \label{eq:cover-fp}
 \tilde M^\alpha = \{ m \in \tilde M \: X_h^{\tilde M}(m) = 0\} 
= q_M^{-1}(M^\alpha).
\end{equation}
As $T_m(q_M)V_+(m) = V_+(q_M(m))$ and $q_M$ is $\alpha$-equivariant, 
we have 
\[ T_m(q_M) C_m^c = C_{q_M(m)}^c,\] 
and thus, for $m \in \tilde M^\alpha$, 
\begin{align*}
 q_M(W_{\tilde M}(h)_m) 
&= q_M(\tilde G^h_e.\Exp_m(C_m^{c,0})) 
= G^h_e.\Exp_{q_M(m)}\big(T_m(q_M)C_m^{c,0}\big)\\
&= G^h_e.\Exp_{q_M(m)}\big(C_{q_M(m)}^{c,0}\big)
= W_M(h)_{q_M(m)}.
\end{align*}
Taking the union over all $\alpha$-fixed points, we arrive with 
\eqref{eq:cover-fp} at 
\begin{equation}
  \label{eq:63}
 q_M(W_{\tilde M}(h)) = W_M(h).
\end{equation}
Applying $q_M$ to the equalities in \eqref{eq:tildem-equl} now leads 
with \eqref{eq:+domcov},  \eqref{eq:kms-dom-cover} and 
\eqref{eq:63} to \eqref{eq:allequal}. 
\end{prf}

\subsection*{Assumptions on the Euler element} 
\mlabel{pers:1}

We assume throughout that the Euler element $h$ 
generating the modular flow on $M = G/H$ 
is contained in $\fh = \g^\tau$ 
(see in particular Proposition~\ref{prop:cc2}(b)). 

For group type spaces $(\g \oplus \g, \tau_{\rm flip})$, 
Euler elements in $\g \oplus \g$ 
are of the form $h = (h_1, h_2)$, where 
$h_j \in \cE(\g)$ are Euler elements.
\begin{footnote}
{If $\g$ is simple hermitian, then  $\Inn(\fg)$ acts transitively 
on $\cE(\fg)$ by Proposition~\ref{prop:cc2}, so that there are only $3$ orbits of 
$\Inn(\g \oplus \g)$ on $\cE(\g \oplus \g)$, represented 
by elements of the form $(h,h)$, $(h,0)$, $(0,h)$ with 
$h \in \cE(\g)$. } 
\end{footnote}
The following proposition shows that, for $h \not\in \fh$, the wedge 
domains may degenerate drastically. 

\begin{prop} \mlabel{prop:h0-grptype}
Let $(\g \oplus \g, \tau_{\rm flix}, C)$ be a 
causal symmetric Lie algebra  of group type. 
If $h_0 \in \cE(\g)$, then $h := (h_0,0) \in \cE(\g \oplus \g)$ and
\[ W_G(h) = W_G^+(h) = W_G^{\rm KMS}(h) = \eset.\] 
\end{prop} 

\begin{prf}  We write 
\[ C = \{ (x,-x) \: x \in C_\g \} \] 
for an invariant cone $C_\g \subeq \g$. 
The modular flow acts on $G$ by left translations 
$\alpha_t(g) = \exp(th_0)g$. 
For $g \in G$, the inclusion 
\[  \exp([0,\pi]i h_0) g \subeq G \exp(iC_\g) \] 
implies $h_0 \in C_\g$, contradicting the fact that 
$\ad h_0$ is diagonalizable and the 
elements of $C_\g$ have purely imaginary spectrum. This shows that 
$W_G^{\rm KMS}(h) = \eset.$ 

Moreover, 
\[ (\g \oplus \g)_{\pm 1}(h) = \g_{\pm 1}(h_0) \oplus \{0\} \] 
implies that $C_\pm = C \cap (\g \oplus \g)_{\pm 1}(h)= \{0\}$, so that 
$W_G(h)= \eset$. 

Finally, we observe that 
\[ p_\fq(\Ad(g_1, g_2)^{-1}(h_0,0))
= \frac{1}{2}(\Ad(g_1^{-1})h_0, - \Ad(g_1)^{-1}h_0)  \in C^0 \] 
is equivalent to $\Ad(g_1^{-1})h \in C_\g^0$, which is never the case 
because $C_\g^0$ consists of elliptic elements. This shows that 
$W_G^+(h) = \eset$. 
\end{prf} 

For right translations one argues similarly. 
On the other hand, the Euler elements of the form $(h,h)$ 
correspond to non-trivial wedge domains, such as 
\[ W_G(h) = G^h \exp(C_\g^{c,0}) \] 
(Theorem~\ref{thm:4.2}).

\section{Nets of standard subspaces} 
\mlabel{sec:6a} 

In this section we turn to representation theoretic aspects of 
compactly causal symmetric spaces. We start with 
introducing standard subspaces, the Brunetti--Guido--Longo construction 
and recall the construction of nets of closed real subspaces 
based on distribution vectors from \cite{NO21a}. 
In Subsection~\ref{subsec:7.4} we show in the general 
Theorem~\ref{thm:7.5} how the methods from \cite{NO21a} can be used 
to construct covariant nets of standard subspaces on 
compactly causal symmetric spaces.

\subsection{Standard subspaces} 
\mlabel{subsec:6a.1}

 In this section we recall the notion of standard subspaces and modular objects and how those
object are related to anti-unitary representation of $\R^\times$. 

Let $\cH$ be a complex Hilbert space. 
A closed real subspace $\sV \subeq \cH$ 
is called {\it standard} if 
\begin{equation}
  \label{eq:stansub}
  \sV \cap i \sV = \{0\} 
\quad \mbox{ and } \quad \cH = \oline{\sV + i \sV}
\end{equation} 
(see  \cite{Lo08} or \cite{NOO21} for basic facts on standard subspaces).   Associated to
every standard subspace is the closed densely defined conjugate linear
operator 
\[ \sigma_\sV \: \sV + i \sV \to \cH, \quad 
x + i y \mapsto x- iy \]
with polar decomposition
\[\sigma_\sV =J_\sV\Delta_\sV^{1/2} .\]
The operator $J_\sV$ is an everywhere defined conjugation (an antiunitary 
involution), and
the {\it modular operator}
$\Delta_\sV$ is a positive selfadjoint operator. These two operators satisfy the {\it modular relation }
$J_\sV \Delta_\sV J_\sV = \Delta_\sV^{-1}$.  

Conversely, every 
pair $(\Delta,J)$, consisting of a positive selfadjoint operator $\Delta$ 
and a conjugation satisfying the modular relation 
$J\Delta J  = \Delta^{-1}$ defines a standard subspace via 
\[ \sV := \Fix(J\Delta^{1/2}).\] 
In this sense standard subspaces are in
bijection with pairs of modular objects.

We also note that, for any pair $(\Delta, J)$ of modular objects, 
\begin{equation}
  \label{eq:modgroupdelta}
 U(e^t) := \Delta^{-it/2\pi} \quad \mbox{ and } \quad 
U(-1) := J
\end{equation}
define a continuous homomorphism 
$U \: \R^\times \to \AU(\cH)$ mapping negative 
numbers to antiunitary operators. Conversely, every 
such homomorphism is obtained for 
\begin{equation}
  \label{eq:modgrpdelta2}
\Delta := e^{2\pi i \cdot\partial U(1)}\quad \mbox{  and } \quad J := U(-1).  
\end{equation}

\subsection{The BGL construction} 
\mlabel{subsec:bgl}

In this section we generalize the above construction
for $\R^\times$ to more general Lie groups.  For an involutive automorphism $\tau^G$ on a connected 
Lie group $G$, we consider 
the group 
\[ G_\tau  = G \rtimes \{\id_G,\tau^G\}.\] 
An {\it antiunitary representation} of $G_\tau$ is a homomorphism 
$U \: G_\tau \to \AU(\cH)$ (the group of unitary and antiunitary operators on 
$\cH$) such that $U$ is strongly continuous and 
$J := U(\tau^G)$ is a conjugation. 
 We then have
\[ J U(g) J = U(\tau^G(g))\quad \mbox{  for } \quad g \in G.\] 
For  every $h \in \g$ fixed by $\tau := \L(\tau^G)$, we then obtain a 
standard subspace 
\[ \sV := \sV_{(h,\tau^G,U)}\subeq \cH,\]
 specified by 
\begin{equation}
  \label{eq:bgl2}
J_\sV = U(\tau^G) \quad \mbox{ and } \quad 
\Delta_{\sV}^{-it/2\pi} = U(\exp th) \quad \mbox{ for } \quad t \in \R
\end{equation}
as in \eqref{eq:modgroupdelta}. Here  the fact that 
$U(\tau^G)$ commutes with $U(\exp \R h)$ implies the modular relation, 
so that 
$\sV = \Fix(J_\sV \Delta_\sV^{1/2})$ 
is a standard subspace of $\cH$.
This assignment is called the {\it Brunetti--Guido--Longo 
(BGL) construction} (see \cite{BGL02}).

\subsection{Nets of real subspaces on homogeneous spaces} 
\mlabel{subsec:7.3}

With the BGL construction, standard subspaces can be associated to antiunitary 
representations in abundance, but only a few of them 
carry interesting geometric 
information. In particular, we would like to understand 
when a standard subspace of the form $\sV_{(h,\tau^G,U)}$ arises from a 
natural family $\sV(\cO)$ of real subspaces associated to open 
subsets of a homogeneous space $M = G/P$ and for which domains $\cO \subeq M$ 
the subspace $\sV(\cO)$ is standard. 

Covariant families of real subspaces of $\cH$ 
which are not necessarily standard are easy to construct 
in any unitary representation $(U,\cH)$ of $G$ 
from distribution vectors 
(see Appendix~\ref{app:smovec} for definitions and basic properties). 
To a real linear subspace  $\sE \subeq \cH^{-\infty}$ and 
an open subset $\cO \subeq G$,  we associate a 
closed real subspace of $\cH$ by 
\begin{equation}
  \label{eq:he1}
\sH_\sE(\cO) := \oline{\spann_\R \big(U^{-\infty}(C^\infty_c(\cO,\R))\sE\big)} 
\subeq \cH.
\end{equation}

With the projection map $q_P \: G \to G/P$, we then obtain 
on $G/P$ a net of real subspaces by 
\[ \sH^{G/P}_\sE(\cO) := \sH_\sE(q_P^{-1}(\cO)).\] 
If $\sE$ is $P$-invariant, then we actually have 
$\sH_\sE(\cO) = \sH_\sE(\cO P)$ (\cite[Lemma~2.11]{NO21a}), so that any subspace 
$\sH_\sE(\cO)$ attached to an open subset of $G$ also corresponds to 
an open subset of~$G/P$. 

\subsection{Connecting wedge domains and standard subspaces}
\mlabel{subsec:7.4}

We now connect the BGL construction with the construction 
based on real subspaces of distribution vectors. 
For an Euler element $h\in \g$ we consider the involution 
$\tau_h = e^{\pi i \ad h}$ and assume that the 
Lie algebra involution 
$\tau_h $ integrates to an involutive automorphism 
$\tau_h^G$ of the connected Lie group $G$ with Lie algebra~$\g$ 
(cf.\ Remark~\ref{rem:2.12}).
Suppose that 
$(U,\cH)$ is an antiunitary representation of $G_{\tau_h} = 
G \rtimes \{\1,\tau_h^G\}$ with 
$C_U$ pointed and that $\sV \subeq \cH$ is the standard subspace 
specified by the triple $(h,\tau_h^G,U)$ as in 
\eqref{eq:bgl2} by the BGL construction. We also assume that 
$\g = \g_C + \R h$ for $\g_C = C_U - C_U$ (cf.\ \cite[\S 3]{NO21a}). 
Then \cite[Thms.~2.16, 3.4]{Ne19} imply that, for $C_\g = C_U$, 
\begin{equation}
  \label{eq:sv2}
S_\sV = G_{\sV} \exp(C_{\g,+} + C_{\g,-}), 
\quad \mbox{ where } \quad 
\L(G_\sV) = \g^{\tau_h}, \qquad C_{\g,\pm} = \pm C_U \cap \g_{\pm 1}(h).
\end{equation}

\begin{defn} For an antiunitary representation 
$(U,\cH)$ of $G_{\tau_h}$ and $J := U(\tau_h^G)$, we write 
$\cH^{-\infty}_{{\rm ext},J}$ for the real linear space of all 
distribution vectors $\eta \in \cH^{-\infty}$, for which the orbit map 
\[ \alpha^\eta \: \R \to \cH^{-\infty}, \quad 
\alpha^\eta(t) := \eta \circ U(\exp -th) \] 
extends to a map on the closed strip $\oline{\cS_\pi}$, 
which is continuous with respect to the weak-$*$-topology on 
the space of distribution vectors, weak-$*$-holomorphic 
on the open strip  $\cS_\pi$ and satisfies 
\begin{equation}
  \label{eq:kms-cond-eta}
 \alpha^\eta(\pi i) = J_\sV \eta 
\end{equation}
(cf.~\cite[Def.~3.6]{NO21a}). 
\end{defn}

\begin{thm} \mlabel{thm:3.6} {\rm(\cite[Thm.~3.5, Prop.~3.13]{NO21a})}
Let $\sE \subeq \cH^{-\infty}$ 
be a real subspace invariant 
under $U(\exp(\R h))$ such that $\sH_\sE(G)$ is total in $\cH$. 
Then the following assertions hold: 
\begin{itemize}
\item[\rm(a)] If $\eset\not=\cO\subeq G$ is open, 
then $\sH_\sE(\cO)$ is total in $\cH$ {\rm(Reeh--Schlieder property)}. 
\item[\rm(b)] If $\sE \subeq \cH^{-\infty}_{{\rm ext},J}$, 
then $\sH_\sE(S_\sV^0) = \sV$  
is the standard subspace from \eqref{eq:bgl2}. 
\end{itemize}
\end{thm}

We shall use the following consequence of this theorem: 

\begin{cor} \mlabel{cor:app.d} 
For any non-empty open subset 
$\eset \not= \cO \subeq S_{\sV}$ with $\exp(\R h) \cO = \cO$, we have 
\[ \sH_\sE(\cO) = \sV.\] 
In particular, $\sH_\sE(S_{\sV,e}^0) = \sV$. 
\end{cor}

\begin{prf} From Theorem~\ref{thm:3.6} we obtain $\sH_\sE(\cO) 
\subeq \sH_\sE(S_\sV^0) = \sV$, 
so that 
\[ \sH_\sE(\cO) \cap i \sH_\sE(\cO) \subeq \sV \cap i \sV = \{0\}.\]
Further, Theorem~\ref{thm:3.6}(a) implies that $\sH_\sE(\cO)$ is total, hence 
standard. Now the invariance of $\sH_\sE(\cO)$ under the modular 
group $U(\exp \R h)$ of $\sV$ shows equality (\cite[Lemma~3.4]{NO21a}). 
\end{prf}

\begin{cor} \mlabel{cor:7.4} 
Let $P \subeq G$ be a closed subgroup with 
$h \in \L(P)$ and 
$\eta \in \cH^{-\infty}$ fixed by $U(P)$ and $J$ which is cyclic in the sense 
that $\sH_{\R \eta}(G)$ is total in $\cH$. Then 
\[ \sH_{\R \eta}(S_{\sV,e}^0) = \sH^{G/P}_{\R \eta}(S_{\sV,e}^0 P/P) =\sV .\]
\end{cor}

\begin{prf} We put $\sE := \R \eta$. As $\eta$ is fixed by $P$, 
it is $U(P)$-invariant, hence in particular invariant under 
$U(\exp \R h)$. As $J\eta = \eta$, we have 
$\sE \subeq \cH^{-\infty}_{{\rm ext}, J}$, so that 
Corollary~\ref{cor:app.d} 
shows that $\sH_\sE(S_{\sV,e}^0) = \sV$. For the second assertion we 
use   \cite[Lemma~2.11]{NO21a} to obtain 
\[ \sV 
= \sH_\sE(S_{\sV,e}^0) 
= \sH_\sE(S_{\sV,e}^0 P) 
= \sH_\sE^{G/P}(S_{\sV,e}^0 P/P).\qedhere \] 
\end{prf}

The preceding corollary shows that the orbit 
$S_{\sV,e}^0 P/P = S_{\sV,e}^0.eP$ of the base point $eP \in G/P$ 
under the open subsemigroup $S_{\sV,e}^0$ is a natural domain 
to which the standard subspace can be associated. 

For compactly causal symmetric spaces, we thus obtain the following 
concretization: 

\begin{thm} \mlabel{thm:7.5} 
Let $M = G/H$ be a compactly causal modular symmetric space with 
infinitesimal data $(\g,\tau, C,h)$. 
Further, let $(U,\cH)$ be an anti-unitary representation of $G_{\tau_h}$ 
whose positive cone $C_U$ is pointed, 
$\g = C_U - C_U + \R h$, and $C = C_U \cap \fq$. 
Let $\eta \in \cH^{-\infty}$ be fixed by $U(H)$ and $J = U(\tau_h^G)$ 
such that $\sH_{\R \eta}(G)$ is total in $\cH$. Then 
\[  \sH^{G/H}_{\R \eta}(W_M(h)_{eH}) = \sV.\] 
\end{thm}

\begin{prf} In view of Corollary~\ref{cor:7.4}, 
it suffices to show that $S_{\sV,e}^0.eH = W_M(h)_{eH}$. 
As $C_U$ is pointed, the kernel of $U$ is discrete. 
Therefore \eqref{eq:sv2} and 
$G^h_e \subeq G^{h,\tau_h^G} \subeq G_\sV$ imply 
\[ S_{\sV,e} = G^h_e \exp(C_U^c).\] 
Now Proposition~\ref{prop:5.2}(b) leads to 
\[ S_{\sV,e}^0.eH 
= (G^h)_e.\Exp_{eH}(C^{c,0})  = W_M(h)_{eH}. \]
This completes the proof. 
\end{prf}

\section{Nets of standard subspaces 
on compactly causal symmetric spaces} 
\mlabel{sec:6}

The main result of the preceding section (Theorem~\ref{thm:7.5}) describes 
how certain 
antiunitary representations of $G_{\tau_h}$ lead to interesting nets of 
standard subspaces on compactly causal symmetric spaces. 
These representations where supposed to have a positive cone $C_U$ 
which is pointed and satisfies $\g = C_U - C_U + \R h$,
and a cyclic distribution vector in 
$\cH^{-\infty}$, fixed by $U(H)$ and $J= U(\tau^G_h)$. 
In this section we construct such representations explicitly 
in spaces of Hilbert--Schmidt operators $\cH_\rho \subeq B_2(\cK)$, where 
$(\rho, \cK)$ is an antiunitary representation of $G_{\tau_h}$ which is a 
finite sum of irreducible representations. This is done in three steps: 
First we recall from \cite{Ne00, Ne19} some results on 
the representations $(\rho,\cK)$ of $G_{\tau_h}$, 
then we use these representations to construct nets of standard subspaces 
on the symmetric space of group type $G \cong (G \times G)/\Delta_G$, 
and finally we use the twisted embedding 
$G \to G \times G, g \mapsto (g, \tau^G(g))$ to obtain 
pullback representations of $G_{\tau_h}$ that can be used to obtain with 
Theorem~\ref{thm:7.5} nets of standard subspaces on $M = G/H$.

In this section we fix the following notation: $\g$ is a 
{\it semisimple} Lie algebra 
and $C_\g \subeq \g$ is a pointed generating invariant cone. 
Further, $G$ is a connected Lie group with Lie algebra~$\g$, 
$\tau^G$~an involutive automorphism of $G$, and 
the corresponding automorphism 
$\tau = \L(\tau^G)$ of $\g$ satisfies $-\tau(C_\g) = C_\g$. 
This implies that $C_\fq := C_\g \cap \fq$, $\fq = \g^{-\tau}$,  
has interior points, so that $(\g,\tau, C_\fq)$ is a compactly causal symmetric Lie algebra. 
We also fix an Euler element $h \in \fh = \g^\tau$ 
and assume that the Lie algebra involution 
$\tau_h = e^{\pi i \ad h}$ integrates to an involutive automorphism 
$\tau_h^G$ of~$G$ and that $-\tau_h(C_\g) = C_\g$ and 
 $\tau_h^G(H) = H$. 

\subsection{$C_\g$-positive representations and 
standard subspaces} 
\mlabel{subsec:6.1}

To prepare our construction of nets of standard subspaces on 
compactly causal symmetric spaces, we first collect some information on 
$C_\g$-positive representations of $G$. 
We thus obtain an interface to \cite{Ne19},  \cite{NOO21}  and \cite{NO21a} from which some 
results and constructions will be used below.  

Let $(\rho, \cK)$ be an antiunitary $C_\g$-positive representation 
of $G_{\tau_h} := G \rtimes \{\id_G, \tau_h^G\}$ and write 
$J_\cK := \rho(\tau_h^G)$ for the corresponding conjugation on $\cK$. 
Then there exists a unique standard subspace $\sV_\cK \subeq \cK$ with 
\[ J_{\sV_\cK} = J_\cK \quad \mbox{ and } \quad 
\Delta_{\sV_\cK}^{-it/2\pi} = \rho(\exp th)\quad \mbox{ for }  \quad t \in \R.\] 
Let
\[G_{\sV_\cK} = \{g\in G\: \rho (g) \sV_\cK = \sV\}. \]
Assume that $\rho$ has a discrete kernel. Then 
the derived representation $\dd\rho$ is injective, so that 
$\rho(G_{\sV_\cK})$ commutes with $\dd\rho(h)$ and thus 
$G_{\sV_\cK} \subeq G^h$. More specifically,  
\begin{equation}
  \label{eq:sc3}
 G_{\sV_\cK} 
= \{ g \in G^h \: \rho(g) J_\cK \rho(g)^{-1} = J_\cK\} 
= \{ g \in G^h \: g \tau_h^G(g)^{-1} \in \ker(\rho)\}. 
\end{equation}
This is an open subgroup of $G^h$ that only depends on the 
kernel of $\rho$. 

Furthermore, the discreteness of the kernel of $\rho$ implies
that its positive cone 
$C_\rho \subeq \g$ is pointed. It is also generating because it contains 
$C_\g$. Now \cite[Thms.~2.16, 3.4]{Ne19} imply that 
\begin{equation}
  \label{eq:sc1}
S_{\sV_\cK} 
:= \{ g \in G \: \rho(g) \sV_\cK \subeq \sV_\cK \} 
= G_{\sV_\cK} \exp(C_\rho^c),
\end{equation}
where 
\[ C_\rho^c := C_{\rho,+} + C_{\rho,-} \quad \text{with}\quad 
C_{\rho,\pm} := \pm C_\rho \cap \g_{\pm 1}(h).\] 

Note that $C_\rho^c$ is a $G_{\sV_\cK}$-invariant hyperbolic
convex cone in $\g^{-\tau_h}$.
Since $C_\rho^c$ coincides with $C_\g^c$ by 
Corollary~\ref{cor:red-minmax}, and 
$-\tau_h(C_\rho) = C_\rho$, the interior of this cone 
intersects~$\g^{-\tau_h}$, so that 
$C_\rho^c$ generates $\g^{-\tau_h}$. Therefore 
\begin{equation}
  \label{eq:sc2}
 S_{\sV_\cK} = G_{\sV_\cK} \exp(C_\g^c),
\end{equation}
and this semigroup only depends on $\rho$ through its unit group~$G_{\sV_\cK}$. 
We conclude that the semigroup 
$S_{V_\cK}$ only depends on the representation 
$\rho$ through the discrete subgroup $\ker(\rho)$.

For a more detailed discussion of irreducible antiunitary representations 
of $G_{\tau_h}$, we refer to Proposition~4.6 and Remark~4.10 in 
\cite{NO21a}. They restrict to representations of $G$ which are 
either irreducible or a direct sum of two inequivalent representations 
exchanged by twisting with $\tau_h^G$. 

For $1\le p  <\infty$ we denote by $B_p (\cH )$ the space of linear operators with finite $p$-Schatten norm:
\[B_p (\cH ) = \{T\in B(\cH )\: \|T\|_p = \tr (|T_p|)^{1/p}<\infty\}.\]
Then $B_1 (\cH)$ is the space of trace class operators and $B_2 (\cH)$ is the space of Hilbert--Schmidt operators.

\begin{prop}  \mlabel{prop:8.1}
There exists an injective antiunitary representation 
$(\rho,\cK)$ of $G_{\tau_h}$ such that 
\begin{itemize}
\item[\rm(a)] $C_\rho \supeq C_\g$ 
\item[\rm(b)] $\rho$ is a  finite direct sum 
of irreducible representations. 
\item[\rm(c)] $\rho$ is a trace-class representation, 
i.e., $\rho(C^\infty_c(G)) \subeq B_1(\cH)$. 
\end{itemize}
\end{prop}

\begin{prf} From \cite[Thm.~XI.5.2]{Ne00} we know that 
there exist finitely many irreducible representations 
$(\pi_j, \cK_j)_{j \in F}$ of $G$ with $C_{\pi_j} \supeq C_\g$ and 
$\bigcap_{j \in F} \ker(\pi_j) = \{e\}$. 
Then $\pi := \oplus_{j \in F} \pi_j$ is injective and 
\[ \rho := \pi \oplus (\pi^* \circ \tau_h^G) \] 
 extends to an injective antiunitary representation of $G_{\tau_h}$. As 
$-\tau_h(C_\g) = C_\g$, we also have $C_\g \subeq C_\rho$. 
This proves (a) and (b). 
That $\rho$ is a trace class representation follows from 
the fact that all representations $\pi_j$  
are trace class (\cite[Thm.~X.4.10]{Ne00}).
\end{prf}

\subsection{The group case: biinvariant nets on Lie groups} 
\mlabel{subsec:6.2}

On $G$ we obtain by  the assignment
\[g\mapsto V_+(g) := g.C_\g^0 \subeq T_g(G) \quad \mbox{ for } \quad g \in G,\] 
a biinvariant cone field, turning $G$ into a compactly causal 
symmetric space on which the group $G \times G$ acts by left and 
right translations. Its infinitesimal 
data is represented by the compactly causal symmetric Lie algebra 
$(\g \oplus \g, \tau , C_{\g \oplus \g})$ of group type, where 
\begin{equation}
  \label{eq:tildecg}
\tau(x,y) = (y,x)\quad\text{and} \quad C= \{ (y,-y) \: y \in C_\g \}\subset (\g \oplus \g)^{-\tau}
\end{equation}
We also fix an Euler element $h \in \cE(\g)$.  

We want to apply \cite{NO21a} to representations of $G \times G$, 
so we assume that $\g$ is {\it semisimple}.\begin{footnote}
{ For generalization of some constructions in \cite{NO21a} 
to non-reductive Lie algebras, we refer to \cite{Oeh21}.}
\end{footnote}

To this end, we consider in $\g \oplus \g$ the 
invariant cone 
\[ C_{\g \oplus \g} := C_\g \oplus - C _\g \quad \mbox{ which satisfies  } \quad 
-\tau(C_{\g \oplus \g}) = C_{\g \oplus \g}.\]

\begin{rem} (Irreducible representations) 
\mlabel{rem:6.2} 
As $C_\g$-positive representations of $G$ 
are type I (\cite[Thm.~X.6.21]{Ne00}), 
irreducible $C_\g$-positive representations of $G \times G$ have the form 
\[ \pi = \rho_1 \boxtimes \rho_2^*, \quad 
\pi(g_1,g_2) = \rho_1(g_1) \otimes \rho_2^*(g_2)\] 
where $(\rho_1, \cK_1)$ and $(\rho_2, \cK_2)$ 
are irreducible $C_\g$-positive representations of $G$. 

For the existence of an antiunitary extension to 
\[ (G \times G)_{\tau_h \times \tau_h}
= (G \times G) \rtimes \{\id, \tau_h^G \times \tau_h^G\},\] 
we need that 
\begin{equation}
  \label{eq:33}
 \pi^* \cong \pi \circ (\tau_h^G \times \tau_h^G) 
= (\rho_1 \circ \tau_h^G) \boxtimes (\rho_2 \circ \tau_h^G)^*, 
\end{equation}
and this is equivalent to 
\begin{equation}
  \label{eq:6.1}
\rho_j^* \cong \rho_j \circ \tau_h^G \quad \mbox{ for } \quad j =1,2.
\end{equation}
If these two conditions are satisfied, then $\rho_1$ and $\rho_2$ extend 
to antiunitary representations of $G_{\tau_h}$ and $\pi$ extends 
to $(G \times G)_{\tau_h \times \tau_h}$ by 
\[ \pi(\tau_h^G \times \tau_h^G) := 
\rho_1(\tau_h^G) \otimes \rho_2(\tau_h^G).\] 
If \eqref{eq:6.1} is not satisfied, then 
the canonical antiunitary extension of the representation 
\[ \pi \oplus (\pi^* \circ (\tau_h^G \times \tau_h^G)) \] 
is irreducible, 
although its restriction to $G \times G$ decomposes into two 
inequivalent irreducible constituents. 
\end{rem}

Let $(\rho_j, \cK_j)$, $j =1,2$, be two antiunitary 
$C_\g$-positive representations of $G_{\tau_h}$ which need not 
be  irreducible (cf.~Subsection~\ref{subsec:6.1}). 
Then 
\[ \cH := B_2(\cK_2, \cK_1) \cong \cK_1 \otimes \cK_2^* \] 
 carries 
the unitary representation $(\cH,\pi)$ of $G \times G$ defined by 
\[ \pi(g_1, g_2) A  = \rho_1(g_1) A \rho_2(g_2)^{-1},\] 
and we extend this representation to an antiunitary 
representation of $(G \times G)_{\tau_h \times \tau_h}$ by 
\[ \pi(\tau_h^G \times \tau_h^G) := J, \quad 
J(A) :=  J_{\cK_1} A J_{\cK_2}, \quad \mbox{ where } \quad 
J_{\cK_j} := \rho_j(\tau_h^G).\] 

\begin{lem} Let $\sV_j \subeq \cH_j$, $j =1,2$, be standard subspaces 
with the modular objects   $(\Delta_j, J_j )$, $j =1,2$  where
$\Delta_j=\Delta_{\sV_j}, J_j=J_{\sV_j}$. Then 
\[  \sV_1 \otimes \sV_2 \subeq \cH_1 \otimes \cH_2 \] 
is the standard subspace $\sV$ corresponding to the modular group defined by 
\begin{equation}
  \label{eq:motens}
\Delta^{it} = \Delta_1^{it} \otimes \Delta_2^{it}\quad \mbox{ for } \quad 
t \in \R, \quad \mbox{ and the conjugation } \quad J = J_1 \otimes J_2.
\end{equation}
\end{lem}

\begin{prf} The tensor product $\sH := \sV_1 \otimes \sV_2$ is a closed 
real subspace of $\cH = \cH_1 \otimes \cH_2$ which we want to compare 
with the standard subspace $\sV := \Fix(J \Delta^{1/2})$ defined by 
the pair $(\Delta, J)$ from~\eqref{eq:motens}. 
For $\xi_1 \in \sV_1$ and $\xi_2 \in \sV_2$, we have 
\[ \Delta_1^{1/2} \xi_1 = J_1 \xi_1 \quad \mbox{ and } \quad 
  \Delta_2^{1/2} \xi_2 = J_2 \xi_2.\] 
For $\xi := \xi_1 \otimes \xi_2$ the orbit map 
\[ \R \to \cH, \quad t \mapsto \Delta^{-it/2\pi}\xi 
= \Delta_1^{-it/2\pi} \xi_1 \otimes \Delta_2^{-it/2\pi} \xi_2 \] 
extends holomorphically to the closure of the strip 
$\cS_\pi \subeq \C$ by 
\[ z \mapsto \Delta_1^{-iz/2\pi} \xi_1 \otimes 
\Delta_2^{-i z/2\pi} \xi_2 \] 
(cf.\ \cite[Prop.~2.1]{NOO21}). This implies that 
\[ \Delta^{1/2}\xi 
= \Delta_1^{1/2} \xi_1 \otimes \Delta_2^{1/2} \xi_2 
= J_1 \xi_1 \otimes J_2 \xi_2 = J \xi.\] 
Therefore $\xi \in \sV$ (\cite[Prop.~2.1]{NOO21}). This argument shows that 
$\sH \subeq \sV$. As $\sH$ is generating in $\cH$ 
and invariant under the modular group 
$\Delta^{-it/2\pi} 
= \Delta_1^{-it/2\pi}  \otimes \Delta_2^{-it/2\pi}$,  
\cite[Cor.~2.1.8]{Lo08} implies that $\sH = \sV$. 
\end{prf}

\begin{cor} Write $\cK_2^{\rm op} \cong \cK_2^*$ for $\cK_2$, endowed 
with the opposite complex structure. 
Then the standard subspace $\sV \subeq \cK_1 \otimes \cK_2^{\rm op}$ 
specified by the Lie algebra element $(h,h) \in \g \oplus\g$ and $J$  
is the tensor product $\sV = \sV_1 \otimes \sV_2'$, 
constructed from the corresponding standard subspaces $\sV_j \subeq \cK_j$. 
\end{cor}

\begin{prf} In view of the preceding lemma, we only have to observe 
that, changing the complex structure on $\cK_2$, corresponds to replacing 
the modular operator $\Delta_{\sV_2} = e^{2\pi i \partial \rho_2(h)}$ by its inverse. 
\end{prf}

We now assume that 
$\rho := \rho_1 = \rho_2$ is a  sum of {\it finitely many irreducible} 
$C_\g$-positive antiunitary representations of~$G_{\tau_h}$. 
Then 
\[ \pi(g_1, g_2) A  = \rho(g_1) A \rho(g_2)^{-1} \quad \mbox{ and } \quad 
J(A) =  J_\cK A J_\cK. \] 
By \cite[Thm.~X.4.10]{Ne00} 
$\rho$ is a trace class representation, i.e., 
$\rho(C^\infty_c(G)) \subeq B_1(\cH)$. 
This implies that the integrated representation 
defines a continuous linear map $\rho \: C_c^\infty(G) \to B_1(\cK)$  
(\cite[Thm.~1.3]{DNSZ16}). 

\begin{lem} \mlabel{lem:b2rep}
The following assertions hold: 
  \begin{itemize}
  \item[\rm(a)] For the left multiplication 
representation of $G$ on $B_2(\cK)$, 
\begin{equation}
  \label{eq:leftsmooth}
 B_2(\cK)^\infty = \{ A \in B(\cK) \: A\cK \subeq \cK^\infty\} 
\subeq B_1(\cK) 
\end{equation}
and the inclusion $B_2(\cK)^\infty \to B_1(\cK)$
is continuous with respect to the norm topology on $B_1(\cK)$ 
and the natural Fr\'echet topology on $B_2(\cK)^\infty$. 
  \item[\rm(b)] We have a $(G \times G)$-equivariant map 
\[  \rho^* \: B(\cK) \to C^{-\infty}(G), \quad \rho^*(A)(\xi) := \tr(\rho(\xi)^*A),\quad  \xi \in C_c^\infty (G).\] 
  \item[\rm(c)] The space of smooth vectors for the 
representation $\pi$ of $G \times G$ on $B_2(\cK)$ is 
\[ B_2(\cK)^{\infty,\pi}
= \{ A \in B(\cK) \: 
A\cK \subeq \cK^\infty, 
A^*\cK \subeq \cK^\infty \}. \] 
\item[\rm(d)] The map 
$B(\cK) \to B_2(\cK)^{-\infty,\pi}, A \mapsto \eta_A$, 
with $\eta_A(B) = \tr(B^*A)$ is  injective 
and equivariant with respect to the representation 
of $(G \times G)_{\tau_h \times \tau_h}$ on $B(\cH)$, given by 
\[(g_1,g_2).A := \rho(g_1) A \rho(g_2)^{-1}\quad \mbox{ for } 
\quad g_1, g_2 \in G 
\quad \mbox{ and } \quad (\tau_h^G \times \tau_h^G).A = J_\cK A J_\cK.\] 
  \end{itemize}
\end{lem}

Note that (a) implies in particular that 
each bounded operator on $\cK$ defines a distribution 
vector for the left multiplication representation on $B_2(\cK)$, 
hence also for the representation $\pi$ of $G \times G$. 
The smooth vectors for $G \times G$ 
coincide with the algebra of so-called {\it Schwartz operators}, 
i.e., operators which remain bounded when composed on the left and 
the right with elements of the enveloping algebra acting on smooth vectors 
(cf.~\cite{DNSZ16}, \cite{KKW16}).

\begin{prf} (a) The first equality in \eqref{eq:leftsmooth} 
follows from \cite[Prop.~1.7]{DNSZ16}. 

Let $A \in B_2(\cK)$ be a smooth vector for the 
left multiplication representation of~$G$. 
By the Dixmier--Malliavin Theorem, $A$ is a finite sum 
of operators of the form $\rho(\xi_j)A_j$, $A_j \in B_2(\cK)$, 
$\xi_j \in C^\infty_c(G)$, and these are trace class because
$(\rho,\cK)$ is a trace class representation.

For the last assertion, we have to show that the inclusion map 
\[ \Gamma \: B_2(\cK)^\infty \to B_1(\cK) \] 
is continuous. This is a linear map from a Fr\'echet space to a Banach space. 
Therefore it suffices to write it as a pointwise limit of a 
sequence $(\Gamma_n)_{n \in \N}$ of continuous linear maps 
(\cite[Lemma~1.2]{DNSZ16}). 
So let $\xi_n\in C^\infty_c(G)$ be a $\delta$-sequence 
and $\Gamma_n(B) :=  \rho(\xi_n) B$. 
These maps are continuous from $B_2(\cK) \to B_1(\cK)$, 
hence in particular from $B_2(\cK)^\infty \to B_1(\cK)$. 
Since the left multiplication representation of $G$ on 
the Banach space $B_1(\cK)$ is continuous, 
$\|\Gamma_n(B)-B\|_1 \to 0$  holds for every $B \in B_1(\cH) 
\supeq B_2(\cH)^\infty$. 

\nin (b) From the continuous linear map 
$\rho \: C^\infty_c(G) \to B_1(\cK)$, we obtain the map $\rho^*$ 
by taking adjoints. The equivariance for the action of 
$G \times G$ on both sides follows from 
\begin{align*}
 \rho^*\big(\rho(g_1) A \rho(g_2)^{-1}\big)(\xi) 
&= \tr\big(\rho(\xi)^*\rho(g_1) A \rho(g_2)^{-1}\big)
= \tr\big(\big(\rho(g_1)^{-1}\rho(\xi)\rho(g_2)\big)^* A\big)\\
&= ((g_1,g_2).\rho^*(A))(\xi).
\end{align*}

\nin (c) follows from \cite[Cor.~1.8]{DNSZ16}. 

\nin (d) That every  functional $\eta_A$ defines a distribution 
vector follows from (a). For the injectivity of the assignment, 
we first observe that  
 all rank one-operators $P_{\xi,\eta} = |\xi\ra \la \eta|$,  
defined by smooth vectors $\xi,\eta \in \cK^\infty$, 
are contained in $B_2(\cK)^\infty$. 
For these we have 
\[ \eta_A(P_{\xi,\eta})    
= \tr(P_{\xi,\eta}^* A)= \tr(P_{\eta,\xi} A) = \la \xi, A \eta\ra, \] 
and since $\cK^\infty$ is dense in $\cK$, any bounded operator $A$ is determined 
by $\eta_A$.

For the equivariance, we calculate 
\[ \eta_A(\pi(g_1,g_2)^{-1}B) 
= \tr((\rho(g_1)^* B \rho(g_2))^*A) 
= \tr( \rho(g_2)^*B^* \rho(g_1) A) 
= \eta_{(g_1,g_2).A}(B).\] 
Further 
\[ \eta_{J_\cK A J_\cK}(B) 
= \tr(B^* J_\cK A J_\cK)
= \oline{\tr(J_\cK B^* J_\cK A)}
= \oline{\tr( (J_\cK B J_\cK)^* A)}
= \pi^{-\infty}(\tau_h^G \times \tau_h^G)(\eta_A)(B).\qedhere\] 
\end{prf}

We want exhibit concrete contexts, where the assumptions of 
Theorem~\ref{thm:7.5} are satisfied. This will be achieved in Theorem~\ref{thm:6.3} 
and the ground is prepared by the following proposition. 
Recall that we assume in this section that $\g$ is  semisimple. 
In the following proposition this is crucial to use 
Corollary~\ref{cor:red-minmax} in the proof (d).

\begin{prop} \mlabel{prop:6.6} Assume that $\ker \rho$ is discrete. Let 
$\cH_\rho \subeq \cH = B_2(\cK)$ denote the closed subspace generated by 
the image $\rho(C^\infty_c(G))$ of the integrated representation of the 
convolution algebra of test functions and denote the corresponding 
subrepresentation of $\pi$ by $\pi_\rho$. Then the following 
assertions hold: 
\begin{itemize}
\item[\rm(a)] $\eta^0 := \tr\res_{\cH_\rho^\infty} \in \cH_\rho^{-\infty}$ 
is a cyclic distribution vector invariant under the diagonal subgroup \break 
$\Delta_G = \{ (g,g) \:  g \in G \}$  and $J(A) = J_\cK A J_\cK$. 
\item[\rm(b)] $\cH_\rho$ is $J$-invariant and 
$\sE := \R \eta^0 \subeq (\cH^{-\infty}_\rho)_{{\rm ext},J}$. 
\item[\rm(c)] $\ker(\pi_\rho)$ is discrete. 
\item[\rm(d)] Let $\sV_\rho \subeq \cH_\rho$ be the standard subspace 
with the modular data 
$J_{\sV_\rho} = J\res_{\cH_\rho}$ and 
$\Delta_{\sV_\rho} = e^{2\pi i \cdot \partial \pi(h,h)}$. 
Then 
\[ S_{\sV_\rho} = \{ g \in G \times G \: \pi_\rho(g) \sV_\rho \subeq \sV_\rho\} 
= (G \times G)_{\sV_\rho} \exp(C_\g^c \oplus -C_\g^c).\] 
\end{itemize}
\end{prop}

\begin{prf} (a), (b) By Lemma~\ref{lem:b2rep}(d), the trace $\eta_\1 = \tr 
\in \cH^{-\infty} = B_2(\cK)^{-\infty}$ is 
invariant under $J$ and the diagonal subgroup~$\Delta_G$. 
In particular, it is fixed under the modular group 
$\exp(\R(h,h))$, which implies (b). 

If we consider $\1 \in B(\cH)$ as a distribution vector, we have 
for $\xi_1, \xi_2 \in C^\infty_c(G)$ the relation 
 
\begin{align}
  \label{eq:74}
\pi^{-\infty}(\xi_1 \otimes \xi_2) \1 
&= \rho(\xi_1) \int_G \xi_2(g) \rho(g^{-1})\, dg 
= \rho(\xi_1) \int_G \xi_2(g^{-1}) \rho(g)\, dg\nonumber \\
&= \rho(\xi_1) \rho(\xi_2^\vee) 
= \rho(\xi_1 * \xi_2^\vee).
\end{align} 
As tensor products of test functions on $G$ span a dense subspace 
of $C^\infty_c(G \times G)$, 
this calculation shows that the closed subspace of $\cH$ 
generated by $\pi^{-\infty}(C^\infty_c(G \times G)) \sE$ coincides with 
$\cH_\rho$. This subspace is also $J$-invariant because 
 
\begin{align*} J \rho(\xi) J 
&= \int_G \oline{\xi(g_1, g_2)} \pi(\tau_h^G(g_1), \tau_h^G(g_2))\, dg_1 \, dg_2
= \int_G \oline{\xi(\tau_h^G(g_1), \tau_h^G(g_2))} \pi(g_1, g_2)
\, dg_1 \, dg_2.
\end{align*}

\nin (c) As $\cH_\rho$ is generated by $\sE$, the kernel of $\pi_\rho$ is 
the largest normal subgroup of $G \times G$ acting trivially on~$\sE$. 
A pair $(g_1, g_2) \in G \times G$ fixes $\1 \in B(\cK)$ if and only if 
$\rho(g_1 g_2^{-1}) = \1$, i.e., $g_1 g_2^{-1} \in \ker(\rho)$. 
If $(g_1,g_2) \in \ker(\pi_\rho)$, then we have for all $a,b \in G$ the relation 
\[ ag_1 a^{-1} b g_2^{-1} b^{-1} \in \ker(\rho), \] 
and since $G$ is connected and $\ker(\rho)$ discrete, it follows that 
$g_1, g_2 \in Z(G)$. This shows that 
\begin{equation}
  \label{eq:kerrho}
\ker(\pi_\rho) 
= \{ (g_1, g_2) \in Z(G) \times Z(G) \: g_1 g_2^{-1} \in \ker(\rho)\} 
= \Delta_{Z(G)}\cdot (\ker(\rho) \times \{e\}),
\end{equation}
which is a discrete subgroup of $G \times G$. 
%

\nin (d) As the kernel of $\pi_\rho$ is discrete by (c), the description 
of the semigroup $S_{\sV_\rho}$ follows from 
\cite[Thms.~2.16, 3.4]{Ne19}. 
Note that the discreteness of $\ker(\pi_\rho)$ implies that 
$C_{\pi_\rho}^c = C_{\g \oplus \g}^c = C_\g^c \oplus - C_\g^c$ 
(Corollary~\ref{cor:red-minmax}). 
\end{prf}

The following theorem shows that, for symmetric spaces of group 
type, there exist biinvariant nets of real subspaces assigning 
to the wedge domain $W_G(h)_e = G^h_e \exp(C_\g^{c,0})$ the standard subspace 
associated to the pair $(h,\tau_h^G)$ by the BGL construction 
(Subsection~\ref{subsec:bgl}).

\begin{thm} \mlabel{thm:6.3}
Let $(\rho, \cK)$ be a $C_\g$-positive antiunitary trace class representation 
of $G_{\tau_h}$ with discrete kernel 
and consider the representation $(\pi_\rho, \cH_\rho)$ 
of $G \times G$  on the closed subspace 
\[ \cH_\rho 
:= \oline{\rho(C^\infty_c(G))} \subeq B_2(\cK), \quad \mbox{ given by } \quad 
\pi_\rho(g_1,g_2)A = \rho(g_1) A \rho(g_2)^{-1},\] 
extended to an antiunitary representation of $(G \times G)_{\tau_h \times \tau_h}$ 
by 
\[ J(A) := \pi_\rho(\tau_h^G \times \tau_h^G)(A) := \rho(\tau_h^G)A \rho(\tau_h^G). \]
For $\sE := \R \eta^0$,  $\eta^0(A) := \tr(A)$, and 
$q \: G \times G \to G,q(g_1, g_2) := g_1 g_2^{-1}$, 
\begin{equation}
  \label{eq:biinv-net}
 \sH^G_\sE(\cO) := \sH_\sE^{G \times G}(q^{-1}(\cO)) 
= \oline{\rho(C^\infty_c(q^{-1}(\cO),\R))}\subeq \cH_\rho 
\end{equation}
defines a $(G \times G)$-covariant net of closed real subspaces in $\cH_\rho$ on 
$G \cong (G \times G)/\Delta_G$. The wedge domain 
\[ W_G(h)_e = G^h_e \exp(C_\g^{c,0}) 
\] 
satisfies 
\begin{equation}
  \label{eq:wedge-stand}
 \sH^G_\sE(W_G(h)_e) = \sV_\rho 
= \oline{\Spann_\R \rho(C^\infty_c(W_G(h)_e,\R))}, 
 \end{equation}
where $\sV_\rho$ is the standard subspace corresponding to the pair 
$(e^{2\pi i \cdot \partial \pi_\rho(h,h)}, J)$ of modular data. 
\end{thm}

\begin{prf}  By Proposition~\ref{prop:6.6}(c), the kernel
of $\pi_\rho$ is discrete. Therefore the 
invariance of $\eta^0$ under the modular group and $J$ and  
Theorem~\ref{thm:3.6}(b) imply that 
\[ \sV_\rho = \sH^{G \times G}_\sE(S_{\sV_\rho}^0) := 
\oline{\pi_\rho^{-\infty}(C^\infty_c(S_{\sV_\rho}^0,\R)) \eta^0}  \] 
holds for the open subsemigroup 
\[  S_{\sV_\rho}^0 = (G \times G)_{\sV_\rho} \exp(C_\g^c \oplus - C_\g^c), 
\quad \mbox{ where } \quad 
C_\g^c = C_{\g,+} - C_{\g,-}\] 
(Proposition~\ref{prop:6.6}(d)). 
With Corollary~\ref{cor:app.d}, we even obtain 
\begin{equation}
  \label{eq:vrho}
\sV_\rho = \sH_\sE(S_{\sV_\rho,e}^0) \quad \mbox{ with }\quad 
S_{\sV_\rho,e} = (G^h_e \times G^h_e) \exp(C_\g^c \oplus - C_\g^c) 
= S(C_\g,h)_e \times S(C_\g,h)_e^{-1},
\end{equation}
where 
\[ S(C_\g,h)^0_e = S(C_\g^0,h)_e = W_G(h)_e \] 
(Theorem~\ref{thm:4.2}). 
With  
\[ \pi_\rho^{-\infty}(S(C_\g^0,h) \times S(C_\g^0,-h)) \eta^0= \rho(S(C_\g^0,h)),\] 
we further derive from \eqref{eq:vrho} 
\begin{equation}
  \label{eq:v-vis-srep}
\sV_\rho = \oline{ \rho(C^\infty_c(S(C_\g^0,h),\R))}. 
\end{equation}

As $\eta^0$ is $\Delta_G$-invariant (Proposition~\ref{prop:6.6}), \cite[Lemma~2.11(a)]{NO21a} shows that the net defined in 
\eqref{eq:biinv-net} satisfies 
\begin{equation}
  \label{eq:42}
 \sV_\rho 
= \sH_\sE^{G \times G}(S_{\sV_\rho,e}^0)
= \sH_\sE^{G \times G}(S_{\sV_\rho,e}^0 \Delta_G).
\end{equation}
Now $q^{-1}(W_G(h)_e) = (W_G(h)_e \times W_G(h)_e^{-1}) \Delta_G 
= S_{\sV_\rho,e}^0\Delta_G $ 
leads to 
  \begin{align*}
    \sH^G_\sE(W_G(h)_e) 
&= \sH_\sE^{G \times G}(q^{-1}(W_G(h)_e)) 
= \sH_\sE^{G \times G}\big(S_{\sV_\rho,e}^0\Delta_G\big) 
\ {\buildrel \eqref{eq:42}\over=}\ \sV_\rho. 
\qedhere\end{align*}
\end{prf}

\begin{rem}
For the stabilizer group of $\sV_\rho$, we have 
\[ (G \times G)_{\sV_\rho} 
= \{ (g_1,g_2) \in G^h \times G^h \: 
\pi_\rho(g_1, g_2) J = J \pi_\rho(g_1, g_2)\} 
\supeq G_{\sV_\cK} \times G_{\sV_\cK}.\] 
By \eqref{eq:kerrho}, 
the relation $\pi_\rho(g_1, g_2) J = J \pi_\rho(g_1, g_2)$ 
is equivalent to 
\begin{equation}
  \label{eq:37}
 (g_1 \tau_h^G(g_1)^{-1}, g_2 \tau_h^G(g_2)^{-1}) 
\in \ker(\pi_\rho) = \Delta_{Z(G)} (\ker(\rho) \times \{e\}) 
\subeq Z(G) \times Z(G).
\end{equation}
From Remark~\ref{rem:7.6new} we know that 
\[ \Ad(G)^h = \Ad(G^h) \subeq \Ad(G)^{\tau_h} \] 
is of index two if $\g$ is simple 
(otherwise of index $2^k$, where $k$ is the number of simple 
hermitian ideals of $\g$). Equation~\eqref{eq:37} implies in particular that 
\[  g_1, g_2 \in \Ad^{-1}(\Ad(G)^{\tau_h}) 
= \{ g \in G \: g \tau_h^G(g)^{-1} \in Z(G) \}.\] 
\end{rem}

\begin{lem} $\cH_\rho = B_2(\cK)$ if and only if $(\rho,\cK)$ is an 
irreducible representation of $G$.
\end{lem} 

\begin{prf} If $\cH_\rho = B_2(\cK)$, then $\rho(C^\infty_c(G))$ is dense 
in $B_2(\cK)$, so that the von Neumann algebra generated by this subspace 
is $B(\cK)$, which by Schur's Lemma implies that $(\rho,\cK)$ is irreducible. 

If, conversely, $(\rho,\cK)$ is irreducible, then 
$\rho(C^\infty_c(G)) \subeq B_1(\cK) \subeq B_2(\cK)$ 
is weakly dense in $B(\cK)$, 
i.e., its annihilator in $B_1(\cK)$ is trivial. 
Therefore the orthogonal space $\rho(C^\infty_c(G))^\bot \subeq B_2(\cK)$ 
has trivial intersection with $B_1(\cK)$. As it is invariant 
under multiplication with $\rho(C^\infty_c(G))$, which consists of trace class 
operators, this can only happen if $\rho(C^\infty_c(G))^\bot = \{0\}$, 
i.e., if $\cH = B_2(\cK)$. 
\end{prf}

\begin{rem} (a) If $\rho$ is reducible 
and extends to an irreducible antiunitary representation 
of $G_{\tau_h}$, then $G$ preserves a decomposition 
$\cK = \cK_1 \oplus \cK_2$ and $\rho(C_c^\infty(G))$ 
generates $B(\cK_1) \oplus B(\cK_2)$. Further 
$J = \rho(\tau_h^G)$ is a conjugation with 
$J\cK_1 = \cK_2$. Accordingly, $\cH_\rho \subeq B_2(\cK)$ 
is isomorphic to $B_2(\cK_1) \oplus B_2(\cK_2)$ (diagonal matrices). 

\nin (b) The trace class representation $(\rho,\cK)$ 
of $G$ on $\cK$ decomposes with finite multiplicities as a direct sum 
\[ (\rho,\cK) \cong \sum_{j \in J} (\rho_j, \cK_j)^{\oplus m_j}.\] 
Identifying $B(\cK_j)$ with diagonal operators in 
$M_{m_j}(B(\cK_j)) \cong B(\cK_j^{\oplus m_j})$, we then have 
\[ \rho(C^\infty_c(G)) \subeq \hat\bigoplus_{j \in J} B_2(\cK_j) \cong \cH_\rho.\] 
\end{rem}

\subsection{The general case: Nets on compactly causal symmetric spaces}
\mlabel{subsec:6.3}

Let $(\g,\tau, C)$ be a semisimple compactly causal symmetric Lie algebra 
for which $\fh = \g^\tau$ contains no non-zero ideal. 
We recall from \eqref{eq:grptypemb} the 
canonical embedding 
\begin{equation}\label{eq:Iota}
 \iota \: \g \to \g \oplus \g, \quad 
\iota(x) := (x,\tau (x) ) \quad \mbox{ for } \quad x \in\g 
\end{equation}
into a symmetric Lie algebra of group  type.  
The Extension Theorem~\ref{thm:extend} 
implies that  the elliptic cone $C \subeq \fq$ extends to a $-\tau$-invariant  
cone $C_\g \subeq \g$, 
so that 
\[ \tilde C_\g := \{ (y,-y) \: y \in C_\g\} \] 
leads to the causal symmetric Lie algebra 
 $(\g\oplus \g, \tau_{\rm flip}, \tilde C_\g)$ of group type. 
On the level of global symmetric spaces 
$M = G/H$, $H \subeq G^{\tau^G}$ open, the infinitesimal embedding 
$\iota$ corresponds to the quadratic representation 
\[ Q \: M \to M_G = G^\sharp_e\subeq G, \quad gH \mapsto gg^\sharp,\]   
which is a covering morphism of symmetric spaces. 

The following theorem is an important tool to 
analyze smooth vectors, resp., distribution vectors for subgroups. 

\begin{thm} {\rm(Zellner's Smooth Vector Theorem; \cite[Thm.~3.1]{NSZ17})} 
\mlabel{thm:zthm}
If $(\pi, \cH)$ is a unitary representation of a Lie group 
and $x_0 \in C_\pi^0$ in the interior of its positive cone, then the inclusion 
$\cH^\infty \to \cH^\infty(\partial \pi(x_0))$ is an isomorphism 
of Fr\'echet spaces. In particular every smooth vector for the single operator 
$\partial \pi(x_0)$ is smooth for $G$. 
\end{thm}

We now consider a $C_\g$-positive antiunitary trace class representation 
$(\rho,\cK)$ of $G_{\tau_h}$ with discrete kernel 
and use Theorem~\ref{thm:6.3} to obtain a 
representation $(\pi_\rho, \cH_\rho)$ of $(G \times G)_{\tau_h \times \tau_h}$ 
 on $\cH_\rho \subeq B_2(\cK)$ by 
\[ \pi_\rho(g_1,g_2)A = \rho(g_1) A \rho(g_2)^{-1} 
\quad \mbox{ and } \quad 
\pi_\rho(\tau_h \times \tau_h) A = J_\cK A J_\cK.\]

We introduce the morphism of graded Lie groups 
\[ \iota_G \: G_{\tau_h} \to (G \times G)_{\tau_h \times \tau_h}, \qquad 
g \mapsto (g, \tau^G(g)), \qquad 
\tau_h^G \mapsto \tau_h^G \times \tau_h^G, \] 
corresponding to the inclusion $\iota :\fg \to \fg\times \fg$ from \eqref{eq:Iota}.
We thus obtain the antiunitary representation of $G_{\tau_h}$, defined by 
\begin{equation}
  \label{eq:U-sym}
 U := \pi_\rho \circ \iota_G, \qquad 
U(g)A = \rho(g)A \rho(g)^\sharp 
\end{equation}

\begin{thm} \mlabel{thm:8.12}
The representations $U$ and $\pi_\rho$ on $\cH_\rho$ have the same 
smooth vectors. 
Let 
\[ \sE = \R \eta^0\subeq \cH^{-\infty}(\pi_\rho) = \cH^{-\infty}(U)
 \quad \mbox{ with } \quad \eta_0 = \tr \] be 
as in {\rm Proposition~\ref{prop:6.6}}.
If $q_M \: G \to M, g \mapsto gH$ denotes the canonical projection, then 
\begin{equation}
  \label{eq:netonM}
 \sH_\sE^M(\cO) 
:= \sH_\sE^{G,U}(q_M^{-1}(\cO)) := 
\oline{U^{-\infty}(C^\infty_c(q_M^{-1}(\cO),\R)) \sE}  
\end{equation}
defines a $G$-covariant net of closed real subspaces of $\cH_\rho$ on 
the symmetric space~$M = G/H$ which satisfies 
\begin{equation}
  \label{eq:standonM}
 \sV_{(h,\tau_h^G,U)} = \sH^M_\sE(W_M(h)_{eH}).
\end{equation} 
\end{thm}

\begin{prf} For any 
$x_0 \in C^0 \subeq \fq$ we have 
$\iota(x_0) = (x_0, -x_0) \in \tilde C_\g^0 \subeq C_\g^0 \oplus - C_\g^0$ 
(see \eqref{eq:tildecg}). 
Therefore Zellner's Theorem~\ref{thm:zthm} implies that 
smooth vectors of the single operator 
$\partial U(x_0) = \partial \pi_\rho(x_0,-x_0)$ are 
smooth for $G \times G$. This shows that 
\begin{equation}
  \label{eq:smovecdoub}
 \cH^\infty(\pi_\rho) = \cH^\infty(U) \quad \mbox{  and therefore } \quad 
\cH^{-\infty}(\pi_\rho) = \cH^{-\infty}(U).
\end{equation}
We conclude in particular that 
$\sE \subeq \cH^{-\infty}(U)$, so that 
\eqref{eq:netonM} defines a $G$-covariant net 
of closed real subspaces of $\cH_\rho$ on the symmetric space~$M = G/H$. 

The positive cone 
$C_U$ of $U$ contains  $\iota^{-1}(C_{\pi_\rho}) \supeq \iota^{-1}(C_\g \oplus - C_\g)$, 
hence is a closed generating 
invariant cone whose intersection with $\fq$ contains~$C$. 
As $\iota(h) = (h,h)$, the standard subspace $\sV_\rho \subeq \cH_\rho$ 
from Theorem~\ref{thm:6.3} 
coincides with the standard subspace $\sV_U := \sV_{(h,\tau_h^G,U)}$ 
associated to the pair $(h,\tau_h^G)$ and the 
representation $U$ of $G$ by the BGL construction. 
  
The discreteness of the kernel of $\pi_\rho$ (Proposition~\ref{prop:6.6}) 
implies that 
$U$ has discrete kernel. Therefore \cite[Thms.~2.16, 3.4]{Ne19} implies that 
\begin{equation}
  \label{eq:svrho}
 S_{\sV_\rho} = (G^h)_{\sV_\rho} \exp(C_U^c).
\end{equation}
We claim that 
\begin{equation}
  \label{eq:conesequal}
 C_U^c = C_\g^c.
\end{equation}
In fact, we have 
\[ C_U \cap \fq 
= \{ x \in \fq \: (x,-x) \in C_{\pi_\rho} \} 
\supeq  C_\g \cap \fq, \] 
so that $C_\g$ and $C_U$ are invariant cones in $\g$ with 
intersecting interiors, and this implies that 
$C_\g^c = C_U^c$ by Corollary~\ref{cor:red-minmax}. 
With  \eqref{eq:svrho} and \eqref{eq:conesequal} we now obtain 
\[ q_M(S_{\sV_{\rho,e}}^0) 
= G^h_e.\Exp(C_\g^{c,0}\cap \fq) 
= W_M(h)_{eH},\] 
and thus by \eqref{eq:wedge-stand} in Theorem~\ref{thm:6.3}
\begin{align*}
 \sV_{\rho} 
&= \sH^{G,U}_\sE(S_{\sV_\rho,e}^0)
= \sH^{G,U}_\sE(S_{\sV_\rho,e}^0 H)
= \sH^{G,U}_\sE(q_M^{-1}(q_M(S_{\sV_\rho,e}^0))) \\
&= \sH^{G,U}_\sE(q_M^{-1}(W_M(h)_{eH}))
= \sH^M_\sE(W_M(h)_{eH}).
\qedhere\end{align*}
\end{prf}

\begin{rem} (The kernel of  $U$) 
(a) For a representation $(\rho,\cK)$ of $G$, 
the representation $\pi_\rho$ on $\cH_\rho \subeq B_2(\cK)$ 
has the property that $\Delta_{Z(G)}$ acts trivially. 
In fact, the subspace $\cH_\rho$ is generated by 
$\rho(C^\infty_c(G,\R)) \subeq B_2(\cK)$ which commutes with 
$\rho(Z(G))$. We therefore have 
for $A \in \cH_\rho$ and $z \in Z(G)$ 
\[ \pi_\rho(z,z) A 
= \rho(z) A \rho(z)^{-1} = A. \] 
For the representation $U(g) := \pi_\rho(g,\tau^G(g))$ of $G$, 
we obtain for $z \in Z(G)$ that 
\[ U(z) A =  \rho(z) A \rho(z^\sharp) = \rho(zz^\sharp)A, \] 
so that $zz^\sharp \in \ker(\rho)$ is equivalent to $z \in \ker(U)$. 
We conclude that 
\begin{equation}
  \label{eq:kerU}
 \ker U = \{ z \in Z(G) \: zz^\sharp \in \ker\rho\}.
\end{equation}

\nin (b) If $\ker(\rho) = Z(G)$ (Proposition~\ref{prop:8.1}), 
the representation  $(U,\cH_\rho)$ factors through $G/Z(G)$.
We then have 
\[ G_\sV 
= \{ g \in G^h \: g g^\sharp\in\ker(U) \} 
= \{ g \in G^h \: gg^\sharp  \in Z(G) \} 
 = \Ad^{-1}(\Ad(G)^{\tau}).\] 

\nin (c) If $\rho$ is faithful (Proposition~\ref{prop:8.1}), 
then 
\[  \ker U = \{ z \in Z(G) \: zz^\sharp =e \} 
= Z(G)^{\tau^G}\]
and thus 
\begin{align} 
  \label{eq:gv-faithf}
 G_\sV 
&= \{ g \in G^h \: g g^\sharp\in \ker(U) \} 
= \{ g \in G^h \: gg^\sharp  \in Z(G)^{\tau^G} \} \notag\\
&= \{ g \in G^h \: gg^\sharp  \in Z(G), (gg^\sharp)^2 = e \}.
\end{align}
\end{rem}

\section{The wedge space as an ordered symmetric 
space} 
\mlabel{sec:9}

In this section we return to a geometric topic. 
We show that, under rather natural assumptions, 
the wedge space $\cW(M,h)$, consisting of all $G$-translates 
of the connected  component $W_M(h)_{eH}$ of the wedge domain in $M$, 
carries the structure of an ordered symmetric space. 
We start in Subsection~\ref{subsec:9.1} with the order 
structure and the determination of the stabilizer group of 
$W(M,h)_{eH}$. This turns out to be rather hard to do directly, 
so we use the representation theoretic results 
from Subsection~\ref{subsec:6.3} on the correspondence between wedge domains 
and standard subspaces to obtain all required information. 
The symmetric space structure is discussed in Subsection~\ref{subsec:9.2}. 
This result establishes a bridge to the abstract wedge spaces 
introduced in \cite{MN21}.

\subsection{The order structure} 
\mlabel{subsec:9.1}

Let $(G,\tau^G)$ be a connected symmetric Lie group corresponding 
to the modular compactly causal symmetric Lie algebra 
$(\g,\tau,C,h)$, let $H \subeq G^{\tau^G}$ be an open subgroup, 
$M = G/H$ the corresponding compactly causal symmetric space 
and 
\[ W := W_M(h)_{eH} \subeq M \] 
 the connected component of the 
wedge domain specified by $h$ corresponding to the base point~$eH$. 
As in Theorem~\ref{thm:6.5}, we assume that 
$C = C_\g \cap \fq$ for a pointed invariant cone $C_\g \subeq \g$ with 
$-\tau(C_\g) = C_\g$.  
{\it We further assume the context of} Theorem~\ref{thm:8.12}, where 
$(U,\cH)$ is an antiunitary representation of $G_{\tau_h}$ with discrete kernel 
and there exists a distribution vector $\eta^0$ with 
\[ \sH^M_{\R \eta^0}(W) = \sV := \Fix(U(\tau^G_h) e^{\pi i \cdot \partial U(h)}).\] 

The {\it wedge space} 
\[ \cW(M,h) := \{ g W \: g \in G \}  \] 
is the $G$-orbit of $W = W_M(h)_{eH}$ in the set of subsets of $M$. 
The order on the wedge space is determined by the subsemigroup 
\begin{equation}
  \label{eq:wvid}
 S_W := \{ g \in G \: gW \subeq W \}
\quad \mbox{ with } \quad 
S_W \cap S_W^{-1} = G_W.
\end{equation}

Below we abbreviate $G_e^h := (G^h)_e$.

\begin{thm} \mlabel{thm:9.1}
The compression semigroup of $W$ is 
\[ S_W = G_W \exp(C_\g^c) \quad \mbox{  with  } \quad G_W = G^h_e H^h.\]
Furthermore, $G_W$ is open in $G^h$.
\end{thm}

\begin{prf} By definition, 
\[ W = W_M(h)_{eH} = G^h_e.\Exp(C_+^0 - C_-^0) 
= S(C_\g^0,h)_e.eH \] 
is the orbit of the base point under the interior of the semigroup 
$S(C_\g,h)_e = G^h_e \exp(C_\g^c)$ 
(cf.\ Proposition~\ref{prop:5.2}(b)). 
This shows that 
\begin{equation}
  \label{eq:semincl}
S(C_\g,h)_e\subeq S_W,
\end{equation}
and in particular $G^h_e \subeq G_W$.
For $s \in S_W$, the inclusion $s W \subeq W$ implies 
\[ U(s)\sV 
= U(s)\sH_\sE^M(W)
= \sH_\sE^M(sW) \subeq  \sH_\sE^M(W) = \sV,\] 
so that $S_W  \subeq S_\sV.$ 
We thus obtain with \eqref{eq:semincl}
\begin{equation}
  \label{eq:95}
G^h_e \exp(C_\g^c) \subeq  S_W \subeq S_\sV 
\ {\buildrel\eqref{eq:sc1} \over \subeq}\ G^h \exp(C_\g^c).
\end{equation}
It follows in particular that 
$G^h_e \subeq G_W  \subeq G^h$. This implies in particular that $G_W$ is
open in~$G^h$. As 
$G^h_e.eH = M^\alpha_{eH}$ by Proposition~\ref{prop:5.3} 
and this submanifold equals $M^\alpha \cap \oline{W}$,
it follows that  
\[ G_W = \{ g \in G^h \: g M^\alpha_{eH} = M^\alpha_{eH} \} 
= G^h_e (G^h \cap H) = G^h_e H^h.\] 

By \eqref{eq:95}, $\exp(C_\g^c) \subeq S_W$, so that 
$G_W \exp(C_\g^c) \subeq S_W$. For the converse, let
$s = g \exp x  \in S_W \subeq G^h \exp(C_\g^c)$. 
Then 
\[ s.eH = g \exp(x).eH \in g W_M(h)_{eH} \subeq W_M(h) \]  
implies $g \in  G_W$. 
\end{prf}

The preceding proof makes heavy use of representation 
theoretic results to determine $S_W$. It would be nice to have a more 
direct geometric proof, but we presently do not see how to do that.

\subsection{The symmetric space structure} 
\mlabel{subsec:9.2}

We now turn to the structure of a symmetric space on 
the wedge space $\cW(M,h)$ . 

\begin{prop} \mlabel{prop:9.2} Let $(\g,\tau)$ be a semisimple 
symmetric Lie algebra, $G$ a connected Lie group with 
Lie algebra $\g$, $\tau^G$ an involutive automorphisms of $G$ 
integrating $\tau$ and $H \subeq G$ a closed subgroup 
satisfying 
\[ \L(H) = \fh \quad \mbox{ and } \quad 
\tau^G(H) = H.\] 
Then $M = G/H$ carries the structure of a symmetric space. 
\end{prop}

\begin{prf} In view of Remark~\ref{rem:9.2}(b), 
we may factor the subgroup $H\cap Z(G)$ which acts trivially on $M$ and 
thus assume that 
$\Ad$ is injective on $H$. Then Remark~\ref{rem:9.2}(d) 
implies that $H \subeq G^{\tau^G}$, and since 
$\L(H) = \fh =\g^\tau =  \L(G^{\tau^G})$, the subgroup $H$ is open in $G^{\tau^G}$. 
Therefore $M = G/H$ is a symmetric space. 
\end{prf}

\begin{prop} \mlabel{prop:9.3}  
{\rm(The wedge space as a symmetric space)} 
 Let $(\g,\tau,C,h)$ be a non-compactly causal semisimple 
symmetric Lie algebra, $G$ a connected Lie group with Lie algebra 
$\g$ and assume that 
\begin{itemize}
\item $\tau$ integrates to an involutive automorphism 
$\tau_h^G$ of $G$, 
\item $H \subeq G^{\tau^G}$ is an open subgroup, 
\item $\tau_h = e^{\pi i \ad h}$ integrates to an involutive automorphism 
$\tau_h^G$ of $G$, 
\item $\tau_h^G(H) = H$, i.e., 
$\tau_h$ induces an involution $\tau_h^M$ on 
$M := G/H$. 
\end{itemize} 
Then the wedge space 
\[ \cW(G,h) = G.W \cong G/G_W \] 
is a symmetric space. 
\end{prop}

\begin{prf} Recall that $G_W = G^h_e H^h \subeq G^h$ 
(Theorem~\ref{thm:9.1}). 
First we observe that 
$\tau^G_h(G^h) = G^h$ follows from $\tau(h) = h$. 
With $\tau_h^G(H) = H$, we therefore  get 
\[ \tau^G_h(G_W) = \tau^G_h(G^h_e H^h) = G^h_e H^h = G_W.\] 
Hence the assertion follows  from Proposition~\ref{prop:9.2}. 
\end{prf}

\begin{rem} In the context of Subsection~\ref{subsec:9.1}, we 
have four subgroups with the Lie algebra $\g_0(h) =  \ker(\ad h)$: 
\[ G^h_e \subeq G_W \subeq G_\sV \subeq G^h.\] 
Accordingly, we have a sequence of coverings of homogeneous spaces 
\[ G/G^h_e \to \cW(M,h) \cong G/G_W \to U(G)\sV \cong G/G_\sV 
\to \cO_h = \Ad(G)h \cong G/G^h.\] 
Here $\cO_h \cong \Ad(G).h$ always is a symmetric space because 
$\tau_h$ defines an involution on $\Ad(G)$ and 
$\Ad(G)^h \subeq \Ad(G)^{\tau_h}$ is an open subgroup. 

Likewise, $G/G_\sV \cong \cO_\sV := U(G)\sV$ 
carries a natural symmetric space structure. 
The point reflection in the base point $\sV \in \cO_\sV$ 
is given by $\sW \mapsto J_\sV \sW'$.
We refer to \cite{Ne19} for more 
on the geometry of the space $\Stand(\cH)$ of standard subspaces 
of $\cH$. 
\end{rem}

\begin{rem} (a) The three  conditions in Proposition~\ref{prop:9.3} 
are in particular satisfied for 
spaces of group type $M = G \cong (G\times G)/\Delta_G$, where 
$\tau^{G \times G}(g_1, g_2) = (g_2, g_1)$ 
and $\tau^G_h$ exists on $G$, so that $\tau^{G \times G}_h = \tau^G_h \times \tau^G_h$.   

\nin  (b) For Cayley type spaces we have 
$\tau = \tau_h$, so that 
$G_W = G^h_e H^h \subeq G^{\tau^G_h}$ follows 
from $H \subeq G^{\tau^G_h}$. 

\nin (c) If $Z(G) = \{e\}$ and $G \cong \Ad(G)$, then 
\[ M_G \cong \{ \phi \in \Aut(G)_e \: \tau \phi \tau = \phi^{-1} \}, \quad 
\tau_h^G(\phi) = \tau_h \phi \tau_h, \qquad 
\tau^G(\phi) = \tau \phi \tau.\] 
Further, $G_W \subeq G^h \subeq G^{\tau^G_h}$ 
(Remark~\ref{rem:7.6new}), so that $\cW(M_G,h)$ is a symmetric space. 

\nin (d) For $M = M_G$ we have $H = G^{\tau^G}$, which is $\tau_h^G$-invariant. 
If $G$ is simply connected, then $H = G^{\tau^G}$ is connected, so that we have 
in particular $M = G/H \cong M_G$. 

\nin (e) If the representation $U$ in Theorem~\ref{thm:8.12} is faithful, then 
\[ G_W \subeq G_\sV   G^{h, \tau_h^G} \subeq G^{\tau_h^G},\] 
which obviously implies the invariance of $G_W$ under $\tau^G_h$. 
\end{rem}

\section{Open problems} 
\mlabel{sec:10}

In this section we briefly discuss some interesting open problems that we plan to 
address in the future. 

\begin{prob} For a compactly causal symmetric space $M = G/H$ with the 
infinitesimal data $(\g,\tau,C,h)$, 
characterize the antiunitary representations 
$(U,\cH)$ of $G$ for which there exists a real subspace $\sE \subeq \cH^{-\infty}$ 
satisfying 
\[ \sH_\sE^{G/H}(W_M(h)_{eH}) = \sV_{(h,\tau^G,U)}.\] 
We conjecture  that $C_\g$-positivity should be enough, but this is not easy 
to see.
\end{prob}

\begin{prob} The symmetric spaces $\cW(M,h)$ of wedge domains 
in $M$ are particular examples of non-compactly causal 
parahermitian symmetric spaces. The corresponding causal 
symmetric Lie algebra is 
$(\g,\tau_h, C_\g^c)$ (Theorem~\ref{thm:9.1}). 
Up to coverings, we thus obtain interesting realizations of 
non-compactly causal symmetric spaces of Cayley type. 

However, there are much more parahermitian symmetric spaces 
(cf.\ \cite{Kan00}, \cite{MN21}), 
and we expect that they can likewise be realized as wedge spaces 
$\cW(M,h)$ of a non-compactly causal symmetric spaces $M = G/H$ 
(see \cite{NO21b} for details and results in this direction). 
If $(\g,\tau,C)$ is non-compactly causal and $\g$ is simple, 
$h \in \cE(\g)$, then $(\g,\tau_h)$ can be any simple parahermitian 
symmetric Lie algebra by \cite[Lem.~2.6]{O91}. 
In this case $G/G^h_e$ carries no natural order structure, which is also 
reflected in the geometry of $\cW(M,h)$,  respectively 
the fact that the inclusion order on this set is trivial.

For example, if $M = G_\C/G$ is irreducible 
non-compactly causal of complex type, 
corresponding to $(\g_\C,\tau, i C_\g, h)$, where $\g$ is simple hermitian 
(not necessarily of tube type),  and $h$  a causal Euler element, 
then $G_\C^h \cong K_\C$ is a connected subgroup, and 
\[ \cW(G_\C/G, h) \cong G_\C/K_\C\] 
is the complexification of the Riemannian symmetric space $G/K$. 
\end{prob}

\section{Anti-de Sitter space $\AdS^d$} 
\mlabel{subsec:a.4}

In this section we discuss the $d$-dimensional anti-de Sitter space 
$\AdS^d$ as an example of a compactly causal symmetric space. 
As we verify all assertions for this example directly, it can be read 
independently of the other sections. It illustrates our results for this 
concrete example and it is a guiding example for the general 
theory. From the perspective of physics it is of particular importance 
because anti-de Sitter spaces are precisely the {\it irreducible} 
compactly causal symmetric spaces which are Lorentzian. There are more 
reducible examples, even compact ones, such as the conformal compactification 
of Minkowski space, considered as a symmetric space of the compact group 
$\SO_2(\R) \times \SO_d(\R)$. 

We endow $V := \R^{d+1}$ with  the  symmetric bilinear form 
\[ \beta(x,y) := x_0 y_0 + x_1 y_1 - x_2 y_2 - \cdots - x_d y_d.\] 
We consider the connected group $G := \SO_{2,d-1}(\R)_e$ and consider 
the {\it anti-de Sitter space} 
\[ M := \AdS^d := \{ x \in V \: \beta(x,x) = 1\} = G.e_0 \cong G/G^{e_0} \] 
and note that 
\begin{equation}
  \label{eq:hconnected}
 H := G^{e_0}  \cong  \SO_{1,d-1}^\up(\R)
\end{equation}
is connected. 
By  \cite[Prop.~B.3]{NO21b}, the 
exponential function of the symmetric space $\AdS^d$ is given by 
\begin{equation}
  \label{eq:expfunct-ads}
 \Exp_p(v) = C(\beta(v,v)) p + S(\beta(v,v)) v,
\end{equation}
where the
functions $C, S \: \C \to \C$ are defined by 
\begin{equation}
  \label{eq:CandS}
 C(z) := \sum_{k = 0}^\infty \frac{(-1)^k}{(2k)!} z^{k} \quad \mbox{ and } \quad 
 S(z) := \sum_{k = 0}^\infty \frac{(-1)^k}{(2k+1)!} z^{k}. 
\end{equation} 
Then  $C(z^2) = \cos (z)$ and $S(z^2)= \frac{\sin z}{z}$.
It follows in particular that $\beta$-positive vectors $v \in T_p(M)$ 
(timelike vectors) 
generate closed geodesics. 

The domain 
\[ V_{> 0} := \{ v \in V \: \beta(v,v) > 0\} \] 
intersects $T_{e_0}(\AdS^d) \cong e_0^\bot$ in a double light cone, 
and we write 
\[ V_+(e_0) := \{ x \in T_{e_0}(\AdS^d) = e_0^\bot 
\: \beta(x,x) > 0, \beta(x, e_1) > 0\}\] 
for the ``positive'' cone containing~$e_1$. As this cone is 
$G^{e_0}$-invariant by \eqref{eq:hconnected}, we obtain a {\it causal structure} on $M$ by 
\[ V_+(g.e_0) := g.V_+(e_0) \quad \mbox{ for } \quad g \in G.\]

The tangent bundle of $M$ is 
\[ T(M) = \{ (x,y) \in V \times V \: \beta(x,x) = 1, \beta(x,y) = 0\}.\] 
As $G$ preserves the causal orientation on $M$, we call 
a pair $(x,y) \in T(M)$ 
\begin{itemize}
\item {\it positive} if $G.(x,y) \cap  V_+(e_0) \not=\eset$, 
and 
\item {\it negative} if $G.(x,y) \cap  -V_+(e_0) \not=\eset$. 
\end{itemize}
Then 
\[ T(M)_\pm = \bigcup_{m \in M} \pm V_+(m)\] 
is the open subset of positive/negative vectors in the tangent bundle.
The reflection $r_{01}(x) = (-x_0, -x_1, x_2, \ldots, x_d)$ is 
contained in $G$, so that 
\begin{equation}
  \label{eq:v+base}
V_+(-e_0) = r_{01}(V_+(e_0)) = - V_+(e_0).
\end{equation}

We consider the boost vector field $X_h(v) = hv$, where 
$h \in \so_{2,d-1}(\R)$ is defined by 
\[ hx = (0, x_2, x_1, 0,\ldots, 0).\] 
It generates the flow 
\[ \alpha_t(x) = e^{th} x 
= (x_0, \cosh t \cdot x_1 + \sinh t \cdot x_2, 
\cosh t \cdot x_2 + \sinh t \cdot x_1, x_3, \ldots, x_d) .\] 
It also leads to the involutions 
\[ \tau_h(x) = e^{\pi i h}x = (x_0, -x_1, -x_2, x_3, \ldots, x_d) \mbox{ on } 
\R^{d+1}  \quad \mbox{ and }\quad \tau^G_h (g) = \tau_h g \tau_h 
\mbox{ on } G,\] 
and the Wick rotation 
\[ \kappa_h(x) = e^{-\frac{\pi i}{2}  h}x
 = (x_0, -i x_2, -i x_1, x_3, \ldots, x_d)\]
with $\tau_h = \kappa_h^2$. 
Note that $\tau_h \not\in G$ by \eqref{eq:hconnected}. This element 
reverses the causal structure on $\AdS^d$.
  
\subsection{Subgroup data}
\mlabel{subsec:subgrpdata}

We now describe the groups $G^{\tau_h^g}$ and $G^h$ in more detail. 
The subgroup $G^h$ consists of all elements $g$ in $G$ such that $gh = hg$.
Hence $G^h$   can also be described as the group of all elements $g\in G$ leaving all 
$h$-eigenspaces invariant: 
\[ G^h = \{ g \in G \: g(e_1 \pm e_2) \in \R (e_1 \pm e_2), 
g\{e_1,e_2\}^\bot = \{e_1,e_2\}^\bot\}.\]
Writing $V_h := hV = \R e_1 \oplus \R e_2 \cong \R^{1,1}$, we have 
\begin{align*}
 G^{\tau_h} 
&= \{ g \in G \: g V_h = V_h \}  
=G \cap (\OO_{1,1}(\R) \times \OO_{1,d-2}(\R))\\
&\supeq G^h \supeq G^h_e = \SO_{1,1}^\up(\R) \times \SO_{1,d-2}^\up(\R).
\end{align*}
A maximal compact subgroup $K$ of $G$ is 
\[ K = G\cap \OO_{d+1}(\R) \cong \SO_2(\R) \times \SO_{d-1}(\R).\] 
As  
\[ K^{\tau_h} 
= \{ \pm \1\} \times {\rm S}(\OO_1(\R) \times \OO_{d-2}(\R))
\cong \{ \pm \1\} \times \OO_{d-2}(\R) \] 
has two connected component, polar decomposition implies 
\begin{equation}
  \label{eq:ghcomp}
  \pi_0(G^{\tau_h}) \cong 
  \begin{cases}
    \Z_2 \times \Z_2 & \text{ for } d > 2 \\ 
    \Z_2  & \text{ for } d = 2. 
  \end{cases}
\end{equation}
The group $G^{\tau_h}$ acts on the Minkowski plane $V_h\cong \R^{1,1}$ by a 
homomorphism 
\[ \gamma \: G^{\tau_h} \to \OO_{1,1}(\R) \] 
which is surjective for $d > 2$, but not for $d = 2$. We have  
\[  G^{\tau_h} 
= G^h\cup r_{01} G^h,\quad \mbox{ where }\quad r_{01}  = \diag(-\1_2,\1_{d-1}). \]  
We also note that 
\[ \Ad(G^{\tau_h})h = \{ \pm h \} \quad \mbox{ and } \quad 
G^{\tau_h}\setminus G^h = \{ g \in G \: \Ad(g) h = -h\}.\] 
It follows in particular that $G^h$ is connected if 
$d=2$ and has two connected components if $d>2$:
\begin{equation}\label{eq:Gh}
G^h = G^h_e\cup rG^h_e, \quad\text{where } r = \diag(-\1_4,\1_{d-3})
\end{equation}

The subgroup $H = G^{e_0} \cong \SO_{1,d-1}^\up(\R)$ satisfies 
\[ H^{\tau_h} = {\rm S}(\OO^\up_{1,1}(\R) \times  \OO_{d-2}(\R)),\] 
which has two connected components for  $d > 2$ and one for $d = 2$.

The submanifold of $\alpha$-fixed points in $M$ is 
\[ M^\alpha = \{ (x_0, 0, 0, x_3, \ldots, x_d) \: 
x_0^2 - x_3^2 - \cdots - x_d^2 = 1 \} 
\cong \bH^{d-2} \cup - \bH^{d-2}. \] 
Proposition \ref{prop:5.3} implies that  $G^h_e. (\pm \bH^{d-2})= \pm \bH^{d-2}$. 
The element  $r\in G^h$ from \eqref{eq:Gh} satisfies $r(e_0)=-e_0$, 
so that $G^h$ acts transitively on $M^\alpha$.

\subsection{The tube domain} 

The tube domain 
\[ \cT_V := V + i V_{> 0} \subeq V_\C \] 
is connected, but its intersection with 
\[ M_\C := \{ z = x + i y \: 
\beta(x,x) - \beta(y,y) = 1, \beta(x,y) = 0\} \] 
decomposes into two connected components. In fact, 
$z \in M_\C \cap \cT_V$ implies 
\[ \beta(x,x) = 1 + \beta(y,y) > 1,\]
 so that  
\[ \Gamma(z)  := \frac{1}{\sqrt{\beta(x,x)}}(x,y) \in T(M),\] 
and the map $\Gamma : M_\C \cap \cT_V \to T(M)$ 
is $G$-equivariant. The tangent vectors 
$\Gamma(z) = (\tilde x, \tilde y)$ are $\beta$-positive 
with $\beta(\tilde y,\tilde y) < 1$. If, conversely, 
$(\tilde x, \tilde y) \in T(M)$ satisfies 
$0 < \beta(\tilde y, \tilde y) < 1$, then 
\[ c := \beta(\tilde x, \tilde x) - \beta(\tilde y, \tilde y) 
= 1- \beta(\tilde y, \tilde y) > 0,\] 
so that $z := c^{-1/2}(\tilde x + i \tilde y) \in M_\C \cap \cT_V$. 
This shows that $\Gamma$ is a diffeomorphism 
\[ \Gamma \: M_\C \cap \cT_V \to \{ (x,y) \in T(M) \: 0 < \beta(y,y) < 1 \}.\] 
In particular, $M_\C \cap \cT_V$ has two $G$-invariant 
connected components 
\begin{equation}
  \label{eq:tmpm}
\cT_M^\pm := \Gamma^{-1}(T(M)_\pm),
\end{equation}
which are called the {\it tube domains of $M$}. 

The following lemma describes how the positive tube domain 
$\cT_M^+$ can be parametrized by the exponential function of $M_\C$.

\begin{lem} \mlabel{lem:a.42}
We have 
\begin{equation}
  \label{eq:desit-tube2}
\cT_M^+= G.\Exp_{e_0}(i V_+(e_0)) = \Exp(i T(M)_+), 
\end{equation}
and the fixed point set of the antilinear extension $\oline\tau_h$ of $\tau_h$ 
in the complex manifold $\cT_M^+$ is 
\[ (\cT_M^+)^{\oline\tau_h} 
= G^{\tau_h}.\Exp_{e_0}(i V_+(e_0)^{-\tau_h}) 
= G^{\tau_h}.\Exp_{e_0}(i \R_+ e_1).\] 
\end{lem}

\begin{prf} The second equality in \eqref{eq:desit-tube2} 
follows from the $G$-equivariance of 
the exponential function of $M_\C$. 
For the first equality in \eqref{eq:desit-tube2}, 
we start with the observation that 
both sides are $G$-invariant. 

\nin ``$\supeq$'':  We have to show that any $x \in V_+(e_0)$ satisfies 
$z := \Exp_{e_0}(ix) \in \cT_M^+$. 
Since the orbit of $x$ under 
$G^{e_0}$ contains an element in $\R_+ e_1$, we may assume that 
$x = x_1 e_1$ with $x_1 > 0$. Then 
\[ z= \Exp_{e_0}(ix) = \Exp_{e_0}(ix_1 e_1) 
= C(-x_1^2) e_0 + S(-x_1^2)  x_1 i e_1 
= \cosh(x_1) e_0 + \sinh(x_1) i e_1 \] 
and $x_1 > 0$ imply 
\[ \Gamma(z) = (e_0, \tanh(x_1) e_1) \in T(M)_+,\] hence that 
$z \in \cT_M^+$. 

\nin ``$\subeq$'': Conversely, let $z = x + i y\in \cT_M^+$. 
Acting with $G$, we may thus assume that 
$y = y_1 e_1$ with $y_1 > 0$. Then $x \in y^\bot = e_1^\bot$ follows from 
$\Im \beta(z,z)  = 0$, so that (acting with $G^{e_1}$) we 
may further assume that $x = x_0 e_0$ for some $x_0 \not=0$. 
Now $z = x_0 e_0 + i y_1 e_1\in \AdS^d_\C$ implies that 
\[ 1 = \beta(z,z) = x_0^2 - y_1^2 .\] 
Since 
\[ \Gamma(z) = \frac{1}{|x_0|}(x_0 e_0, y_1 e_1) \in T(M)_+,\] 
and $y_1 > 0$, we also have $x_0 > 0$. 
Hence there exists a $t > 0$ with 
$y_1 = \sinh t$ and $x_0 = \cosh t$. We then have 
\[ \Exp_{e_0}(i t e_1) = 
\cosh(t) e_0 + \sinh(t) i e_1 
= x_0 e_1 + y_0 i e_1 = z.\] 
This proves \eqref{eq:desit-tube2}. 

For the second assertion, we start from 
$\cT_M^+ = G.\Exp_{e_0}(i V_+(e_0)),$ which immediately shows that 
\[ (\cT_M^+)^{\oline \tau_h} 
\supeq G^{\tau_h}.\Exp_{e_0}(i V_+(e_0)^{-\tau_h}), \] 
where 
\[ V_+(e_0)^{-\tau_h} 
= V_+(e_0) \cap (\R e_1 + \R e_2) 
= \{ x_1 e_1 + x_2 e_2 \:  x_1 > 0,  x_1 > | x_2| \}.\] 

We now show that we actually have equality. 
So let $z = x + i y \in (\cT_M^+)^{\oline\tau_h}$. 
Then 
$y = y_1 e_1 + y_2 e_2 \in V^{-\tau_h}$ is $\beta$-positive, i.e., 
$y_1^2 > y_2^2$. Acting with 
$G^h$, we may thus assume that 
$x = x_0 e_0$ and $y = y_1 e_1$. 
From $1 = \beta(z,z) = x_0^2 - y_1^2$ and 
$z = x + i y = x_0 e_0 + i y_1 e_1 \in \cT_M^+$ it now follows that 
$x_0 y_1 > 0$. 
Acting with the rotation $r_{01}(x) = (-x_0,-x_1, x_2, \ldots, x_d)$, 
which is contained in $G^{\tau_h}$, we may further assume that $x_0 > 0$, so that 
$y_1 >0$ holds as well. 
Hence there exists a $t > 0$ with 
$y_0 = \sinh(t)$ and $x_0= \cosh(t)$. As above, 
we now obtain 
\[ \Exp_{e_0} (ite_1) = \cosh(t) e_0 + \sinh(t) i e_1 = z.\] 
This implies the second assertion. 
\end{prf}

\subsection{The wedge domains} 

For the closed convex $H$-invariant cone 
$C := \oline{V_+(e_0)} \subeq T_{e_0}(M) \cong \fq$, we obtain 
\[ C_+ = [0,\infty) (e_2 + e_1) \quad \mbox{ and } \quad 
C_- = [0,\infty) (e_2 - e_1). \] 

\begin{lem} \mlabel{lem:b.2new} We have 
\[ W_M(h)_{e_0} = G^h_e.\Exp_{e_0}(C_+^0 + C_-^0) 
\quad \mbox{ and }  \quad 
W_M(h)_{-e_0} = G^h_e.\Exp_{-e_0}(-C_+^0 - C_-^0).\] 
\end{lem}

\begin{prf} 
First we recall that 
\[  W_M(h)_{e_0}  
= G^h_e.\Exp_{e_0}(C_+^0 + C_-^0)
=G^h_e.\Exp_{e_0}(\R_+(e_1 + e_2) - \R_+(e_1 - e_2)).\] 
For $g = r_{01} = \diag(-1,-1,1,\ldots, 1)\in G$ we have 
$\Ad(g)h = -h$, $g e_0 = -e_0$ and 
$gC_+ = C_-$. With Remark~\ref{rem:6.5} we further obtain 
\begin{align*}
 W_M(h)_{-e_0}
&=  W_M(h)_{g.e_0}
= g.W_M(-h)_{e_0}
= gG^h_e.\Exp_{e_0}(-(C_+^0 + C_-^0))\\
&= G^h_e.g\Exp_{e_0}(-(C_+^0 + C_-^0))
= G^h_e.\Exp_{-e_0}(-g.(C_+^0 + C_-^0))\\
&= G^h_e.\Exp_{-e_0}(-(C_+^0 + C_-^0)).
\qedhere\end{align*}
\end{prf}

\begin{lem}
  \mlabel{lem:a.4.1} 
The positivity domain 
\[ W_M^+(h) = \{ x \in M \: X_h(x) \in V_+(x) \} \] 
of the vector field $X_h$ satisfies 
\begin{equation}
  \label{eq:posdomdesc}
 W_M^+(h) = G^h.\{ x_0 e_0 + x_2 e_2 \in M \: x_0 x_2 > 0\}.  
\end{equation}
\end{lem}

\begin{prf} If $X_h(x) \in V_+(x)$, then 
$x_2^2 > x_1^2$ implies that 
\[ x_0^2 
= 1 + \underbrace{x_2^2 - x_1^2}_{>0} + x_3^2 + \cdots + x_d^2 
> 1 + x_3^2 + \cdots + x_d^2 > 0.\] 
Acting with the centralizer $G^h$, which contains the group 
$\SO_{1,d-2}^\up(\R)$, acting on 
\[ \ker(h) = \Spann \{ e_0, e_3, \ldots, e_d\},\] 
we may assume that $x_3 = \cdots = x_d = 0$, so that 
$x = (x_0, x_1, x_2, 0,\ldots, 0)$ is an element of 
the $2$-dimensional anti-de Sitter space~$\AdS^2$. 
 \begin{footnote}
{As a 
manifold $\AdS^2$ identifies naturally with $\dS^2$, but the causal 
structure is different.}   
\end{footnote}
Next we observe that the modular flow $\alpha_t$ moves 
any element $x$ with $x_2^2 > x_1^2$ to an element with $x_1 = 0$. 
For $x = (x_0, 0, x_2,0,\ldots,0) \in \AdS^d$ we have 
$x_0^2 > x_2^2$. 
Now $X_h(x) = x_2 e_1 \in V_+(x)$ is equivalent to 
$x_0 x_2 > 0$ (see \eqref{eq:v+base}). 
As $W_M^+(h)$ is invariant under $G^h \supeq \exp(\R h)$, the assertion 
follows.   
\end{prf}

\begin{lem}
  \mlabel{lem:a.4.2} 
The subset 
\[ W_M^{\rm KMS}(h) 
:= \{ x \in M \: (\forall t \in (0,\pi))\ 
\alpha_{it}(x) \in \cT_M^+ \}\] 
coincides with $W_M^+(h)$. 
\end{lem}

\begin{prf}
For $x \in W_M^{\rm KMS}(h)$ and $0 < t < \pi$ we have 
\[ \alpha_{it}(x) 
= (x_0, \cos t \cdot x_1 + \sin t \cdot i x_2, 
\cos t \cdot x_2 + \sin t \cdot i  x_1, x_3, \ldots, x_d),\] 
so that 
\[ \Im(\alpha_{it}(x)) = 
\sin t \cdot(0,  x_2,  x_1, 0, \ldots, 0).\] 
This element is $\beta$-positive for every $t \in (0,\pi)$ 
if and only if $x_2^2 > x_1^2$. 
Using the $G^h$-invariance of both sides, 
we see as in the proof of Lemma~\ref{lem:a.4.1} 
that we may assume that $x = (x_0, 0, x_2, 0,\ldots, 0)$, 
so that $x_0^2 - x_2^2 = 1$ and $x_2 \not=0$. Then 
\[z_t := x_t + i y_t := \alpha_{it}(x) 
= (x_0,  \sin t \cdot i x_2, \cos t \cdot x_2,0, \ldots, 0) 
=x_0 e_0 + \cos t\cdot x_2 e_2 
+ \sin t \cdot  x_2 i e_1\] 
is $G$-conjugate to 
\[\sgn(x_0) \sqrt{x_0^2 - \cos^2(t) x_2^2} \cdot 
e_0 + \sin t \cdot i x_2 e_1.\] 
Hence $\Gamma(z_t) \in T(M)_+$ for $0 < t < \pi$ is equivalent to 
$x_0 x_2 > 0$. In view of \eqref{eq:posdomdesc}, this is equivalent to  
$x \in W_M^+(h)$.  We conclude that 
$W_M^+(h) = W_M^{\rm KMS}(h).$
\end{prf}

\begin{lem}
  \mlabel{lem:a.4.3} 
$W_M^{\rm KMS}(h) = \kappa_h\big((\cT_M^+)^{\oline\tau_h}\big).$
\end{lem}

\begin{prf} For $m \in W_M^{\rm KMS}(h)$ we have 
\[ 
\oline\tau_h(\alpha_{\frac{\pi i}{2}}(m))
= \alpha_{-\frac{\pi i}{2}}(\tau_h(m))
= \alpha_{-\frac{\pi i}{2}}\alpha_{\pi i}(m)
= \alpha_{\frac{\pi i}{2}}(m),\] 
so that 
\[ \alpha_{\frac{\pi i}{2}}(m) \in (\cT_M^+)^{\oline\tau_h},
\quad \mbox{ and thus } \quad 
 W_M^{\rm KMS}(h) \subeq \kappa_h\big((\cT_M^+)^{\oline\tau_h}\big).\] 

To show that we actually have equality, 
let $z = x + i y \in (\cT_M^+)^{\oline\tau_h}$. 
Then 
$y = y_1 e_1 + y_2 e_2 \in V^{-\tau_h}$ is $\beta$-positive, i.e., 
$y_1^2 > y_2^2$. Further, $\beta(x,x) > \beta(y,y) > 0$. 
Acting with 
$G^h$, we may thus assume that 
$x = x_0 e_0$ and $y = y_1 e_1$. 
From $1 = \beta(z,z) = x_0^2 - y_1^2$ and 
$z = x + i y = x_0 e_0 + i y_1 e_1 \in \cT_M^+$ it now follows that 
$x_0 y_1 > 0$. So 
\[ \kappa_h(z) = (x_0, 0, y_1,0,\ldots, 0) \in W_M^+(h)\] 
follows from $x_0 y_1 > 0$ and \eqref{eq:posdomdesc}. 
This proves the assertion.
\end{prf}

The following theorem represents four different characterizations of the 
wedge domain in~$M$.

\begin{thm} For anti-de Sitter space $M = \AdS^d$, we have the equalities 
\[ W_M^+(h) = W_M^{\rm KMS}(h) 
= \kappa_h\big((\cT_M^+)^{\oline\tau_h}\big) 
= W_M(h). \]
\end{thm}

\begin{prf}
The proof follows from the preceding five lemmas. 
It only remains to show that 
\[ \kappa_h((\cT_M^+)^{\oline\tau_h}_{e_0}) = W_M(h)_{e_0}
\quad \mbox{ and } \quad 
\kappa_h((\cT_M^+)^{\oline\tau_h}_{-e_0}) = W_M(h)_{-e_0}.\] 
The component 
\[ (\cT_M^+)^{\oline\tau_h}_{e_0} 
= G^h_e .\Exp_{e_0}(i V_+(e_0)^{-\tau_h}) 
= G^h_e .\Exp_{e_0}(i(C_+^0 - C_-^0)) \] 
is mapped by $\kappa_h$ to 
\[ G^h_e .\Exp_{e_0}(C_+^0 + C_-^0) = W_M(h)_{e_0}.\]

For $g = (-1,-1,1,\ldots, 1)$ we have $g.e_0 = -e_0$, 
$\Ad(g)h = -h$ and $\Ad(g) C^c = C^c$. Now
\[ (\cT_M^+)^{\oline\tau_h}_{-e_0} 
= G^h_e g_0.\Exp_{e_0}(i V_+(e_0)^{-\tau_h}) \] 
is mapped by $\kappa_h$ to 
\begin{align*}
G^h_e \kappa_h(g_0.\Exp_{e_0}(i V_+(e_0)^{-\tau_h}))
&= G^h_e g_0 \kappa_h^{-1}(\Exp_{e_0}(i V_+(e_0)^{-\tau_h}))\\
&= G^h_e g_0 \Exp_{e_0}( - C_+^0 - C_-^0)
= W_M(h)_{-e_0},
\end{align*}
where the last equality follows from Lemma~\ref{lem:b.2new}. 
This completes the proof.   
\end{prf}

\appendix

\section{Some facts on symmetric Lie groups} 
\mlabel{app:1} 
In this section we collect some basic facts on
symmetric subgroups of a Lie group $G$ and on Olshanski semigroups.

\begin{lem}
  \mlabel{lem:a.1} 
Let $(G,\tau)$ be a symmetric Lie group. 
Then $G^{\sharp} := \{ g \in G \: g^\sharp = g\}$ 
is a submanifold of $G$ containing $e$. 
Its identity component is 
\[ G^{\sharp}_e = \{ g g^\sharp \: g \in G_e \} \cong G_e/(G_e)^\tau.\] 
For $g_0 \in G^{\sharp}$ and the involutive 
automorphism $\sigma(g) := g_0 \tau(g) g_0^{-1}$, the connected component 
of $G^{\sharp}$ containing $g_0$ is 
\[ G^{\sharp}_{g_0} 
= \{ g g_0 g^\sharp \: g \in G_e \}
= \{ g \sigma(g)^{-1}  g_0  \: g \in G_e \}.\] 
\end{lem}

\begin{prf}
In the group $G_\tau := G \rtimes \{\id_G, \tau\}$ 
the involutions form a submanifold $\Inv(G_\tau)$ 
and $G_e$ acts transitively on its connected components by conjugation 
(\cite[Ex.~4.6.12]{GN}). 
From 
\[ \Inv(G_\tau) \cap (G \times \{\tau\}) = G^{\sharp} \rtimes \{\tau\},\] 
we thus derive that  $G^{\sharp}$ is a submanifold of $G$ and that 
\[ G_e.\tau = \{ g g^\sharp \: g \in G_e \} \times \{\tau\}\] 
is a connected component of $\Inv(G_\tau)$. 
Therefore $G^{\sharp}_e$ is the $G_e$-orbit of $e$ with 
respect to the action given by $g.x := g x g^\sharp$.  
The same argument yields the description of $G^{\sharp}_{g_0}$. 
\end{prf}

The following lemma could also be drawn from \cite[\S IV]{KNO97}, but 
for the sake of completeness we provide the complete arguments. 

\begin{lem}
  \mlabel{lem:real-olsh} 
Let $(G,\sigma)$ be a connected symmetric Lie group, 
let $L \subeq G^\sigma$ be an open subgroup, and let 
$C \subeq \g^{-\sigma}$ be an open $\Ad(L)$-invariant hyperbolic 
convex cone. Suppose further that $\exp_G$ is injective on 
$\fz(\g) \cap \fq$, so that we can form the real Olshanski semigroup 
\[ S :=  L \exp(C) \subeq G\]  
for which the polar map $L \times C \to S, (g, x) \mapsto g \exp x$ is a 
homeomorphism. 
Let $\tau$ be an involutive automorphism of $G$ commuting with $\sigma$ 
such that $-\tau(C) = C$ and put $g^\sharp := \tau(g)^{-1}$. 
Then the following assertions hold: 
\begin{itemize}
\item[\rm(1)] $S$ is $\sharp$-invariant and the subset 
$S^\sharp := \{s \in S \: s^\sharp = s\}$ is invariant under the action of $S$ on 
itself by $s.t := st s^\sharp$. 
\item[\rm(2)] The projection 
$p \: S \to L, p(g \exp x) := g$ 
satisfies 
\begin{equation}
  \label{eq:projL}
p(s^\sharp) = p(s)^\sharp \quad \mbox{ and } \quad 
p(g.s) = g p(s) g^\sharp \quad \mbox{ for } \quad s \in S, g \in L.
\end{equation}
\item[\rm(3)] The connected component $S^\sharp_e$ 
of $S^\sharp$ containing $\exp(C^{-\tau})$ coincides 
with 
\[ L_e.\exp(C^{-\tau}) = \{ g \exp(x) g^\sharp \: g \in L_e, x \in C^{-\tau}\}\] 
and the map 
\begin{equation}
  \label{eq:lpolar}
 L_e \times_{(L_e)^\tau} C^{-\tau} \to S^\sharp_e,\quad 
[g, x] \mapsto g \exp(x) g^\sharp 
\end{equation}
is a diffeomorphism. 
\item[\rm(4)] $S^\sharp_e = \{ ss^\sharp \: s \in S_e \}$. 
\item[\rm(5)] For each connected component $\cC$ of $S^\sharp$ there exists an 
element $\ell_0 \in L^\sharp \cap \cC$ such that the involution 
$\gamma(g) := \ell_0^{-1} \tau(g) \ell_0$  of $G$ satisfies
\begin{itemize}
\item[\rm(a)] $\gamma$ commutes with $\sigma$. 
\item[\rm(b)] $-\gamma(C) = C$. 
\item[\rm(c)] The connected component $\cC$ is of the form 
\[ S^\sharp_{\ell_0} 
= \ell_0 \{ s \gamma(s)^{-1} \: s \in S_e\} 
= \bigcup_{\ell_1 \in L_e} \ell_0 \ell_1 \exp(C^{-\gamma}) \ell_0^{-1} 
\ell_1^\sharp \ell_0,\]  
which is diffeomorphic to $L_e \times_{(L_e)^\gamma} C^{-\gamma}.$
\end{itemize}
\end{itemize}
\end{lem}

\nin If $\g$ is semisimple, then $\fz(\g) = \{0\}$,
so that the above condition on the exponential function is automatically 
satisfied. It is also satisfied if $G$ is simply connected.

\begin{prf} For the existence of real Olshanski semigroups, we refer to 
\cite[Thm.~3.1]{La94}.  
To see that the assumptions of this theorem are satisfied, 
note that, as  $Z(G)$ is a $\sigma$-invariant closed subgroup, so is 
$Z(G)^{-\sigma}_e = \exp(\fz(\g) \cap \fq)$. Therefore the injectivity of 
$\exp$ on $\fz(\g) \cap  \fq$ implies in particular that, for every 
$x \in \fq \cap \fz(\g)$, the one-parameter subgroup $\exp(\R x)$ is closed 
and isomorphic to $\R$. 

\nin (1) follows immediately from $L^\sharp = L$ and 
$(\exp x)^\sharp = \exp(-\tau(x)) \in \exp(C)$ for $x\in C$. 

\nin (2) From the bijective polar decomposition 
(cf.\ Lawson's Theorem \cite[XI.1.7]{Ne00} and \cite{La94}), 
we obtain the map $p \: S \to L$ and 
\eqref{eq:projL} is easily verified. 

\nin (3) By (2), the subset $S^\sharp_e$ is $L_e$-invariant and contains $\exp(C^{-\tau})$. From \eqref{eq:projL} we get $p(S^\sharp) \subeq L^\sharp$, 
and by continuity $p(S^\sharp_e) \subeq L^\sharp_e$. 
As 
\[ L^\sharp_e= \{ gg^\sharp \: g \in L_e \} \cong L_e/(L_e)^\tau \]
by Lemma~\ref{lem:a.1}, the group $L_e$ acts transitively on $L^\sharp_e$. 
With \eqref{eq:projL} we  conclude that 
\[ S^\sharp_e = L_e.p^{-1}(e), \quad \mbox{ where } \quad 
p^{-1}(e) = \exp(C^{-\tau}).\] 
This proves the first part of (3). For the second part we observe that 
$g_1 \exp(x_1) g_1^\sharp  = g_2 \exp(x_2) g_2^\sharp$ is equivalent to 
$g_2^{-1} g_1 \in L^\tau$ (which implies 
$(g_2^{-1} g_1)^\sharp = (g_2^{-1} g_1)^{-1}$), and $\Ad(g_2^{-1} g_1)x_1 = x_2$. 
Therefore the map in \eqref{eq:lpolar} is a bijection and the 
invertibility of its differential implies that 
it is a diffeomorphism. 

\nin (4) Clearly $ss^\sharp \in S^\sharp_e$ for every $s \in S_e$. 
Hence it suffices to show that every $g \in (S^\sharp)_e$ 
is of the form $ss^\sharp$ for some $s \in S_e$. Since both sides are 
$L_e$-invariant, we may assume that $g = \exp(2x)$ for some $x \in C^{-\tau}$. 
Then the assertion follows with $s := \exp x$. 

\nin (5) Let $S^\sharp_c \subeq S^\sharp$ be a connected component 
and $\ell_0 \in p(S^\sharp_c)$. 
As $\ell_0 \in L \subeq G^\sigma$, the automorphism 
$\gamma$ of $G$ is an involution which commutes with $\sigma$. 
Further $-\gamma(C) = C$ follows from  the fact 
that $C$ is $\Ad(\ell_0)$-invariant. Therefore $\gamma$ satisfies 
both assumptions made for $\tau$. Next we observe that 
$\ell_0 = \ell_0^\sharp$ leads to 
\[
(\ell_0 s)^\sharp = s^\sharp \ell_0^\sharp 
= \ell_0 (\ell_0^{-1} s^\sharp \ell_0) 
= \ell_0 \gamma(s)^{-1}\] 
implies that 
\[ (\ell_0 S)^\sharp = \ell_0 \{ s \in S \: \gamma(s)^{-1} = s\}.\]  
The remaining assertions now follow from (3), applied to the involution 
$\gamma$ instead of~$\tau$. 
\end{prf}

\begin{rem} \mlabel{rem:a.3}
 The preceding lemma extends to connected  real Olshanski semigroups 
\[ \Gamma_L(C)= L \exp(C) \] 
not contained in a group $G$ (cf.~\cite[Def.~XI.1.11]{Ne00}). 
Here $(\g,\sigma)$ is a symmetric Lie algebra, $L$ is a connected Lie group 
with Lie algebra $\fl = \g^{\sigma}$ 
to which the adjoint action of $\fl$ on $\g^{-\sigma}$ integrates, 
and $C \subeq \g^{-\sigma}$ is a pointed generating $\Ad(L)$-invariant 
(weakly) hyperbolic cone. 

If $G$ is the simply connected Lie group with Lie algebra $\g$, 
these semigroups are constructed by first obtaining 
$\Gamma_{\tilde L}(C)$, where $q_L \: \tilde L \to L$ is the universal covering,
 as a the simply connected covering of $\Gamma_{G^\sigma}(C) 
= G^\sigma \exp(C) \subeq G$, and then factoring the discrete central 
subgroup $\ker(q_L) \cong \pi_1(L) \subeq Z(\tilde L)$ 
(see Definition~\ref{def:cplxols}  and \cite[\S XI.1]{Ne00} for more details).   
\end{rem}

\section{The topology on $\cH^\infty$ and the 
space $\cH^{-\infty}$} 
\mlabel{app:smovec}

Let $(U,\cH)$ be a unitary representation of $G$.
A {\it smooth vector} is an element $\eta\in\cH$ for which the orbit map 
$U^\eta : G\to \cH, g \mapsto U(g)\eta$ 
is smooth. We write~$\cH^{\infty} = \cH^\infty(U)$ for the space 
of smooth vectors. It carries the {\it derived representation} 
$\dd U $ of the Lie algebra $\fg$ given by
\begin{equation}
  \label{eq:derrep}
\dd U(x)\eta =\lim_{t\to 0}\frac{U(\exp t x)\eta -\eta}{t}.
\end{equation}
We extend this representation to a homomorphism 
$\dd U \:  \cU(\g) \to \End(\cH^\infty),$ 
where $\cU(\g)$ is the complex enveloping algebra of $\g$. This algebra 
carries an involution $D \mapsto D^*$ determined uniquely by 
$x^* = -x$ for $x \in \g$.
For $D \in \cU(\g)$, we obtain a seminorm on $\cH^\infty$ by 
\[p_D(\eta )=\|\dd U(D)\eta\|\quad \mbox{ for } \quad \eta \in \cH^\infty.\] 
These seminorms define a topology on $\cH^\infty$ which turn the injection 
\begin{equation}
  \label{eq:tophinfty}
 \eta \: {\cal H}^\infty \to {\cal H}^{{\cal U}(\g_\C)}, \quad 
\xi \mapsto (\dd U(D)\xi)_{D \in {\cal U}(\g_\C)}
\end{equation}
into a topological embedding, where the right hand side carries the product 
topology (cf.\ \cite[3.19]{Mag92}).  
Then $\cH^\infty$ becomes a complete locally convex space 
for which the linear operators $\dd U(D)$, $D \in \cU(\g)$, are continuous. 
Since $\cU(\g)$ has a countable basis, 
countably many such seminorms already determine the topology, so that 
$\cH^\infty$ is metrizable, hence a Fr\'echet space.
We also observe that 
the inclusion $\cH^\infty\hookrightarrow \cH$ is continuous.

The space $\cH^\infty$ of smooth vectors is $G$-invariant 
and we denote the corresponding representation by~$U^\infty$. We have 
the intertwining relation 
\[ \dd U(\Ad (g)x)= U^\infty(g) \dd U(x) U^\infty(g)^{-1} \quad \mbox{ for } \quad 
g \in G, x \in \g.\] 

If $\varphi \in C_c^\infty (G)$ and $\xi\in \cH$ then
$U(\varphi )\xi = \int_G U(g)\varphi(g)\, dg\in \cH^\infty$.
A sequence $(\varphi_n)_{n \in \N}$ in $C^\infty_c(G)$ is called a 
{\it $\delta$-sequence} if $\varphi_n \ge 0 $ 
and $\int_G \varphi_n(g)\, dg = 1$  for all $n\in \N$ 
and for every $e$-neighborhood $V \subeq G$, we have 
$\supp(\varphi_n) \subeq V$ if $n$ is sufficiently large. 
If $(\varphi_n)_{n \in \N}$ is a 
{$\delta$-sequence}, then $U(\varphi_n)\xi \to \xi$, so that $\cH^\infty$
is dense in~$\cH$.

We write $\cH^{-\infty}$ for the space 
of continuous anti-linear functionals on $\cH^\infty$. 
Its elements are called \textit{distribution vectors}. 
The group $G$, $\cU(\g)$ and $C^\infty_c(G)$ act on $\eta \in \cH^{-\infty}$ by
\begin{itemize}
\item $(U^{-\infty}(g)\eta ) (\xi ):= \eta  (U(g^{-1})\xi)$, $g \in G, 
\xi \in \cH^\infty$. \\ If $U \: G \to \AU(\cH)$ is an antiunitary 
representation and $U(g)$ is antiunitary, then we have to modify 
this definition slightly by 
$(U^{-\infty}(g)\eta) (\xi ):= \oline{\eta(U(g^{-1})\xi)}$. 
\item $(\dd U^{-\infty}(D) \eta ) (\xi ):= \eta(\dd U(D^*) \xi)$, $D \in \cU(\g), 
\xi \in \cH^\infty$. 
\item $U^{-\infty}(\varphi) \eta =\eta \circ U^\infty(\varphi^*)$, 
$\varphi \in C_c^\infty (G).$ 
\end{itemize} 
We have natural $G$-equivariant linear embeddings 
\begin{equation}
  \label{eq:embindistr}
\cH^\infty \into \cH
\mapright{\xi \mapsto  | \xi \ra} \cH^{-\infty}. 
\end{equation}

For each $\phi \in C^\infty_c(G)$, the map 
$U(\phi) \: \cH \to \cH^{\infty}$ 
is continuous, so that its adjoint defines a weak-$*$-continuous map 
$U^{-\infty}(\phi^*) \: \cH^{-\infty} \to \cH$. We actually have 
$U^{-\infty}(\phi)\cH^{-\infty} \subeq \cH^\infty$ 
as a consequence of the 
Dixmier--Malliavin Theorem \cite[Thm.~3.1]{DM78}, 
which asserts that every $\phi \in C^\infty_c(G)$ is a sum of functions 
of the form $\phi_1 * \phi_2$ with $\phi_j \in C^\infty_c(G)$.

\bigskip

 \noindent
 \footnotesize{Address of K-H. Neeb: Department of Mathematics,
 Friedrich-Alexander-University Erlangen--Nuremberg,
Cauerstrasse 11, 91058 Erlangen, Germany, 
 neeb@math.fau.de\\
 Address of G. \'Olafsson: 
 Department of mathematics, 
Louisiana State University, 
Baton Rouge, LA 70803,
olafsson@math.lsu.edu}
 
\end{document}